\documentclass[11pt,a4paper,english]{article}
\usepackage[T1]{fontenc}
\usepackage[latin1]{inputenc}
\usepackage{babel}
\usepackage{amsthm}
\usepackage{amsmath}
\usepackage{amssymb}
\usepackage{amsfonts}
\usepackage[all]{xy}
\usepackage{hyperref}


\DeclareMathOperator{\Alb}{Alb}

\DeclareMathOperator{\Bilabk}{Bilin_{\Abk}}
\DeclareMathOperator{\CH}{CH}

\DeclareMathOperator{\Divf}{\underline{Div}}

\DeclareMathOperator{\Ext}{Ext^{1}}
\DeclareMathOperator{\Extf}{\underline{Ext}^{1}}
\DeclareMathOperator{\Extabk}{Ext_{\Abk}^{1}}
\DeclareMathOperator{\Extfabk}{\underline{Ext}_{\Abk}^{1}}

\DeclareMathOperator{\Extfcttabk}{\underline{Ext}_{\Compl^{[-1,0]}(\Abk)}^{1}} 
 
\DeclareMathOperator{\Expah}{Exp}
\DeclareMathOperator{\Eah}{E}

\DeclareMathOperator{\Frob}{F}
\DeclareMathOperator{\Gal}{Gal}
\DeclareMathOperator{\Gl}{GL}

\DeclareMathOperator{\Hom}{Hom}
\DeclareMathOperator{\Homf}{\underline{Hom}}

\DeclareMathOperator{\Homka}{Hom_{\fld \mathrm{-alg}}}

\DeclareMathOperator{\Homcka}
                                        {Hom_{\fld \mathrm{-alg}}^{\mathrm{cont}}}

\DeclareMathOperator{\Homabk}{Hom_{\Abk}}

\DeclareMathOperator{\Homfabk}{\underline{Hom}_{\Abk}}
 
\DeclareMathOperator{\Homfcttabk}{\underline{Hom}_{\,\Compl^{[-1,0]}(\Abk)}}

\DeclareMathOperator{\K}{K}

\DeclareMathOperator{\Lie}{Lie}

\DeclareMathOperator{\Mat}{Mat}

\DeclareMathOperator{\Nil}{Nil}

\DeclareMathOperator{\R}{R}

\DeclareMathOperator{\Pic}{Pic}

\DeclareMathOperator{\Picf}{\underline{Pic}}

\DeclareMathOperator{\Spec}{Spec}

\DeclareMathOperator{\Spf}{Spf}

\DeclareMathOperator{\Supp}{Supp}

\DeclareMathOperator{\Ver}{V}

\DeclareMathOperator{\Z}{Z}

\DeclareMathOperator{\alb}{alb}

\DeclareMathOperator{\cl}{cl}
\DeclareMathOperator{\chr}{char}

\DeclareMathOperator{\dv}{div}

\DeclareMathOperator{\hgt}{ht}

\DeclareMathOperator{\id}{id}
\DeclareMathOperator{\im}{im}

\DeclareMathOperator{\nty}{n}
\DeclareMathOperator{\modu}{mod}
\DeclareMathOperator{\modun}{mod^{\Finv}}
\DeclareMathOperator{\moeb}{\mu}
\DeclareMathOperator{\multy}{\mu}

\DeclareMathOperator{\val}{v}

\newcommand{\Alba}[1]{\Alb(#1)}
\newcommand{\Albb}[1]{\Alb^{(#1)} \hspace{-1.5pt}{}}
\newcommand{\Fm}[2]{\fmlG_{#1,#2}}
\newcommand{\Fmo}[2]{\fmlG_{#1,#2}^{0}}
\newcommand{\Fmor}[2]{\fmlG_{#1,#2}^{0,\red}}

\newcommand{\Fmn}[2]{\fmlG^{\Finv}_{#1,#2}}

\newcommand{\Lm}[2]{L\left(#1,#2\right)}
\newcommand{\Jacc}[1]{\Jac\left(#1\right)}
\newcommand{\Jacm}[2]{\Jac\left(#1,#2\right)}
\newcommand{\Jaccm}[3]{\Jac^{(#1)}(#2,#3)}
\newcommand{\Albm}[2]{\Alb(#1,#2)}

\newcommand{\Albmm}[4]{\Alb_{#1,#2}^{#3,#4}}
\newcommand{\Albbm}[3]{\Alb^{(#1)}(#2,#3)}
\newcommand{\albbm}[3]{\alb^{(#1)}_{#2,#3}}

\newcommand{\albbF}[1]{\alb_{\fmlG}^{(#1)}}

\newcommand{\albFb}{\alb_{\fmlGb}}

\newcommand{\AlbF}[1]{\Alb_{\fmlG}(#1)}
\newcommand{\AlbFb}[1]{\Alb_{\fmlGb}(#1)}
\newcommand{\AlbbF}[2]{\Alb_{\fmlG}^{(#1)}(#2)}

\newcommand{\Picor}[1]{\Pic^{0,\red}_{#1}}
\newcommand{\Picorf}[1]{\Picf^{0,\red}_{#1}}

\newcommand{\Divfc}[1]{\wh{\Divf_{#1}}}
\newcommand{\Divorc}[1]{\wh{\Divf^{0,\red}_{#1}}}
\newcommand{\Divor}[1]{\Divf^{0,\red}_{#1}}

\newcommand{\Extt}[1]{\mathrm{Ext}^{#1}}

\newcommand{\Extabkt}[1]{\mathrm{Ext}_{\Abk}^{#1}}
\newcommand{\Extfabkt}[1]{\underline{\mathrm{Ext}}_{\Abk}^{#1}}
 
\newcommand{\Extabs}[1]{\mathrm{Ext}_{\mathcal{A} \mathit{b} / #1}^{1}}

\newcommand{\Homabs}[1]{\mathrm{Hom}_{\mathcal{A} \mathit{b} / #1}}

\newcommand{\Homfabs}[1]
                    {\underline{\mathrm{Hom}}_{\mathcal{A}\mathit{b}/#1}}

\renewcommand{\H}{\mathrm{H}}

\newcommand{\Wcfl}[1]{{}_{#1} \hspace{-1pt} \Wittc}

\newcommand{\WeilRes}[2]{\Pi_{#1/#2}}

\newcommand{\WR}[2]{#1 \left( \llul \tens #2 \right)}

\newcommand{\Rlm}[2]{\fR_0 (#1,#2)}

\newcommand{\ZOo}[1]{\Z_0 (#1)^{0}}

\newcommand{\CHm}[2]{\CH_0 (#1,#2)}
\newcommand{\CHmo}[2]{\CHm{#1}{#2}^{0}}

\newcommand{\Ab}{\mathsf{Ab}}
\newcommand{\Abk}{\mathcal{A} \mathit{b} / \fld}
\newcommand{\AbS}[1]{\mathcal{A} \mathit{b} / #1}

\newcommand{\Alg}[1]{\mathsf{Alg} / #1}

\newcommand{\Algk}{\mathsf{Alg} / \fld}

\newcommand{\Artk}{\mathsf{Art} / \fld}
\newcommand{\Compl}{\mathcal{C}}
\newcommand{\Cat}{\mathsf{C}}

\newcommand{\Gk}{\mathcal{G} / \fld}
\newcommand{\Gak}{\mathcal{G} \hspace{-.1em} \mathit{a} / \fld}
\newcommand{\aGk}{\mathit{a} \hspace{+.06em} \mathcal{G} / \fld}
\newcommand{\aGak}{\mathit{a} \hspace{+.06em} \mathcal{G} 
                              \hspace{-.1em} \mathit{a} / \fld}
\newcommand{\Gfk}{\mathcal{G} \hspace{-.01em} \mathit{f} / \fld}

\newcommand{\dGfk}{\mathit{d} \hspace{+.06em} \mathcal{G} 
                              \hspace{-.01em} \mathit{f} / \fld}

\newcommand{\Fctr}{\mathsf{Fctr}}

\newcommand{\Set}{\mathsf{Set}}

\newcommand{\Mav}{\mathsf{Mav}}
\newcommand{\Mr}{\mathsf{Mr}}
\newcommand{\Mrt}[2]{\mathsf{Mr}_{#1} \pn{#2}}
\newcommand{\MrF}{\mathsf{Mr}_{\fmlG} \pn{X}}
\newcommand{\MrFb}{\mathsf{Mr}_{\fmlGb}}

\newcommand{\Mrm}[2]{\mathsf{Mr} \pn{#1,#2}}

\newcommand{\MrCHm}[2]{\mathsf{Mr}^{\CH} \pn{#1,#2}}

\newcommand{\Nat}{\mathbb{N}}
\newcommand{\Zint}{\mathbb{Z}}

\newcommand{\Afn}{\mathbb{A}}

\newcommand{\et}{\acute{\mathrm{e}} \mathrm{t}}

\newcommand{\fppf}{\mathrm{fppf}}

\newcommand{\mdl}{D}
\newcommand{\mdll}{E}

\newcommand{\num}{r}

\newcommand{\red}{\mathrm{red}}

\newcommand{\der}{\mathrm{d}}
\newcommand{\pder}{\mathrm{\del}}
\newcommand{\gal}{\sigma}

\newcommand{\hgg}[2]{h^{(#1)}_{#2}}

\newcommand{\homm}[1]{h^{(#1)}}

\newcommand{\lin}{\ell}

\newcommand{\morphism}{\psi}

\newcommand{\pr}{\mathrm{pr}}

\newcommand{\repr}{\rho}

\newcommand{\trafo}{\tau}

\newcommand{\Pop}{\Phi}

\newcommand{\clfld}{\overline{\fld}}
\newcommand{\fld}{\mathit{k}}

\newcommand{\Kab}{\K^{\mathrm{ab}}}

\newcommand{\bring}{R}

\newcommand{\Algbrc}{\sA}

\newcommand{\Ring}{R}
\newcommand{\Rfin}{R}

\newcommand{\Sfin}{S}

\newcommand{\Pb}{\mathrm{P}}
\newcommand{\Gp}{\mathrm{G}}

\newcommand{\fmlE}{\mathcal{E}}

\newcommand{\fmlG}{\mathcal{F}}
\newcommand{\fmlGb}{\overline{\mathcal{F}}}
\newcommand{\fmlGd}{\fmlG^{\vee}}

\newcommand{\fmlGr}{\mathcal{G}}

\newcommand{\Ga}{\mathbb{G}_{\mathrm{a}}}
\newcommand{\Gac}{\widehat{\mathbb{G}}_{\mathrm{a}}}

\newcommand{\Gm}{\mathbb{G}_{\mathrm{m}}}
\newcommand{\GmS}[1]{\mathbb{G}_{\mathrm{m},#1}}
\newcommand{\Jac}{J}
\newcommand{\Kot}{\Omega}
\newcommand{\Lin}{\mathbb{L}}

\newcommand{\Trr}[1]{G^{(#1)}}
\newcommand{\Trs}{T}
\newcommand{\Trsd}{\Trs^{\vee}}

\newcommand{\Upt}{U}

\newcommand{\Upf}{\mathbb{U}}

\newcommand{\Witt}{W}
\newcommand{\Wittc}{\wh{\Witt}}
\newcommand{\Wpl}{\Lambda}

\newcommand{\fil}{\mathrm{fil}}
\newcommand{\filf}{\mathrm{fil}^{\Frob}}
\newcommand{\filb}{{}^{\flat} \mathrm{fil}}
\newcommand{\filfb}{{}^{\flat} \mathrm{fil}^{\Frob}}
\newcommand{\qfil}{\fd}
\newcommand{\Qfil}{\fD}
\newcommand{\qfill}{\ol{\fd}}
\newcommand{\Qfill}{\ol{\fD}}
\newcommand{\Qfillb}{{}^{\flat} \Qfill}

\newcommand{\pnt}{q}
\newcommand{\pntt}{p}

\newcommand{\Es}{E}

\newcommand{\Crv}{C}

\newcommand{\Xo}{X}
\newcommand{\Xb}{\overline{X}}

\newcommand{\Ct}{\widetilde{C}}

\newcommand{\Lamd}{\Lam^{\vee}}

\newcommand{\Ed}{E^{\vee}}

\newcommand{\Nd}{N^{\vee}}
\newcommand{\Md}{M^{\vee}}

\newcommand{\Ad}{A^{\vee}}

\newcommand{\Gd}{G^{\vee}}
\newcommand{\Ld}{L^{\vee}}
\newcommand{\ld}{l^{\vee}}
\newcommand{\hd}{h^{\vee}}

\newcommand{\Mdd}{M^{\vee \vee}}

\newcommand{\alp}{\alpha}

\newcommand{\gam}{\gamma}
\newcommand{\del}{\delta}

\newcommand{\tha}{\vartheta}

\newcommand{\lam}{\lambda}
\newcommand{\sig}{\sigma}

\newcommand{\phe}{\varphi}
\newcommand{\oma}{\omega}

\newcommand{\Gam}{\Gamma}

\newcommand{\Lam}{\Lambda}

\newcommand{\fd}{\mathfrak{d}}

\newcommand{\fm}{\mathfrak{m}}

\newcommand{\fD}{\mathfrak{D}}

\newcommand{\fR}{\mathfrak{R}}

\newcommand{\sA}{\mathcal{A}}

\newcommand{\sD}{\mathcal{D}}

\newcommand{\sF}{\mathcal{F}}

\newcommand{\sK}{\mathcal{K}}
\newcommand{\sL}{\mathcal{L}}

\newcommand{\sO}{\mathcal{O}}

\newcommand{\sQ}{\mathcal{Q}}

\newcommand{\sU}{\mathcal{U}}

\newcommand{\bF}{\mathbb{F}}

\newcommand{\bt}[1]{[#1]}
\newcommand{\lrbt}[1]{\left[#1\right]}
\newcommand{\bigbt}[1]{\big[#1\big]}

\newcommand{\pn}[1]{(#1)}
\newcommand{\lrpn}[1]{\left(#1\right)}
\newcommand{\bigpn}[1]{\big(#1\big)}
\newcommand{\Bigpn}[1]{\Big(#1\Big)}
\newcommand{\biggpn}[1]{\bigg(#1\bigg)}

\newcommand{\lrpair}[1]{\left\langle #1\right\rangle}

\newcommand{\lrst}[1]{\left\{#1\right\}}

\newcommand{\ra}{\rightarrow}
\newcommand{\la}{\leftarrow}
\newcommand{\dra}{\dashrightarrow}

\newcommand{\lra}{\longrightarrow}
\newcommand{\lla}{\longleftarrow}

\newcommand{\Lra}{\Longrightarrow}
\newcommand{\Lla}{\Longleftarrow}
\newcommand{\Llra}{\Longleftrightarrow}
\newcommand{\lmt}{\longmapsto}

\newcommand{\iso}{\overset{\sim}\longrightarrow}
\newcommand{\inj}{\hookrightarrow}
\newcommand{\linj}{\lra}
\newcommand{\sur}{\twoheadrightarrow}
\newcommand{\lsur}{-\hspace{-.7em} -\hspace{-.7em} \twoheadrightarrow}

\newcommand{\wt}[1]{\widetilde{#1}}
\newcommand{\wh}[1]{\widehat{#1}}
\newcommand{\ol}[1]{\overline{#1}}
\newcommand{\ul}[1]{\underline{#1}}
\newcommand{\cut}{\cdot}

\newcommand{\dsum}{\bigoplus}

\newcommand{\tens}{\otimes}

\newcommand{\tensc}{\hspace{+.15em} \widehat{\otimes} \hspace{+.15em}}

\newcommand{\tms}{\times}

\newcommand{\blank}{\_}
\newcommand{\lul}{?}                    %%%  {\_}
\newcommand{\llul}{\hspace{+.05em} ?}   %%%  {\_}
\newcommand{\lull}{? \hspace{+.05em}}   %%%  {\_}
\newcommand{\see}{}
\newcommand{\seecite}{\textrm{see }}
\newcommand{\laurin}{. }                     %%% {.} 
\newcommand{\laurink}{, }                   %%% {,} 
\newcommand{\laurind}{: }                   %%% {:} 
\newcommand{\bDpl}{\;$}
\newcommand{\eDpl}{$\;}

\newcommand{\vs}{12pt}

\newcommand{\kukwad}{0pt}

\newcommand{\bThm}{\begin{theorem}}
\newcommand{\eThm}{\end{theorem}}
\newcommand{\bAck}{\begin{acknowledgement}}
\newcommand{\eAck}{\end{acknowledgement}}
\newcommand{\bAlg}{\begin{algorithm}}
\newcommand{\eAlg}{\end{algorithm}}
\newcommand{\bAxm}{\begin{axiom}}
\newcommand{\eAxm}{\end{axiom}}
\newcommand{\bCas}{\begin{case}}
\newcommand{\eCas}{\end{case}}
\newcommand{\bClm}{\begin{claim}}
\newcommand{\eClm}{\end{claim}}
\newcommand{\bCcl}{\begin{conclusion}}
\newcommand{\eCcl}{\end{conclusion}}
\newcommand{\bCdn}{\begin{condition}}
\newcommand{\eCdn}{\end{condition}}
\newcommand{\bCjc}{\begin{conjecture}}
\newcommand{\eCjc}{\end{conjecture}}
\newcommand{\bCor}{\begin{corollary}}
\newcommand{\eCor}{\end{corollary}}
\newcommand{\bCrt}{\begin{criterion}}
\newcommand{\eCrt}{\end{criterion}}
\newcommand{\bDef}{\begin{definition}}
\newcommand{\eDef}{\end{definition}}
\newcommand{\bExm}{\begin{example}}
\newcommand{\eExm}{\end{example}}
\newcommand{\bExc}{\begin{exercise}}
\newcommand{\eExc}{\end{exercise}}
\newcommand{\bLem}{\begin{lemma}}
\newcommand{\eLem}{\end{lemma}}
\newcommand{\bNot}{\begin{notation}}
\newcommand{\eNot}{\end{notation}}
\newcommand{\bPar}{\begin{para}}
\newcommand{\ePar}{\end{para}}
\newcommand{\bPnt}{\begin{point}}
\newcommand{\ePnt}{\end{point}}
\newcommand{\bPrb}{\begin{problem}}
\newcommand{\ePrb}{\end{problem}}
\newcommand{\bPrp}{\begin{proposition}}
\newcommand{\ePrp}{\end{proposition}}
\newcommand{\bRmk}{\begin{remark}}
\newcommand{\eRmk}{\end{remark}}
\newcommand{\bSol}{\begin{solution}}
\newcommand{\eSol}{\end{solution}}
\newcommand{\bSmr}{\begin{summary}}
\newcommand{\eSmr}{\end{summary}}
\newcommand{\bVar}{\begin{variant}}
\newcommand{\eVar}{\end{variant}}
\newcommand{\bPf }{\begin{prooof}}
\newcommand{\ePf }{\end{prooof}}

\theoremstyle{plain}
\newtheorem{theorem}{Theorem}[section]
\newtheorem{axiom}[theorem]{Axiom}
\newtheorem{conjecture}[theorem]{Conjecture}
\newtheorem{corollary}[theorem]{Corollary}
\newtheorem{criterion}[theorem]{Criterion}
\newtheorem{lemma}[theorem]{Lemma}
\newtheorem{problem}[theorem]{Problem}
\newtheorem{proposition}[theorem]{Proposition}

\theoremstyle{definition}
\newtheorem{acknowledgement}[theorem]{Acknowledgement}
\newtheorem{algorithm}[theorem]{Algorithm}
\newtheorem{case}[theorem]{Case}
\newtheorem{claim}[theorem]{Claim}
\newtheorem{condition}[theorem]{Condition}
\newtheorem{conclusion}[theorem]{Conclusion}
\newtheorem{definition}[theorem]{Definition}
\newtheorem{example}[theorem]{Example}
\newtheorem{exercise}[theorem]{Exercise}
\newtheorem{notation}[theorem]{Notation}
\newtheorem{point}[theorem]{Point}
\newtheorem{remark}[theorem]{Remark}
\newtheorem{solution}[theorem]{Solution}
\newtheorem{summary}[theorem]{Summary}

\newtheorem{para}[subsubsection]{}

\newenvironment{prooof}[1][Proof]{\textbf{#1.} }
{\ \rule{0.5em}{0.5em}} 
\newenvironment{variant}[1][Variant]{\textbf{#1.} }
{\ \rule{0.5em}{0.5em}}

\begin{document}

%\title{Albanese with Modulus\\
%         over a Perfect Field} 
%\author{Henrik Russell}
%\maketitle

\centerline{ } 
\vspace{10pt} 
\centerline{\LARGE{Albanese varieties with modulus}} 
\vspace{5pt} 
\centerline{\LARGE{over a perfect field}} 
\vspace{15pt}
\centerline{\large{Henrik Russell}} 
\vspace{15pt} 

%%%%%%%%  Abstract  %%%%%%%%%%%%%%%%%%%%%%%
\begin{abstract} 
Let $X$ be a smooth proper variety over a perfect field $k$ 
of arbitrary characteristic. 
Let $D$ be an effective divisor on $X$ with multiplicity. 
We introduce an Albanese variety $\Albm{X}{D}$ of $X$ of modulus $D$ 
as a higher dimensional analogue of 
the generalized Jacobian of Rosenlicht-Serre 
with modulus for smooth proper curves. 
%The Albanese variety with modulus appears as a special case of 
%a broader notion of universal objects of 
%categories of rational maps from $X$ to commutative algebraic groups. 
%defined by a universal mapping property for 
Basing on duality of 1-motives with unipotent part 
(which are introduced here), 
we obtain explicit and functorial descriptions 
of these generalized Albanese varieties and their dual functors. 

We define a relative Chow group of zero cycles $\CHm{X}{\mdl}$ 
of modulus $\mdl$ and show that $\Albm{X}{\mdl}$ can be viewed 
as a universal quotient of $\CHmo{X}{\mdl}$. 

As an application %of Albanese varieties with modulus 
we can rephrase Lang's class field theory 
of function fields of varieties over finite fields 
in explicit terms. 
\end{abstract}

%\newpage
%\quad 
%\newpage

\setcounter{tocdepth}{2}
\tableofcontents{} 
%\newpage 

\setcounter{section}{-1}

\newpage

%%%%%%%%%%%%  Introduction  %%%%%%%%%%%%%%%%%%
\section{Introduction} 

The generalized Jacobian variety with modulus 
of a smooth proper curve $X$ over a field 
is a well-established object in algebraic geometry and number theory 
and has shown to be of great benefit 
e.g.\ for the theory of algebraic groups, ramification theory 
and class field theory.
In this work we extend this notion from \cite[V]{S} 
to the situation of a higher dimensional smooth proper variety $X$ 
over a perfect field $k$. 
The basic idea of this construction comes from \cite{Ru1} 
and is accomplished in \cite{KR1}, 
both only for the case that $k$ is of characteristic 0. 
Positive characteristic however requires distinct methods 
and turns out to be the difficult part of the story. 

To a rational map $\phe: X \dra P$ 
from $X$ to a torsor $P$ under a commutative algebraic group $G$ 
we assign an effective divisor $\modu\pn{\phe}$, the \emph{modulus of $\phe$} 
(Def.~\ref{DefMod}). 
Our definition from \cite{KR2} 
coincides with the classical definition in the curve case as in \cite[III, No.~1]{S}. 
For an effective divisor $D$ on $X$ 
the generalized Albanese variety $\Albbm{1}{X}{\mdl}$ of $X$ of modulus $\mdl$ 
and the Albanese map $\albbm{1}{X}{\mdl}:X\dra\Albbm{1}{X}{\mdl}$ 
are defined by the following universal property: 
for every torsor $P$ under a commutative algebraic group $G$ 
and every rational map $\phe$ from $X$ to $P$ of modulus $\leq\mdl$ 
there exists a unique homomorphism of torsors 
$h:\Albbm{1}{X}{\mdl}\lra P$ such that $\phe=h\circ\albbm{1}{X}{\mdl}$. 
%up to translation by a constant $g \in G(k)$. 
Every rational map to a torsor for a commutative algebraic group 
admits a modulus, 
and the effective divisors on $X$ form an inductive system. 
Then the projective limit $\varprojlim \Albbm{1}{X}{\mdl}$ 
over all effective divisors $\mdl$  on $X$ 
yields a torsor for a pro-algebraic group 
that satisfies the universal mapping property 
for all rational maps from $X$ 
to torsors for commutative algebraic groups. 

The Albanese variety with modulus 
(Thm.~\ref{AlbMod}) %and Thm.~\ref{descentAlbm(X,D)}) 
arises as a special case of a broader notion 
of generalized Albanese varieties defined by a universal mapping property 
for categories of rational maps from $X$ 
to torsors for commutative algebraic groups. 
As the construction of these universal objects is based on duality, 
a notion of duality for smooth connected commutative algebraic groups 
over a perfect field $k$ of arbitrary characteristic is required. 
For this purpose we introduce so called \emph{1-motives with unipotent part} 
(Def.~\ref{Def 1-motive}), 
which generalize Deligne 1-motives \cite[D\'efinition~(10.1.2)]{D2} 
and Laumon 1-motives \cite[D\'efinition~(5.1.1)]{L}. 
In this context, we obtain explicit and functorial descriptions 
of these generalized Albanese varieties and their dual functors 
(Thm.~\ref{univ_object}). %and Thm.~\ref{descent}). 

In a geometric way we define a relative Chow group of 0-cycles 
$\CHm{X}{\mdl}$ with respect to the modulus $\mdl$ 
(Def.~\ref{DefCHoMod}). 
Then we can realize $\Albbm{1}{X}{\mdl}$ as a universal quotient of 
$\CHmo{X}{\mdl}$, the subgroup of $\CHm{X}{\mdl}$ of cycles of degree 0
(Thm.~\ref{CHoMod}), 
in the case that the base field is algebraically closed. 
The relation of $\CHm{X}{\mdl}$ 
to the K-theoretic id\`ele class groups from \cite{KS} 
gives rise to some future study, 
but is beyond the scope of this paper. 
Using these id\`ele class groups, 
\"Onsiper \cite{On} proved the existence 
of generalized Albanese varieties for smooth proper surfaces 
in characteristic $p > 0$. 

Lang's class field theory of function fields of varieties over finite fields 
\cite[V]{S} is written in terms of so called \emph{maximal maps}, 
which appeared as a purely theoretical notion, 
apart from their existence very little seemed to be known about which. 
The Albanese map with modulus allows 
to replace these black boxes by concrete objects 
(Thm.~\ref{geomCFT}).  

We present the main results %of this paper 
by giving a summary of each section. 

%\newpage 

\subsection{Leitfaden} 

%\medskip 
\hspace{11pt} 
\textbf{Section \ref{1-Motives}} 
is devoted to the following generalization of 1-motives: 
%from \cite{Ru2} 
A \emph{1-motive with unipotent part} 
($\see$Definition \ref{Def 1-motive}) 
is roughly a homomorphism $\lrbt{\fmlG \ra G}$ 
in the category of sheaves of abelian groups 
over %an algebraically closed field $k$ 
a perfect field $k$ 
from a dual-algebraic commutative formal group $\fmlG$ 
to an extension $G$ of an abelian variety $A$ 
by a commutative affine algebraic group $L$. 
Here a commutative formal group $\fmlG$ is called \emph{dual-algebraic} 
if its Cartier-dual $\fmlGd = \Homf\pn{\fmlG,\Gm}$ is algebraic. 
1-motives with unipotent part admit duality ($\see$Def.~\ref{Def:dual1-mot}). 
The dual of $\left[0\ra G\right]$
is given by $\left[\Ld\ra\Ad\right]$, 
where $\Ld=\Homf\pn{L,\Gm}$ is the Cartier-dual of $L$ 
and $\Ad=\Pic_{A}^{0}=\Extf\pn{A,\Gm}$ is the dual abelian variety of $A$, 
and the homomorphism between them is the connecting homomorphism 
associated to $0 \ra L \ra G \ra A \ra 0$. 
In particular, every smooth connected commutative algebraic group over $k$ 
has a dual in this category. 
Moreover, these 1-motives may contain torsion. 

\medskip 

\textbf{Section \ref{sec:Univ-Fact-Prbl}}. 
Let $X$ be a smooth proper variety over a perfect field $k$. 
In the framework of 
\emph{categories of rational maps from $X$ 
to torsors for commutative algebraic groups} 
($\see$Def.~\ref{CatMr}), 
we ask for the existence of universal objects 
($\see$Def.~\ref{univObj}) for such categories, 
i.e.\ objects having the universal mapping property 
with respect to the category they belong to. 
Assume for the moment $k$ is an algebraically closed field. 
Then a torsor can be identified with the algebraic group acting on it. 
A necessary and sufficient condition for the existence of such 
universal objects is given in Theorem \ref{Exist univObj}, 
as well as their explicit construction, %($\see$Remark \ref{Alb_constr}), 
using duality of 1-motives with unipotent part. 
(This was done in \cite{Ru1} for $\chr(k) = 0$.) 
We pass to general perfect base field in Theorem \ref{univ_object}. 

In particular we show the following: 
Let $\Divf_X$ be the sheaf of relative Cartier divisors, 
i.e.\ the sheaf of abelian groups that assigns to any $k$-algebra $R$ 
the group $\Divf_X(R)$ of all Cartier divisors on $X \tens_k R$ 
generated locally on $\Spec R$ by effective divisors which are flat over $R$. 
Let $\Picf_X$ be the Picard functor and $\Picor{X}$ the Picard variety of $X$. 
Then let $\Divor{X}$ be the inverse image of $\Picor{X}$ 
under the class map $\cl: \Divf_X \lra \Picf_X$. 
A rational map $\phe:X\dra G$, 
where $G$ is a smooth connected commutative algebraic group 
with affine part $L$,
induces a natural transformation $\trafo_{\phe}:\Ld\lra\Divor{X}$ 
($\see$No.~\ref{sec:indTrafo}). 
If $\fmlG$ is a %dual-algebraic 
formal subgroup of $\Divor{X}$, 
denote by $\MrF$ the category of rational maps for which the
image of this induced transformation lies in $\fmlG$. 
If $\fld$ is an arbitrary perfect base field, %and $\fmlG$ is defined over $k$, 
we define $\MrF$ via base change to an algebraic closure $\clfld$ 
($\see$Def.~\ref{Mr_F}). 

\bThm 
\label{univ_object}
Let $\fmlG$ be a dual-algebraic formal $k$-subgroup of $\Divor{X}$. 
The category $\MrF$ admits a universal object $\albbF{1}:X\dra\AlbbF{1}{X}$. 
Here $\AlbbF{1}{X}$ is a torsor for an algebraic group $\AlbbF{0}{X}$, 
which arises as an extension of the classical Albanese $\Alba{X}$ 
by the Cartier-dual of $\fmlG$. 
The algebraic group $\AlbbF{0}{X}$ 
is dual to the 1-motive $\bigbt{\fmlG\lra\Picor{X}}$,
the homomorphism induced by the class map 
$\cl:\Divf_{X}\lra\Picf_{X}$. 
\eThm 

Theorem \ref{univ_object} results from (the stronger) Theorem \ref{Exist univObj}, 
which says that a category of rational maps to algebraic groups 
(over an algebraically closed field) 
admits a universal object 
if and only if it is of the shape $\MrF$ 
for some dual-algebraic formal subgroup $\fmlG$ of $\Divor{X}$, 
and Galois descent ($\see$Thm.~\ref{descent}). 
%($\see$No.~\ref{descentBaseField}). 
The generalized Albanese varieties $\AlbbF{i}{X}$ (i = 1,0) 
satisfy an obvious functoriality property 
($\see$Thm.~\ref{alb_F(morphism)}). 
%($\see$No.~\ref{Functoriality}). 

\medskip 

\textbf{Section \ref{sec:Modulus}} 
is the main part of this work, 
where we establish a higher dimensional analogue 
to the generalized Jacobian with modulus of Rosenlicht-Serre. 
Let $X$ be a smooth proper variety over a perfect field $k$. 
We use the notion of modulus from \cite{KR2}, 
which associates to a rational map $\phe: X \dra P$ 
%to a commutative algebraic group $G$ 
an effective divisor $\modu\pn{\phe}$ on $X$ ($\see$Def.~\ref{DefMod}). 
If $\mdl$ is an effective divisor on $X$, 
we define a  formal subgroup 
$\Fm{X}{\mdl} = \lrpn{\Fm{X}{\mdl}}_{\et} \tms_k \lrpn{\Fm{X}{\mdl}}_{\inf}$ 
of $\Divf_X$ (cf.\ Def.~\ref{DefFm(X,D)}) 
by the conditions 
\[ \lrpn{\Fm{X}{\mdl}}_{\et} = 
   \left\{ B\in\Divf_X(k) \;\big|\; \Supp(B) \subset \Supp(\mdl) \right\}
\] 
and for $\chr(k)=0$ 
\[ \lrpn{\Fm{X}{\mdl}}_{\inf} = \exp \lrpn{ \Gac \tens_k 
   \Gam\bigpn{X, \sO_{X}\lrpn{\mdl-\mdl_{\red}} \big/ \sO_{X}}} 
\] 
where \;$\exp$\; is the exponential map 
and $\Gac$ is the completion of $\Ga$ at $0$, 

\noindent 
and for $\chr(k)=p>0$ 
\[ \lrpn{\Fm{X}{\mdl}}_{\inf} = 
   \Expah \lrpn{ \sum_{r > 0} \Wcfl{r} \tens_{\Witt(k)} 
   \Gam\Bigpn{X, \filf_{\mdl-\mdl_{\red}} \Witt_r(\sK_{X}) \Big/ \Witt_r(\sO_{X})},1}  
\] 
where $\mdl_{\red}$ is the underlying reduced divisor of $\mdl$, 
$\Expah$ denotes the Artin-Hasse exponential, 
$\Wcfl{r}$ is the kernel of the $r^{\textrm{th}}$ power of the Frobenius 
on the completion $\Wittc$ of the Witt group $\Witt$ at $0$ 
and $\filf_{\mdl} \Witt_r(\sK_X)$ is a filtration of the Witt group 
($\see$Definition \ref{fil_D Witt}). 
Let $\Fmor{X}{\mdl} = \Fm{X}{\mdl} \tms_{\Divf_X} \Divor{X}$ 
be the intersection of $\Fm{X}{\mdl}$ and $\Divor{X}$. 
The formal groups $\Fm{X}{\mdl}$ and $\Fmor{X}{\mdl}$ are dual-algebraic 
($\see$Prop.\ \ref{Fm(X,D)algebraic}). 

Then it holds \;$\modu\pn{\phe} \leq \mdl$\; if and only if 
\;$\im\pn{\trafo_{\phe}} \subset \Fmor{X}{\mdl}$\; 
($\see$Lem.\ \ref{mod(phe)-im(trafo_phe)}). 
This yields (cf.\ Thm.\ \ref{AlbMod-construction} 
and Thm.\ \ref{descentAlbm(X,D)}) 
%No.~\ref{Alb(X,D)descent}) 

\bThm 
\label{AlbMod}
The category $\Mrm{X}{\mdl}$ 
of those rational maps $\phe:X\dra P$ s.t.\ \,$\modu\pn{\phe} \leq \mdl$ 
admits a universal object $\albbm{1}{X}{\mdl}:X\dra\Albbm{1}{X}{\mdl}$, 
called the \emph{Albanese of $X$ of modulus $\mdl$}. 
The algebraic group $\Albbm{0}{X}{\mdl}$ 
acting on $\Albbm{1}{X}{\mdl}$ 
is dual to the 1-motive $\bigbt{\Fmor{X}{\mdl}\lra \Picor{X}}$. 
\eThm 

The Albanese varieties with modulus $\Albbm{i}{X}{\mdl}$ (i = 1,0) 
are functorial ($\see$Prop.\ \ref{alb_X,D(morphism)}). 
%($\see$No.~\ref{Alb(X,D)functorial}). 
%and descend to arbitrary perfect base field 
%($\see$No.~\ref{Alb(X,D)descent}). 
In the case that $X=C$ is a curve, 
our Albanese with modulus $\Albbm{i}{C}{D}$ coincide with the 
generalized Jacobians with modulus $\Jaccm{i}{C}{D}$ of Rosenlicht-Serre 
($\see$Thm.\ \ref{dual_Jacm} and Galois descent). 

A relative Chow group $\CHm{X}{\mdl}$ of modulus $\mdl$ 
is introduced in Definition \ref{DefCHoMod}. 
We say a rational map $\phe:X \dra P$ to a torsor $P$ 
under a commutative algebraic group $G$ 
\emph{factors through $\CHmo{X}{\mdl}$} 
if the associated map on 0-cycles of degree 0 
(where $U$ is the open set on which $\phe$ is defined) \; 
%\begin{eqnarray*}
  $\ZOo{U}  \lra  G(k)$, 
  $\sum l_i \, p_i  \lmt  \sum l_i \, \phe(p_i)$  \;
%\end{eqnarray*} 
factors through a homomorphism of abstract groups 
$\CHmo{X}{\mdl}\lra G(k)$. 
We show (cf.\ Thm.\ \ref{ThmCHoMod}) 

\bThm 
\label{CHoMod}
Assume $k$ is algebraically closed. 
A rational map $\phe:X \dra P$ factors through $\CHmo{X}{\mdl}$ 
if and only if it factors through $\Albbm{1}{X}{\mdl}$. 
In other words, $\Albbm{0}{X}{\mdl}$ is a universal quotient of 
$\CHmo{X}{\mdl}$. 
\eThm 

The theory of Albanese varieties with modulus has an application 
to the class field theory of function fields of varieties over finite fields. 
Let $\Xo$ be a geometrically irreducible projective variety over 
a finite field $\fld = \bF_q$. Let $\clfld$ be an algebraic closure of $\fld$. 
Let $\K_{\Xo}$ denote the function field of ${\Xo}$, 
let $\Kab_{\Xo}$ be the maximal abelian extension of $\K_{\Xo}$. 
%In particular, $\Kab_{\Xo}$ contains $\clfld$. 
From Lang's class field theory one obtains 
\bThm 
\label{geomCFT}
The geometric Galois group $\Gal\bigpn{\Kab_{\Xo}/\K_{\Xo}\clfld}$ 
is isomorphic to the projective limit of the $\fld$-rational points 
of the Albanese varieties of $\Xo$ with modulus $\mdl$ 
\[ \Gal\bigpn{\Kab_{\Xo} \big/ \K_{\Xo}\clfld} \,\cong\, 
   \varprojlim_{\mdl} \, \Albbm{0}{\Xo}{\mdl}\pn{\fld} 
\] 
where $\mdl$ ranges over all effective divisors on $\Xo$ 
rational over $\fld$. 
\eThm 

The proof of Theorem \ref{geomCFT} 
is analogous to the proof of Lang's class field theory 
given in \cite[VI, \S~4, No.~16--19]{S}, 
replacing \emph{maximal maps} by the universal objects 
\;$\albbm{1}{\Xo}{\mdl}: \Xo \dra \Albbm{1}{\Xo}{\mdl}$ 
for the category of rational maps to $\fld$-torsors of modulus $\leq \mdl$ 
from Theorem \ref{AlbMod}. %\ref{descentAlbm(X,D)}. 

\vspace{\vs} 

\textbf{Acknowledgement.} 
I owe many thanks to Kazuya Kato for his hospitality, help and support. 
His influence on this work is considerable. 
I also thank the referee for helpful suggestions. In particular I replaced 
my original %(explicit) 
proof of Theorem \ref{Ext(A,L)=Hom(Ld,Ad)} 
by a shorter argument due to the referee.

%\newpage 

%%%%%%%%%%%  Section:  1-Motives  %%%%%%%%%%%%%%%%%
\section{1-motives} 
\label{1-Motives}

The aim of this section is to construct a category of generalized 1-motives 
that contains all smooth connected commutative algebraic groups over 
a perfect field 
and provides a notion of duality for them. 

\subsection{Algebraic groups and formal groups} 
\label{AlgGroups-FmlGroups}

I will use the language of group functors, 
algebraic groups and formal groups. 
References for algebraic groups are \cite{DG} and \cite{W}, 
for formal groups and Cartier duality are 
\cite[$\mathrm{VII}_\mathrm{B}$]{SGA3}, 
\cite[II]{D} and \cite[I]{Fo}. 
%A key-tool for our description of Cartier duality will be the functor 
%$\Ring \lmt \Lin_{\Ring}$ ($\see$No.~\ref{LinGroup_Ring}), 
%which assigns to a $\fld$-algebra $\Ring$ the Weil restriction 
%$\Lin_{\Ring} := \WeilRes{\Ring}{\fld}\GmS{\Ring}$ 
%%$\Gm\lrpn{\llul\tens_k\Ring}$ 
%of $\GmS{\Ring}$ from $\Ring$ to $\fld$. 

By \emph{algebraic group} resp.\ \emph{formal group} 
I will always mean a \emph{commutative} (algebraic resp.\ formal) group. 

\vspace{\vs}

%%%%%%  Subsection:  Group Functors  %%%%%%%%%%%%%%%%%%
%\subsubsection{Notations on Group Functors} 
\label{GroupFunctors}
%%%%%%%%  k-Functors and Completion  %%%%%%%%%%%%%%%%%
%\subsection{$k$-Functors and Completion} 
%\label{k-Functors}

Let $k$ be a ring (i.e.\ associative, commutative and with unit). 
$\Set$ denotes the category of sets, 
$\Ab$ the category of abelian groups. 
$\Algk$ denotes the category of $k$-algebras, and 
$\Artk$ the category of finite $k$-algebras (i.e.\ of finite length). 
%$\Cat$ denotes an arbitrary category. 
%
A \emph{$k$-functor} %(with values in $\Cat$) 
is by definition a covariant functor from $\Algk$ to $\Set$. 
%The category of $k$-functors is denoted by $\Fctr(\Algk,\Cat)$. 
%The morphisms %in $\Fctr(\Algk,\Cat)$ 
%are given by natural transformations of functors. 
%
A \emph{formal $k$-functor} %(with values in $\Cat$) 
is by definition a covariant functor from $\Artk$ to $\Set$. 
%The category of formal $k$-functors is denoted by $\Fctr(\Artk,\Cat)$. 
%
%The inclusion $\Artk \lra \Algk $ induces the \emph{completion functor} \\ 
%$\wh{\phtm} : \Fctr(\Algk,\Cat) \lra \Fctr(\Artk,\Cat)$. 
%The completion $\wh{\Fc}$ of a $k$-functor $\Fc$ is given by \; 
%$\wh{\Fc}\lrpn{\Ring} = \Fc\lrpn{\Ring}$ \; for $\Ring\in\Artk$. 
%
A (formal) $k$-functor with values in $\Ab$ is called a 
\emph{(formal) $k$-group functor}. 
%The category of $k$-group functors %with values $\Ab$ 
%will be denoted by $\Abkp$. 
%(rather than $\Fctr(\Algk,\Ab)$). 

%%%%%%%% Algebraic Groups %%%%%%%%%%%%%%%%%%%%%

A \emph{$\fld$-group} (or \emph{$\fld$-group scheme}) 
is by definition a $\fld$-group functor with values in $\Ab$ 
%(we only consider commutative group schemes) 
whose underlying set-valued $\fld$-functor 
is represented by a $\fld$-scheme. 
The category of $\fld$-groups is denoted by $\Gk$, 
the category of affine $\fld$-groups by $\Gak$. 
An \emph{algebraic $\fld$-group} 
(or just just \emph{algebraic group}) 
is a $\fld$-group whose underlying scheme 
is separated and of finite type over $\fld$. 
The category of algebraic $\fld$-groups is denoted by $\aGk$, 
the category of affine algebraic $\fld$-groups by $\aGak$.

%%%%%%%%  Formal Groups  %%%%%%%%%%%%%%%%%%%%%%%
%\subsubsection{Formal Groups} 
\label{sub:Formal-Groups}

Now let $k$ be a field. 
A \emph{formal $\fld$-scheme} is by definition 
a formal $\fld$-functor with values in $\Set$ 
which is the limit of a directed inductive system 
of finite $\fld$-schemes. 
%
%A formal $\fld$-functor $\fmlF$ is represented 
%by a formal $\fld$-scheme if there exists a 
%directed projective system $(\Algbr_i)$ of finite $\fld$-rings 
%and %a functorial (in $\Ring\in\Algk$) isomorphisms 
%an isomorphism of functors 
%$\fmlF \cong \varinjlim \Spec\Algbr_i$. 
%
%\bDef 
Let $\Algbrc$ be a \emph{profinite} $\fld$-algebra. 
%i.e.\ a complete topological $\fld$-algebra whose topology 
%has a basis of neighbourhoods of zero consisting of ideals 
%of finite codimension; this means $\Algbrc$ is the projective limit 
%(as a topological ring) of discrete quotients which are finite 
%$\fld$-algebras. 
%
The \emph{formal spectrum of $\Algbrc$} is %by definition 
the formal $\fld$-functor that assigns to $\Ring\in\Artk$ 
the set of continuous homomorphisms of $\fld$-algebras 
from the topological ring $\Algbrc$ to the discrete ring $\Ring$: \quad
$\Spf\Algbrc \lrpn{\Ring} = \Homcka(\Algbrc,\Ring)$. 
%\eDef 

%\bPrp 
%\label{formalScheme}
%For a formal $\fld$-functor $\fmlF$ the following conditions are equivalent: \\
%\begin{tabular}{rl} 
%(i) & $\fmlF$ is represented by a formal $\fld$-scheme. \\
%(ii) & There is a profinite $\fld$-algebra $\Algbrc$ and an isomorphism 
%%of functors \\
%         \, $\fmlF \,\cong\, \Spf\Algbrc$. \\
%(iii) & $\fmlF$ is left-exact, 
%i.e.\ $\fmlF$ commutes with finite projective limits. 
%\end{tabular} 
%\ePrp 
%(See \cite[I, No.~6]{D} or \cite[I, \S~4]{Fo}.)
%
%\vspace{\vs}

A \emph{formal $k$-group} (or just \emph{formal group}) 
is a formal $\fld$-group functor with values in $\Ab$ 
%(we are only interested in commutative formal groups) 
whose underlying set-valued formal $\fld$-functor 
is represented by a formal $\fld$-scheme, 
or equivalently is isomorphic to $\Spf\Algbrc$ 
for some profinite $\fld$-algebra $\Algbrc$. 
The category of formal $\fld$-groups is denoted by $\Gfk$. 

\bRmk
\label{fmlG in Abk}
A formal $\fld$-group $\fmlG: \Artk\lra\Ab$ extends in a natural way 
to a $\fld$-group functor \,$\wt{\fmlG}: \Algk\lra\Ab$, 
by defining $\wt{\fmlG}(\Ring)$ for $\Ring\in\Algk$ 
as the inductive limit of the $\fmlG(\Sfin)$, where $\Sfin$ ranges over the 
finite $\fld$-subalgebras of $\Ring$. 
If $\fmlG = \Spf\Algbrc$ for some profinite $\fld$-algebra $\Algbrc$, 
then $\wt{\fmlG}(\Ring) = \Homcka(\Algbrc,\Ring)$ for every $\Ring\in\Algk$. 
%In the following, we will denote $\fmlG$ and $\wt{\fmlG}$ by the same letter.
\eRmk 

%\bThm 
%\label{Gfk abelian Cat}
%The category $\Gfk$ of formal $\fld$-groups is abelian. 
%\eThm 
%
%\bPf 
%\cite[$\textrm{VII}_{\textrm{B}}$, 2.4.2, p.~521]{SGA3}. 
%\ePf 
%
%(This is equivalent to the fact that the category $\Gak$ is abelian, 
%by Cartier-duality, see \ref{Cartier-Duality}.) 

\bThm 
\label{Str-fml-Grp} 
A formal $\fld$-group $\fmlG$ is canonically an extension 
of an \'etale formal $\fld$-group $\fmlG_{\et}$ 
by an infinitesimal (= connected) formal $\fld$-group 
(i.e.\ the formal spectrum of a local ring) $ \fmlG_{\inf}$. 
%A formal $\fld$-group $\fmlG$ admits a canonical decomposition 
%\bDpl 
%0 \lra \fmlG_{\inf} \lra \fmlG \lra \fmlG_{\et} \lra 0 \laurink 
%\eDpl 
Here $\fmlG_{\et}\lrpn{\Rfin} = \fmlG\lrpn{\Rfin_{\red}}$ and 
$\fmlG_{\inf}\lrpn{\Rfin} = 
 \ker\bigpn{\fmlG\lrpn{\Rfin} \lra \fmlG\lrpn{\Rfin_{\red}}}$ 
for $\Rfin\in\Artk$, where $\Rfin_{\red} = \Rfin/\Nil(\Rfin)$. 
If the base field $\fld$ is perfect, there is a unique isomorphism 
\bDpl 
\fmlG \,\cong\, \fmlG_{\inf} \times_{\fld} \fmlG_{\et} \laurin 
\eDpl 
%$\fmlG_{\et}$ is called \emph{\'etale part} and 
%$\fmlG_{\inf}$ is called \emph{infinitesimal part} of $\fmlG$. 
\eThm 

\bPf 
\cite[I, No.~7, Prop.\ on p.~34]{D} or 
\cite[I, 7.2, p.~46]{Fo}. 
\ePf 

\vspace{\vs}

%%%%%%%% Sheaves of Abelian Groups %%%%%%%%%%%%%%%%%
%\subsubsection{Sheaves of Abelian Groups} 
\label{Sheaves}

Let $\bring$ be a ring. 
%Let $\Cat$ be an arbitrary category. 
%
An \emph{$\bring$-sheaf} %(with values in $\Cat$) 
is a sheaf (of sets) %(of objects of $\Cat$) 
on $\Alg{\bring}$ for the topology fppf. 
An $\bring$-sheaf with values in $\Ab$ is called 
an \emph{$\bring$-group sheaf}. 
The category of $\bring$-group sheaves is denoted by $\AbS{\bring}$. 
%The category $\AbS{\bring}$ is by definition a full subcategory 
%of the category of $\bring$-group functors $\AbSp{\bring}$. 

Let $\fld$ be a field. 
%\bPrp 
\label{Gk emb Abk}
The category of $\fld$-groups $\Gk$ and 
the category of formal $\fld$-groups $\Gfk$ 
are full subcategories of $\Abk$. 
%\ePrp 
This can be seen as follows: 
%\bPf 
A $\fld$-functor that is represented by a scheme 
is a sheaf, see \cite[III, \S~1, 1.3]{DG}. 
This gives the sheaf property for $\Gk$ by definition. 
For $\Gfk$ we can reduce to this case by Remark~\ref{fmlG in Abk} 
and the fact that a formal $\fld$-group is the direct limit of finite $\fld$-schemes. 
%\ePf 

%%%%%%%%  Linear Group associated to a Ring  %%%%%%%%%%%%%%%%%%%%%%%%
\subsubsection{Linear group associated to a ring}
\label{LinGroup_Ring}

Let $k$ be a field. 

\bDef 
\label{Lin_R}
Let $\Ring$ be a $k$-algebra. 
The \emph{linear group associated to $\Ring$} 
is the Weil restriction 
\;$\Lin_{\Ring} := \WeilRes{\Ring}{k} \GmS{\Ring} := \WR{\Gm}{\Ring}$ 
of $\GmS{\Ring}$ from $\Ring$ to $k$. 
\eDef 

If $\Sfin$ is a finite $\fld$-algebra,  then $\Lin_{\Sfin}$ is an affine 
algebraic $\fld$-group, according to \cite[I, \S~1, 6.6]{DG}. 
%Proposition \ref{WeilRes_repr}. 
%Thus the completion of the functor $\Lin_{\ghost}: \Algk \lra \Abk $
%is a formal $\fld$-group functor with values in $\aGak$ 
%(we omit the $\wh{\empt}$ here): 
%\[ 
%\Lin_{\ghost} \::\: \Artk \lra \aGak \laurin 
%\] 

%Suppose now that the base field $\fld$ is perfect. 
%Every finite $\fld$-algebra $\Sfin$ is of the form 
%$\Sfin = \Sfin_{\et} \oplus \Nil(\Sfin)$, 
%where $\Sfin_{\et} \cong \Sfin_{\red} = \Sfin/\Nil(\Sfin)$, 
%i.e.\ we have an injective map $\Sfin_{\red} \inj \Sfin$. 
%Since $\Gm$ is left-exact and tensor-product over a field $\fld$ is exact, 
%the linear group functor $\Lin_{\ghost}$ is left-exact.\footnote{
%Accepting the Full Embedding Theorem \cite[Chp.~7, Thm.~7.14]{Fr} 
%one can rephrase this by saying 
%``$\Lin_{\ghost}$ is a formal group with values in $\aGak$''.} 
%This yields a splitting 
%\[ \Lin_{\Sfin} = \Trz_{\Sfin} \tms_{\fld} \Upf_{\Sfin} 
%\] 
%where $\Trz_{\Sfin} = \WR{\Gm}{\Sfin_{\red}}$ and  
%$\Upf_{\Sfin} = \ker\Big(\WR{\Gm}{\Sfin} \lra \WR{\Gm}{\Sfin_{\red}}\Big)$. 
%One can show that $\Trz_{\Sfin}$ is a torus over $\fld$ 
%and $\Upf_{\Sfin}$ is a smooth connected unipotent algebraic 
%$\fld$-group. 

\bLem 
\label{L subset L_S}
Let $k$ be a perfect field. 
Every affine algebraic $\fld$-group $L$ is isomorphic to a closed subgroup 
of \;$\Lin_{\Sfin}$\; for some $\Sfin\in\Artk$. 
\eLem 

\bPf 
By Galois descent we can reduce to the case that $k$ is algebraically closed. 
Every affine algebraic $\fld$-group $L$ is isomorphic to a closed subgroup 
of $\Gl_{\num}$ for some $\num\in\Nat$ 
($\seecite$\cite[3.4 Thm.\ p.~25]{W}). 
Let $\repr: L \lra \Gl_{\num}$ be a faithful representation. 
%Define $\Sfin$ as the group algebra of $\repr(L)$, 
%i.e.\ $\Sfin := \fld\big[\repr(L)(\fld)\big]$. 
Define $\Sfin$ to be the group algebra of $\repr(L)$, 
i.e\ the $\fld$-subalgebra 
of the algebra of $(\num\tms\num)$-matrices $\Mat_{\num\times\num}(\fld)$ 
generated by $\repr(L)(\fld)$. 
In particular, $\Sfin$ is finite dimensional. 
Here we may assume that $L$ is reduced, 
hence determined by its $\fld$-valued points: 
otherwise embed the multiplicative part into $(\Gm)^t$ for some $t\in\Nat$ 
($\seecite$\cite[IV, \S~1, 1.5]{DG}) 
and the unipotent part into $(\Witt_r)^n$ for some $r,n\in\Nat$ 
($\seecite$\cite[V, \S~1, 2.5]{DG}), 
and replace $L$ by $(\Gm)^t \tms_{\fld} (\Witt_r)^n$. 
Then $\repr(L)(\fld)$ is contained in the unit group of $\Sfin$, 
and $\repr:L \lra \WR{\Gm}{\Sfin} = \Lin_{\Sfin}$ 
is a monomorphism from $L$ to $\Lin_{\Sfin}$. 
\ePf 

%\newpage

%%%%%%%%%%%%  Cartier Duality  %%%%%%%%%%%%%%%%%%%%%%%%%%%%%%
\subsubsection{Cartier duality} 
\label{Cartier-Duality}

Let $k$ be a field. 
We will use the functorial description of Cartier-duality 
as in \cite[II, No.~4]{D}. 
We may consider formal groups as objects of $\Abk$, 
cf.\ Remark \ref{fmlG in Abk}. 
%($\see$Subsection \ref{sub:Formal-Groups}).  
Let $G$ be a $\fld$-group sheaf. 
Let $\Homfabk\lrpn{G,\Gm}$ be the $\fld$-group functor defined by 
\bDpl 
\Ring \lmt \Homabs{\Ring}\big(G_{\Ring},\GmS{\Ring}\big) \laurink 
\eDpl 
which assigns to a $k$-algebra $\Ring$ 
the group of homomorphisms of $\Ring$-group sheaves 
from $G_{\Ring}$ to $\GmS{\Ring}$. 

\bThm 
\label{Cartier-dual} 
If $G$ is an affine group (resp.\ formal group), the $\fld$-group functor 
$\Homfabk\lrpn{G,\Gm}$ is represented by 
a formal group (resp.\ affine group) $\Gd$, 
which is called the \emph{Cartier dual} of $G$. \\
Cartier duality is an anti-equivalence between the category
of affine groups $\Gak$ and the category of formal groups $\Gfk$.
The functors $L\lmt\Ld$ and $\fmlG\lmt\fmlG^{\vee}$ are
quasi-inverse to each other. 
\eThm 

\bPf 
See \cite[$\textrm{VII}_{\textrm{B}}$, 2.2.2]{SGA3} 
or \cite[I, 5.4, p.~37]{Fo} for a description of Cartier duality via bialgebras. 
See \cite[II, No.~4, Thm p.~27]{D} for one direction of the functorial description 
of Cartier duality (it is only one direction since formal groups and affine groups 
are not considered as objects of the same category there). 
According to No.~\ref{Gk emb Abk} and the properties of the group functor 
$\Homfabk\lrpn{G,\Gm}$ as described in \cite[II, \S~1, 2.10]{DG}, 
it is an easy exercise to invert the given direction $L\lmt\Ld$ 
of the functorial description 
(one has to replace affine groups by formal groups, 
$\tens_{\fld}$ by $\tensc_{\fld}$, $\Homcka$ by $\Homka$ 
and $\Spf$ by $\Spec$). 
\ePf 

\bLem 
\label{Ld(R)} 
Let $L$ be an affine group and $\Ring$ a $\fld$-algebra. 
The $\Ring$-valued points of the Cartier-dual of $L$ are given by 
\bDpl 
\Ld(\Ring) \;=\; \Homabk\lrpn{L,\Lin_{\Ring}} \laurin 
\eDpl 
\eLem 

\bPf 
The statement is due to the fact that Weil restriction 
is right-adjoint to base extension: 
%($\see$Proposition \ref{adjoint functors}): 

\hspace{5mm} 
\bDpl 
\Ld(\Ring) = \Homabs{\Ring} \lrpn{ L_{\Ring},\GmS{\Ring} } 
 = \Homabk \big( L,\WR{\GmS{\Ring}}{\Ring} \big) \laurin 
%= \Homabk \lrpn{ L,\Lin_{\Ring} } 
\eDpl 
\hfill 
\ePf

\subsubsection*{Cartier dual of a multiplicative group} 
\label{sub:T dual} 

\bPrp 
\label{DualMultiplicative}
Let $L$ be an affine $\fld$-group. 
Then $L$ is multiplicative if and only if 
the Cartier-dual $\Ld$ is an \'etale formal $\fld$-group. 
\ePrp 

\bPf 
\cite[II, No.~8]{D}. 
\ePf 

\bExm 
In particular, 
the Cartier-dual of a split torus $\Trs\cong\lrpn{\Gm}^{t}$
is a lattice of the same rank: $\Trsd\cong\Zint^{t}$, 
i.e.\ a torsion-free \'etale formal group. 
\eExm 

\bPrp 
\label{DualMltAlg}
Let $L$ be a multiplicative $\fld$-group. 
Then $L$ is algebraic 
if and only if \,$\Ld(\clfld)$\, is of finite type. 
\ePrp 

\bPf 
\cite[IV, \S~1, 1.2]{DG}. 
\ePf

\subsubsection*{Cartier dual of a unipotent group} 
\label{sub:V dual} 

\bPrp 
\label{DualUnipotent}
Let $L$ be an affine $\fld$-group. 
Then $L$ is unipotent if and only if 
the Cartier-dual $\Ld$ is an infinitesimal formal $\fld$-group. 
\ePrp 

\bPf 
\cite[II, No.~9]{D}. 
\ePf 

\bExm[Cartier duality of Witt vectors] 
\label{WittVectors}
Suppose $\chr(\fld)=p>0$. 
Let $\Witt$ denote the $\fld$-group of Witt-vectors, 
$\Witt_r$ the $\fld$-group of Witt-vectors of finite length $r$. 
Let $\Wittc$ be the completion of $\Witt$ at $0$, 
i.e.\ $\Wittc$ is the subfunctor of $\Witt$ 
that associates to $\Ring\in\Algk$ the set of 
$(w_0,w_1,\ldots) \in \Witt(\Ring)$ such that $w_{\nu} \in \Nil(\Ring)$ 
for all $\nu\in\Nat$ and $w_{\nu} = 0$ for almost all $\nu\in\Nat$. 
%\[ \Wittc\pn{R} = 
%   \left\{ \lrpn{w_{\nu}}_{\nu\in\Nat} \in \Witt(R) \;\left|\;  
%   \begin{array}{l} 
%   w_{\nu} \in\Nil(R) \; \textrm{ for all } \nu\in\Nat  \\ 
%   w_{\nu} = 0 \; \textrm{ for almost all } \nu\in\Nat 
%   \end{array}
%   \right.\right\}
%\] 
Moreover let $\Wcfl{r} = \ker\big(\Frob^r:\Wittc\lra\Wittc^{\lrpn{p^r}}\big)$ 
be the kernel of the $r^{\textrm{th}}$ power of the Frobenius $\Frob$. 
%\[ \Frob:\Witt\lra\Witt^{\lrpn{p}} \laurink  \hspace{4em}
%   \lrpn{w_{\nu}}_{\nu\in\Nat} \lmt \lrpn{w_{\nu}^p}_{\nu\in\Nat} 
Let $\Wpl$ denote the affine $k$-group which associates with $R \in \Algk$ 
the multiplicative group $1 + t R[[t]]$ of formal power series in $R$. 
Let $\Eah$ be the series 
\bDpl  
\Eah\lrpn{t} 
 \, = \, \exp \lrpn{ -\sum_{r\geq 0} \frac{t^{p^{r}}}{p^{r}}} 
 \, = \, \prod_{\substack{r\geq 1 \\ (r,p)=1}} \lrpn{1-t^{r}}^{\moeb(r)/r} \laurink
\eDpl 
where $\moeb$ denotes the M\"obius function. 
%i.e.\ 
%\[ \moeb\lrpn{r} = 
%   \left\{ 
%      \begin{array}{cl}
%          0 & \textrm{ if $r$ is divisible by the square of a prime} \laurink \\
%  (-1)^n & \textrm{ if $r=p_1\cdots p_n$ 
%                  and $p_1,\ldots,p_n$ are distinct primes} \laurink \\
%           1& \textrm{ if $r=1$} \laurin 
%      \end{array}
%   \right. 
%\] 
The \emph{Artin-Hasse exponential} is the homomorphism of $k$-groups 
\bDpl 
\Expah: \Witt \lra \Wpl 
\eDpl 
defined by 
\bDpl 
\Expah\pn{w,t} := \Expah\pn{w}\pn{t} := \prod_{r\geq 0} \Eah\lrpn{w_r t^{p^r}} \laurin 
\eDpl 
%It satisfies for Witt vectors $w,v$ 
%\[ \Expah\lrpn{w+v} \,=\, \Expah\lrpn{w} \cdot \Expah\lrpn{v} \laurin 
%\] 
For details see e.g.\ \cite[III, No.s~1 and 2]{D}. 

Then $\Wittc$ (resp.\ $\Wcfl{r}$) is Cartier-dual to 
$\Witt$ (resp.\ $\Witt_r$), 
and the pairing 
\bDpl 
\lrpair{\llul,\lull}: \Wittc \tms \Witt \lra \Gm 
\hspace{0mm} \lrpn{\textrm{resp.} \hspace{0mm} 
\lrpair{\llul,\lull}: \Wcfl{r} \tms \Witt_r \lra \Gm} 
\eDpl 
is given by 
\bDpl 
\lrpair{v,w} 
 = \Expah\lrpn{v \cdot w, 1} 
 = \prod_{\substack{r\geq 0 \\ s\geq 0}} \Eah\lrpn{v_r^{p^s} w_s^{p^r}} 
\laurin 
\eDpl 
See \cite[V, \S~4, Prop.~4.5 and Cor.~4.6]{DG}. 
\eExm 

\bPrp 
\label{DualUptAlg_p} 
Suppose $\chr(k) = p > 0$. 
Let $L$ be an affine $\fld$-group. 
The following conditions are equivalent: 

\begin{tabular}{rl}
(i) & $L$ is unipotent algebraic. \\ 
(ii) & There is a monomorphism $L \,\inj \lrpn{\Witt_r}^n$ 
       for some $r,n \in \Nat$. \\ 
(iii) & There is an epimorphism $\pn{\Wcfl{r}}^n \sur \Ld$ 
        for some $r,n \in \Nat$. 
\end{tabular} 

\noindent 
(Here we use the notation from Example \ref{WittVectors}.) 
\ePrp 

\bPf 
(i)$\Lra$(ii) \cite[V, \S~1, 2.5]{DG}. 

(ii)$\Lra$(i) The underlying $\fld$-scheme of $\Witt_r$ is 
the affine space $\Afn^r$, thus $\Witt_r$ is algebraic. 
\;$ 0 = \Witt_0 \subset \Witt_1 \subset \Witt_2 \subset \ldots \subset \Witt_r 
$\; 
is a filtration of $\Witt_r$ with quotients 
$\Witt_{\nu}/\Witt_{\nu-1} \cong \Witt_1 = \Ga$, 
hence $\Witt_r$ is unipotent, 
according to \cite[IV, \S~2, 2.5]{DG}. 
Products of unipotent groups and closed subgroups of a unipotent group 
are unipotent by \cite[IV, \S~2, 2.3]{DG}. 
%Moreover, closed subschemes of algebraic schemes are algebraic, 
%as a quotient of a finitely generated $\fld$-algebra is 
%finitely generated. 
Since $L$ is isomorphic to a closed subgroup of $\lrpn{\Witt_r}^n$, 
it is unipotent and algebraic. 

(ii)$\Llra$(iii) This is due to Cartier-duality of Witt vectors, 
see Exm.~\ref{WittVectors}. 
%the Cartier-dual of $\Witt_n$ is isomorphic to $\Wcfl{n}$, 
%see \cite[V, \S~4, 4.5]{DG}. 

\hfill 
\ePf

\subsubsection{Dual abelian variety} 
\label{dualAbVar}

Let $k$ be a field. %an algebraically closed field. 
Let $A$ be an abelian variety over $k$. 
The dual of $A$ is given by $\Ad=\Pic^0 A$. 
According to the generalized Barsotti-Weil formula 
($\seecite$\cite[III.18]{O}), 
the dual abelian variety $\Ad$ 
represents the $k$-group sheaf $\Extfabk(A,\Gm)$, associated to 
\bDpl  
\Ring \;\lmt\; \Extabs{\Ring} \lrpn{A_{\Ring},\GmS{\Ring}} \laurin 
\eDpl 

\bPrp 
\label{Xeil}
Let $A$ be an abelian $\fld$-variety and $\Sfin$ a finite $k$-algebra. 
There is a canonical isomorphism 
\bDpl 
%\Xeil{\Sfin}{\fld}: 
\Extabk\pn{A,\Lin_{\Sfin}} \iso \Extabs{\Sfin}\lrpn{A_{\Sfin},\GmS{\Sfin}} \laurin 
\eDpl 
%induced by Weil restriction of the fibres. \\ 

\noindent 
Thus the $\Sfin$-valued points of the dual abelian variety are given by 
\bDpl 
\Ad(\Sfin) = \Extabk(A,\Lin_{\Sfin}) \laurin 
\eDpl 
\ePrp 

\bPf 
Consider the following composition of functors on $\AbS{\Sfin}$: \\ 
\bDpl  \Homabk\pn{A,\lull} \circ \WeilRes{\Sfin}{k}: 
G \lmt \Homabk\bigpn{A,\WR{G}{\Sfin}} = \Homabs{S}\pn{A_{\Sfin},G} \laurin 
\eDpl 

\noindent 
Since 
$\Extabk(A,\Lin_{\Sfin})$ resp.\ $\Extabs{\Sfin}\lrpn{A_{\Sfin},\GmS{\Sfin}}$ 
is identified with the set of primitive elements in 
$\H^1\big(A, \Lin_{\Sfin}(\sO_{A})\big)$ resp.\ 
$\H^1\big(A_{\Sfin}, \Gm(\sO_{A_{\Sfin}})\big)$ 
(see \cite[VII, No.~15, Thm.~5]{S} resp.\ \cite[III.17.6]{O}), 
we may compute these Ext-groups 
using the \'etale site instead of the flat site, 
according to \cite[III, Thm.~3.9]{Mi}. 
As $\Sfin$ is a finite $k$-algebra, the Weil restriction 
$\WeilRes{\Sfin}{k}: G \lmt \WR{G}{\Sfin}$ is exact for the \'etale topology 
($\seecite$\cite[II, Cor.~3.6]{Mi}). 
Then the exact sequence of low degree terms 
of the Grothendieck spectral sequence 
yields a canonical isomorphism 
$\bigpn{\R^1\Homabk\pn{A,\lull}}\bigpn{\WeilRes{\Sfin}{k}(\Gm)} \iso  
\R^1\bigpn{\Homabk\pn{A,\lull} \circ \WeilRes{\Sfin}{k}}(\Gm)$ 
($\seecite$\cite[App.~B, Cor.~2]{Mi}), 
showing the statement. 
\ePf

\subsubsection{Extensions of formal groups} 
\label{Ext fml grps}

Let $k$ be a field. 

\bLem 
\label{Ext(FmlG,Gm)}
\bDpl  \Extfabk(\sF,\Gm) = 0 
\eDpl 
for any dual-algebraic formal $\fld$-group $\fmlG$. 
\eLem

\bPf 
Let $L$ be the affine algebraic group dual to $\fmlG$. 
Let $R$ be a $\fld$-algebra. 
We show that $\Extabs{R}(\Ld_R,\GmS{R}) = 0$ 
flat locally over $\Spec R$. 
As $L$ has a filtration 
\bDpl  0 = L_0 \,\subset\, L_1 \,\subset\, \ldots \,\subset\, L_r = L 
\eDpl  
with quotients equal to $\Gm$, $\Ga$ or a finite $\fld$-group, 
%and an exact sequence $0 \ra \fmlG' \ra \fmlG \ra \fmlG'' \ra 0$ 
%induces an exact sequence 
%$ \Extf^1(\fmlG'',\Gm) \ra \Extf^1(\fmlG,\Gm) \ra \Extf^1(\fmlG',\Gm) $. 
%Therefore it suffices to show that 
%$\Extf^1(\Gm^{\vee},\Gm) = 0$ and $\Extf^1(\Ga^{\vee},\Gm) = 0$. 
it suffices to show the statement for $L = \Gm$, $\Ga$ or a finite group. 

If $L$ is a finite $\fld$-group, 
the Cartier dual $\Ld$ is again a finite $\fld$-group $F$, 
and $\Extabs{R}(F_R,\GmS{R}) = 0$ flat locally 
according to \cite[III, \S~4, Lem.~4.17]{Mi}. 

If $L = \Gm$, then $\Ld = \Zint$, 
and $\Extabs{R}(\Zint,\Gm) = 0$ is clear. 

If $L = \Ga$ and $\chr(\fld) = 0$, then $\Ld = \Gac$. 
It holds $\Extabs{R}\pn{\Gac,\Gm} = \Extabs{R}\pn{\Gac,\Gac} = 0$, 
see \cite[Lem.~A.4.6]{BaBe}. 
%By definition $\Extabs{R}\pn{\Gac,\lull} = \R^1\Homabs{R}\pn{\Gac,\lull}$ 
%is the first right derived functor of $\Homabs{R}\pn{\Gac,\lull}$. 
%Any homomorphism $\Gac \lra \Gm$ factors uniquely through 
%$\wh{\bG}_m = \Gac$, 
%therefore $\Homabs{R}\pn{\Gac,\Gm} = \Homabs{R}\pn{\Gac,\Gac}$. 
%As the category $\GfS{R}$ admits enough injectives, one finds a compatible 
%injective resolution $\wh{I}^{\bullet} \ra J^{\bullet}$ of $\Gac \ra \Gm$, 
%such that $\wh{I}^{\nu}$ is the completion of $J^{\nu}$ for all $\nu \in \Nat$. 
%Then $\Extabs{R}\pn{\Gac,\Gm} = \Extabs{R}\pn{\Gac,\Gac} = 0$, 
%since any extension of $\Gac$ by $\Gac$ is again an infinitesimal formal group, 
%and any infinitesimal formal group is a product of $\Gac$'s in characteristic 0. 
%%($\see$Proposition \ref{infGrp-LieAlg}). 
%See also \cite[Lem.~A.4.6]{BaBe}. 

If $L = \Ga = \Witt_1$ and $\chr(\fld) = p > 0$, then $\Ld = \Wcfl{1}$ 
(notation from Example \ref{WittVectors}). 
Since $\Wcfl{1} = \ker\bigpn{\Frob: \Wittc \ra \Wittc}$ 
is annihilated by the Frobenius $\Frob$ 
and $\ker\bigpn{\Frob: \Gm \ra \Gm} = \mu_p$ 
is the group of $p^{\textrm{th}}$ roots of unity over $\Spec R$ 
(this group is finite, hence both an algebraic group and a formal group), 
any extension $E \in \Extabs{R}\pn{\Wcfl{1},\Gm}$ is the push-out 
of an extension $\fmlG \in \Extabs{R}\pn{\Wcfl{1},\mu_p}$. 
As $\mu_p$ and $\Wcfl{1}$ are base changes of formal $k$-groups, 
the affine algebra $\sO\pn{\mu_p}$ of $\mu_p$ 
is a free $R$-module of finite rank 
and the affine algebra $\sO\pn{\Wcfl{1}}$ of $\Wcfl{1}$ 
is the projective limit of free $R$-modules of finite rank; 
I will refer to those algebras as \emph{free pro-finite-rank $R$-algebras}. 
The underlying $\mu_p$-bundle of $\sF$ is flat locally trivial, 
hence flat locally the affine algebra of $\sF$ is 
$\sO\pn{\sF} = \sO\pn{\mu_p} \tens_R \sO\pn{\Wcfl{1}}$, 
and this is a free pro-finite-rank $R$-algebra as well. 
In this case Cartier duality works in the same way as for formal $k$-groups 
\footnote{Note that the category of flat locally free pro-finite-rank $R$-algebras 
is not abelian. 
The references for Cartier duality listed in the proof of Theorem \ref{Cartier-dual} 
make additional assumptions on the base ring $R$ in order to achieve 
that the category of $R$-formal groups is abelian.}, 
so the exact sequence 
\[ 0 \lra \mu_p \lra \fmlG \lra \Wcfl{1} \lra 0 
\] 
is turned into the following exact sequence of $R$-groups 
\[ 0 \lra \Ga \lra \fmlGd \lra \Zint/p\Zint \lra 0 
\] 
where $\fmlGd = \Homfabs{R}\pn{\fmlG,\Gm}$. 
%For $x \in \fmlGd$ we have $p x \in \ker\bigpn{\fmlGd \lra \Zint/p\Zint} = \Ga$, 
%and $\Ga$ is killed by $p$ since $R$ is of characteristic $p$. 
%Thus $\fmlGd$ is killed by $p^2$. 
Applying %the functor 
$\Homfabs{R}\pn{\llul,\Gm}$ 
to the push-out diagram $\Gm \la \mu_p \ra \fmlG$ of $E$ 
shows that $\Ed := \Homfabs{R}\pn{E,\Gm}$ 
is the pull-back of the diagram $\Zint \ra \Zint/p\Zint \la \fmlGd$. 
In particular, since $\fmlGd \lra \Zint/p\Zint$ is surjective, 
$\Ed \lra \Zint$ is surjective as well. %($\seecite$\cite[Pullback Theorem 2.54]{Fr}). 
Thus we obtain an exact sequence 
\[ 0 \lra \Ga \lra \Ed \lra \Zint \lra 0 
\] 
which is obviously split. 
Dualizing again gives the split exact sequence 
\[ 0 \lra \Gm \lra E^{\vee \vee} \lra \Wcfl{1} \lra 0 
\] 
where $E^{\vee \vee} = \Homfabs{R}\pn{\Ed,\Gm}$. 
The canonical map of any abelian sheaf $\sA$ to its double dual 
$\sA^{\vee \vee}$ yields the following commutative diagram with exact rows 
\[ \xymatrix{ 
0 \ar[r]  & \Gm \ar[r] \ar@{=}[d] & E \ar[r] \ar[d] & \Wcfl{1} \ar[r] \ar@{=}[d] & 0 \\ 
0 \ar[r]  &  \Gm \ar[r]  &  E^{\vee \vee} \ar[r]  &  \Wcfl{1} \ar[r]  &  0 \\ 
   }
\] 
where the vertical arrow in the middle is an isomorphism by the Five Lemma. 
Thus $E \cong E^{\vee \vee}$ is split. 
\ePf

\subsection{1-motives with unipotent part} 
\label{1-MotUpt} 

%Let $\fld$ be an algebraically closed field. 
Let $\fld$ be a field. 
%$\clfld$ an algebraic closure. 
%$\sepfld$ the separable closure in $\clfld$. 

\subsubsection{Definition of a 1-motive with unipotent part} 
\label{Anatomy}

%A \emph{1-motive with unipotent part}, 
%which we will define in this Subsection, 
%is the analogue to a Laumon 1-motive \cite[D\'efinition (5.1.1)]{L} 
%for arbitrary characteristic extended by torsion. 

\bDef 
\label{fmlGrp_dual-algebraic} 
A formal $\fld$-group $\fmlG$ is called \emph{dual-algebraic} 
if its Cartier-dual $\fmlGd$ is algebraic. 
The category of dual-algebraic formal $\fld$-groups is denoted by $\dGfk$. 
\eDef 

\bPrp 
\label{dual-algebraic} 
A formal $\fld$-group $\fmlG$ is dual-algebraic 
if and only if the following conditions are satisfied: 

\begin{tabular}{rl} 
(1) & $\fmlG\pn{\clfld}$ is of finite type$\laurink$ \\ 
(2) & for $\chr(\fld)=0$: \,
      $\Lie\left(\fmlG\right)$ is finite dimensional$\laurink$ \\ 
    & for $\chr(\fld)>0$: \,
      $\fmlG_{\inf}$ is a quotient of $\bigpn{\Wcfl{r}}^n$ 
      for some $r,n\in\Nat$ \\ 
 & \phantom{for $\chr(\fld)>0$: \,} 
 (see Example \ref{WittVectors} 
 %Proposition \ref{DualUptAlg_p} 
 for the definition of $\Wcfl{r}$). 
\end{tabular} 
\ePrp 

\bPf 
The decomposition of $\fmlG$ 
into \'etale part $\fmlG_{\et}$ and infinitesimal part $\fmlG_{\inf}$ 
gives the decomposition of the affine group $\fmlGd$ 
into multiplicative part and unipotent part, according to 
Propositions \ref{DualMultiplicative} and \ref{DualUnipotent}. 
Then that statement follows directly from Propositions 
\ref{DualMltAlg} and \ref{DualUptAlg_p} for $\chr(k)>0$. 
For $\chr(k)=0$, the assertion in (2) is due to the fact that 
the Lie functor yields an equivalence between 
the category of commutative infinitesimal formal $k$-groups 
and the category of $k$-vector spaces, 
see \cite[$\mathrm{VII}_{\mathrm{B}}$, 3.2.2]{SGA3}. 

\hfill 
\ePf 

\bLem 
\label{subFmlGroup_algebraic}
Let $\fmlG$ be a dual-algebraic formal group. 
Then any formal group $\fmlGr$ 
that is a subgroup or a quotient of $\fmlG$ 
is also dual-algebraic. 
\eLem 

\bPf 
By Cartier-duality, this is equivalent to 
the dual statement about affine algebraic groups, 
see \cite[II, No.~6, Cor.~4 of Thm.~2, p.~32]{D}. 
%Lemma \ref{quotAffGroup_algebraic} below. 
\ePf 

%\bLem 
%\label{quotAffGroup_algebraic}
%Let $L$ be an affine algebraic group. 
%Then any affine closed subgroup $K$ of $L$ 
%and any affine quotient $N$ of $L$ 
%is also algebraic. 
%\eLem 
%
%\bPf 
%\cite[II, No.~6, Cor.~4 of Thm.~2, p.~32]{D}. 
%\ePf 

%\bPf 
%Let $N$ be the quotient of the closed embedding $K \lra L$, 
%resp.\ let $K$ be the kernel of the quotient homomorphism $L \lra N$. 
%We obtain an exact sequence of algebraic groups 
%\bDpl 0 \lra K \lra L \lra N \lra 0 \laurin 
%\eDpl 
%As the category $\Gak$ of affine groups is abelian, 
%$K$ is an affine group if and only if $N$ is an affine group. 
%A closed subscheme of an algebraic scheme is algebraic. 
%Therefore $K$ is algebraic, since $L$ is. 
%As the category $\aGk$ of algebraic groups is abelian, 
%$N$ is algebraic. 
%\ePf 

\bDef 
\label{Def 1-motive}
A \emph{1-motive with unipotent part} is a tuple 
$M = \lrpn{\fmlG,L,A,G,\mu}$, 
where 
%is a complex concentrated in degrees $-1$ and $0$, 
%i.e.\ a homomorphism,  
%in the category of $\fld$-group sheaves 
%of the form $M = \left[ \sF \overset{\mu} \lra G \right]$, satisfying 

\begin{tabular}{cl}
(a) & $\fmlG$ is a dual-algebraic formal group 
      ($\see$Definition \ref{fmlGrp_dual-algebraic})$\laurink$ \\ 
(b) & $L$ is an affine algebraic group$\laurink$ \\ 
(c) & $A$ is an abelian variety$\laurink$ \\ 
(d) & $G$ is an extension of $A$ by $L \laurink$ \\ %in $\Abk$ \\ 
%      ($\see$Definition \ref{algGrp_ExtType}) 
(e) & $\mu:\fmlG\lra G$ is a homomorphism in $\Abk \laurin$ 
%    &  in the category of $\fld$-group sheaves 
\end{tabular} 

\noindent 
A \emph{homomorphism between 1-motives with unipotent part} 
$M = (\fmlG,L,A,G,\mu)$ and $N = (\fmlE,\Lam,B,H,\nu)$ 
is a tuple $h = \pn{\phe,\lam,\alp,\gam}$ of homomorphisms 
$\phe:\fmlE\ra\fmlG, \lam:L\ra\Lam, \alp:A\ra B, \gam:G\ra H$, 
compatible with the structures of $M$ and $N$ 
as 1-motives with unipotent part, 
i.e.\ giving an obvious commutative diagram. 

%The category of 1-motives with unipotent part is denoted by $\Mot$. 

For convenience, we will refer to a \emph{1-motive with unipotent part} 
only as a \emph{1-motive}. 
\eDef

If $G$ is a smooth connected algebraic group, 
it admits a canonical decomposition $0 \ra L \ra G \ra A \ra 0$ 
as an extension of an abelian variety $A$ 
by a connected affine algebraic group $L$, 
according to the Theorem of Chevalley. %\ref{Thm Chevalley} 
Thus a homomorphism $\mu:\fmlG\lra G$ in $\Abk$ 
gives rise to a 1-motive $M = (\fmlG,L,A,G,\mu)$ 
that we will denote just by $M = \bigbt{\sF \overset{\mu}\lra G}$.

\subsubsection{Duality of 1-motives} 
\label{dual1-mot}

\bThm 
\label{Ext(A,L)=Hom(Ld,Ad)}
Let $L$ be an affine algebraic group and $A$ an abelian variety. 
There is a canonical isomorphism of abelian groups  
\[ \Pop:\,\Extabk\big(A,L\big)\,\iso\,\Homabk\big(\Ld,\Ad\big) \laurin 
\] 
\eThm 

\bPf 
Consider the following left exact functor on $\Abk$ 
\begin{align*} 
 F: G \lmt \Bilabk\pn{A,\Ld;G} 
              & = \Homabk\bigpn{A,\Homfabk\pn{\Ld,G}} \\ 
              & = \Homabk\bigpn{\Ld,\Homfabk\pn{A,G}} \laurink 
\end{align*} 
where $\Bilabk\pn{A,\Ld;G}$ is the group of $\Zint$-bilinear maps 
$A \tms \Ld \lra G$ of sheaves of abelian groups. 
The two ways of writing $F$ as a composite 
yield the following two spectral sequences: 
\begin{eqnarray*}
\Extabkt{p}\bigpn{A,\Extfabkt{q}\pn{\Ld,G}} & \Lra & \R^{p+q} F(G) \\ 
\Extabkt{p}\bigpn{\Ld,\Extfabkt{q}\pn{A,G}} & \Lra & \R^{p+q} F(G) 
\end{eqnarray*}
For $G = \Gm$, the associated exact sequences of low degree terms are 
\[ 0 \lra \Ext\bigpn{A,\Homf\pn{\Ld,\Gm}} \lra \R^1 F(\Gm) \lra 
\Hom\bigpn{A,\Extf\pn{\Ld,\Gm}} = 0 
\] 
where the last term vanishes due to Lemma \ref{Ext(FmlG,Gm)}, 
and 
\begin{align*} 
0 = \Ext\bigpn{\Ld,\Homf\pn{A,\Gm}} \lra \R^1 F(\Gm) & \lra  \Hom\pn{\Ld,\Ad} \\ 
 & \lra  \Extt{2}\pn{\Ld,\Homf\pn{A,\Gm}} = 0 \laurin 
\end{align*} 
%where the first and the last term vanish since $\Homf\pn{A,\Gm} = 0$. 
Putting these together we obtain isomorphisms 

\hspace{10mm} 
\bDpl 
\Extabk\pn{A,L} \iso \R^1 F(\Gm) \iso \Homabk\pn{\Ld,\Ad} \laurin 
\eDpl 
\hfill 
\ePf 

\bRmk[Explicit description of \;$\Ext(A,L) \iso \Hom(\Ld,\Ad)$] 
\label{Pop=conHom} 
The isomorphism $\Pop$ in Theorem \ref{Ext(A,L)=Hom(Ld,Ad)} 
sends $G \in \Ext(A,L)$ 
to the connecting homomorphism 
$\Homfabk(L,\Gm) \lra \Extfabk(A,\Gm)$ 
in the long exact cohomology sequence obtained 
from applying $\Homfabk(\llul,\Gm)$ 
to the short exact sequence $0 \ra L \ra G \ra A \ra 0$. 
Explicitly, this is the map $\Pop(G): \lam \lmt \lam_*G$, 
which sends $\lam \in \Ld(R) = \Homabk(L,\Lin_R)$ to the push-out 
$\lam_* G \in \Extabk\lrpn{A,\Lin_{\Rfin}} = \Ad(\Rfin)$ 
of the diagram $\Lin_R \overset{\lam} \lla L \lra G$. 

If $L = \Lin_{\Sfin}$ for some $\Sfin\in\Artk$, the inverse map of $\Pop$ 
is given by the map $\Pop^{-1}: \tha \lmt \tha(\id_L)$, 
which sends a homomorphism $\tha: \Ld \lra \Ad$ 
to the image $\tha(\id_L) \in \Extabk(A,L) = \Ad(S)$ 
of the identity $\id_L \in \Homabk(L,L) = \Ld(S)$. 
In general, a given homomorphism $\tha \in \Homabk(\Ld,\Ad)$ 
can be written as an element $\bt{\Ld \ra \Ad} \in \Compl^{[-1,0]}\pn{\Abk}$ 
of the category of two-term complexes in $\Abk$, 
with $\Ld$ placed in degree $-1$ and $\Ad$ in degree $0$. 
Then $\Pop^{-1}(\tha) = \Extfcttabk\bigpn{\bt{\Ld \ra \Ad},\Gm}$, 
and this $k$-group sheaf is represented by an algebraic group: 
the short exact sequence 
$0 \lra \Ad \lra \bt{\Ld \ra \Ad} \lra \Ld[1] \lra 0$ 
gives rise to the exact sequence 
$0 \ra \Homf(\Ld,\Gm) \ra \Extf\bigpn{\bt{\Ld \ra \Ad},\Gm} \ra \Extf(\Ad,\Gm) \ra 0$ 
since $\Extfabk(\Ld,\Gm) = 0$ by Lemma \ref{Ext(FmlG,Gm)}. 
Thus $\Extfcttabk\bigpn{\bt{\Ld \ra \Ad},\Gm}$ 
is an extension of $A$ by $L$. 
\eRmk 

\bDef 
\label{Def:dual1-mot}
The dual of a 1-motive $M = \pn{\fmlG,L,A,G,\mu}$ 
is the 1-motive $\Md = \lrpn{\Ld,\fmlGd,\Ad,H,\eta}$, 
where $H = \Extfcttabk\big([\fmlG\ra A],\Gm\big) = \Pop^{-1}\pn{\ol{\mu}}$ 
for $\ol{\mu}: \fmlG \overset{\mu} \lra G \lra A$ the composite, 
and $\eta: \Ld \lra H$ is the connecting homomorphism 
$\Homfabk\big(L,\Gm\big) \,\lra\, \Extfcttabk\big([\fmlG\ra A],\Gm\big)$ 
in the long exact cohomology sequence associated with 
$0 \lra \bt{0\ra L} \lra \bt{\fmlG\ra G} \lra \bt{\fmlG\ra A} \lra 0$. 
%\[ 0 \lra 
%   \left[ \begin{array}{c} 0 \\ \downarrow \\ L \end{array}\right] \lra 
%   \left[ \begin{array}{c} \fmlG \\ \downarrow \\ G \end{array}\right] \lra 
%   \left[ \begin{array}{c} \fmlG \\ \downarrow \\ A \end{array}\right] \lra 
%   0
%\] 
%Here $\Compl\pn{\Abk}$ is the category of complexes of sheaves of abelian groups, 
%and in $\bt{\fmlG\ra A}$, $\fmlG$ is placed in degree $-1$ 
%and $A$ in degree $0$. 
%(In this note we will only use duality of 1-motives of pure type 
%as in Definition \ref{Def:dual1-mot}.)
\eDef 

\bRmk 
The double dual $\Mdd$ of a 1-motive $M$ is canonically isomorphic to $M$. 
If $M$ is of the form $\bigbt{0\ra G} := \bigpn{0,L,A,G,0}$, then the dual is 
$\bigbt{\Ld \overset{\Pop\pn{G}}\lra \Ad} := 
\bigpn{\Ld,0,\Ad,\Ad,\Pop\pn{G}}$. 
If $M$ is of the form 
$\bigbt{\fmlG \overset{\mu}\lra A} := \bigpn{\fmlG,0,A,A,\mu}$, the dual is 
$\bigbt{0 \lra \Pop^{-1}\pn{\mu}} := 
\bigpn{0,\fmlGd,\Ad,\Pop^{-1}\pn{\mu},0}$. 
This is clear by Theorem \ref{Ext(A,L)=Hom(Ld,Ad)}, 
and these ``pure 1-motives'' are the only ones 
that we are concerned with in this note. 
For the general case, the proof carries over literally from \cite[(5.2.4)]{L}. 
\eRmk 

\bPrp 
\label{duality-functorial}
Duality of 1-motives is functorial, 
i.e.\ duality assigns to a homomorphism of 1-motives $h: M \lra N$ 
a dual homomorphism $\hd: \Nd \lra \Md$. 
\ePrp 

\bPf 
%Functoriality comes from the fact that duality of 1-motives is derived 
%from an application of the functor $\Homfcttabk(\llul,\Gm)$: 
%%cf.\ Remark \ref{Pop=conHom}: 
%
Let $M = (\fmlG,L,A,G,\mu)$ and $M' = (\fmlG',L',A',G',\mu')$ 
be 1-motives %with unipotent parts 
and $h: M \lra M'$ a homomorphism of 1-motives. 
Applying $\Homfcttabk(\llul,\Gm)$ 
to the commutative diagram with exact rows 
\[ \xymatrix{  
0 \ar[r] & [0 \ra L] \ar[d] \ar[r] & [\fmlG \ra G] \ar[d]^{h} \ar[r] & 
   [\fmlG \ra A] \ar[d] \ar[r] & 0 \\ 
0 \ar[r] & [0 \ra L'] \ar[r] & [\fmlG' \ra G'] \ar[r]  &  [\fmlG' \ra A'] \ar[r] & 0 } 
\] 
yields the homomorphism $\hd: [(L')^{\vee} \ra H'] \lra [\Ld \ra H]$ 
%given by 
%\[ \xymatrix{ 
%\Ld \ar[r] & H  \\ 
%(L')^{\vee} \ar[r] \ar[u] & H' \ar[u] 
% }
%\] 
where $H = \Extfcttabk\big([\fmlG \ra A],\Gm\big)$ 
and $H' = \Extfcttabk\big([\fmlG' \ra A'],\Gm\big)$. 
%and $\Ld = \Homfabk(L,\Gm)$, $\Lamd = \Homfabk(\Lam,\Gm)$. 
Applying $\Homfcttabk(\llul,\Gm)$ 
to the commutative diagram with exact rows 
\[ \xymatrix{  
0 \ar[r] & [0\ra A] \ar[d] \ar[r] & [\fmlG \ra A] \ar[d] \ar[r] & 
   [\fmlG \ra 0] \ar[d] \ar[r] & 0 \\ 
0 \ar[r] & [0 \ra A'] \ar[r] & [\fmlG' \ra A'] \ar[r]  &  [\fmlG' \ra 0] \ar[r] & 0 } 
\] 
%yields 
%\[ \xymatrix{  
%0 \ar[r] & \fmlGd \ar[r] & H \ar[r] & \Ad \ar[r] & 0 \\ 
%0 \ar[r] & (\fmlG')^{\vee} \ar[r] \ar[u] & H' \ar[r] \ar[u] & A'^{\vee} \ar[r] \ar[u] & 0 } 
%\] 
%where 
%$\fmlGd = \Homfabk(\fmlG,\Gm) = \Extfcttabk\big([\fmlG \ra 0],\Gm\big)$, \\
%$(\fmlG')^{\vee} = \Homfabk((\fmlG')^{\vee},\Gm) = 
%\Extfcttabk\big([(\fmlG')^{\vee} \ra 0],\Gm\big)$ and \\
%$\Ad = \Extfabk(A,\Gm)$, $(A')^{\vee} = \Extfabk(A',\Gm)$. 
shows that the image of $(\fmlG')^{\vee}$ under $\hd$ is contained in $\fmlGd$, 
%$\im \hd|_{(\fmlG')^{\vee}} \subset \fmlGd$. 
which implies that $\hd: (M')^{\vee} \lra \Md$ is a homomorphism of 1-motives. 
\ePf

%\newpage 

\section{Universal rational maps} 
\label{sec:Univ-Fact-Prbl}

%Let $X$ be a variety  over a field $k$. 
The classical Albanese variety $\Alba{X}$ of a variety $X$ over a field $k$ 
(as in \cite[II, \S~3]{La}) is an abelian variety, 
defined together with the Albanese map $\alb: X \dra \Alba{X}$ 
by the following universal mapping property: 
for every rational map $\phe:X\dra A$ to an abelian variety $A$ 
there is a unique homomorphism $h:\Alba{X} \lra A$ 
such that $\phe = h \circ \alb$ up to translation by a constant $a \in A(k)$. 
Now we replace in this definition the category of abelian varieties 
by a subcategory $\Cat$ of the category 
of commutative algebraic groups. 
A result of Serre \cite[No.~6, Th\'eor\`eme~8, p.~10-14]{S2} 
says that if the category $\Cat$ contains the additive group $\Ga$ 
and $X$ is a variety of dimension $>0$, 
there does not exist an Albanese variety in $\Cat$ that is universal 
for all rational maps from $X$ to algebraic groups in $\Cat$. 
One is therefore led to restrict the class of considered rational maps. 
This motivates the concept of 
\emph{categories of rational maps from $X$ to 
commutative algebraic groups}, 
or more generally, 
\emph{categories of rational maps from $X$ to torsors  
under commutative algebraic groups} 
($\see$Definition \ref{CatMr}), 
and to ask for the existence of universal objects for such categories. 

For $k$ an algebraically closed field with $\chr(k) = 0$, 
in \cite[Section 2]{Ru1} a criterion is given, 
for which categories $\Mr$ 
of rational maps from a smooth proper variety $X$ over $k$ 
to algebraic groups 
there exists a universal object $\Alb_{\Mr}\pn{X}$, 
as well as an explicit construction of these universal objects 
via duality of 1-motives. 
In this section we prove similar results for categories of rational maps, 
defined over a perfect field, 
from a smooth proper variety $X$ to torsors for commutative algebraic groups. 
%Similar results are true for perfect base field of arbitrary characteristic as well, 
%as we will see in this section. 

\subsection{Relative Cartier divisors} 
\label{relDiv}

The construction of such universal objects as above involves  
the functor  
$\Divf_X: \Algk \lra \Ab$ 
of families of Cartier divisors, 
given by 
\[ \Divf_X\pn{R}=\left\{ 
   \begin{array}{c}
     \textrm{Cartier divisors $\sD$ on $X\times_{k}\Spec R$}\\
     \textrm{whose fibres $\sD_p$ define Cartier divisors 
              on $X \tms_k \{p\}$} \\
     \textrm{for all }\; p \in \Spec R 
   \end{array} 
\right\} 
\] 
for each $k$-algebra $R$, 
and for a homomorphism $h:R\lra S$ in $\Algk$ 
the induced homomorphism $\Divf_X\pn{h}:
\Divf_X\pn{R}\lra\Divf_X\pn{S}$ in $\Ab$ 
is the pull-back of Cartier divisors on $X\times_{k}\Spec R$
to those on $X\times_{k}\Spec S$. 
The elements of $\Divf_X\pn{R}$ are called 
\emph{relative Cartier divisors}. 
See \cite[No.~2.1]{Ru1} for more details on $\Divf_X$. 

 %This functor admits a natural transformation to the Picard functor, 
%which describes families of line bundles, 
%i.e.\ families of classes of Cartier divisors. 

We will be mainly concerned with the completion 
$\Divfc{X}: \Artk \lra \Ab$ of $\Divf_X$, 
which is given for every finite $k$-algebra $R$ by 
\[ \Divfc{X}\pn{R} = 
   \Gam\bigpn{X\tens R, \lrpn{\sK_X \tens_k R}^*/\lrpn{\sO_X \tens_k R}^*}  
   \laurin 
\] 

We will regard $\Divfc{X}$ as a subsheaf of $\Divf_X$, 
cf.\ Remark \ref{fmlG in Abk}.  

\bPrp 
\label{Divf formal group} 
$\Divfc{X}$ is a formal $k$-group. 
\ePrp 

\bPf 
According to \cite[I, No.~6]{D} or \cite[I, \S~4]{Fo} 
%Proposition \ref{formalScheme} 
it suffices to show that $\Divfc{X}$ is left-exact 
(i.e.\ commutes with finite projective limits). 
We are going to show that $\Divfc{X}$ is the composition of left-exact functors. 

Let $\Rfin$ %$\in \Artk$ 
be a finite $k$-algebra. 
%Each non unit $r \in \Rfin\setminus\Rfin^*$ is a zero divisor. 
%Thus we have 
%$\sK_{X\tens\Rfin/\Rfin} = \sK_{X\tens\Rfin} = \sK_X \tens_k \Rfin$ 
%and  $\sO_{X\tens\Rfin} = \sO_X \tens_k \Rfin$. 
$\Divfc{X}(\Rfin) = \Gam\big(X,\sQ(\Rfin)\big)$ 
is the abelian group of global sections of the sheaf 
$\sQ(\Rfin):=\lrpn{\pr_X}_*
  \bigpn{\left.\lrpn{\sK_X\tens\Rfin}^* \right/ \lrpn{\sO_X\tens\Rfin}^*}$, 
where $\pr_X:X\tens\Rfin \lra X$ is the projection. 
The global section functor $\Gam(X,\lul)$ is known to be left-exact. 
We show that the formal $k$-group functor  $\sQ:\Artk\lra\AbS{X}$ 
(with values in the category of sheaves of abelian groups over $X$) 
commutes with finite projective limits (hence is left-exact): 

Let $(\Rfin_i)$ be a projective system of local finite $k$-algebras, 
with homomorphisms $h_{ij}:\Rfin_j \lra \Rfin_i$ for $i<j$. 
We have projections $\pr_j:\varprojlim \Rfin_i \lra \Rfin_j$ for each $j$, 
which commute with the $h_{ij}$. 
Functoriality of $\sQ$ in $\Rfin\in\Artk$ induces homomorphisms 
$\sQ(h_{ij}):\sQ(\Rfin_j)\lra\sQ(\Rfin_i)$ and 
$\sQ(\pr_j):\sQ(\varprojlim \Rfin_i)\lra\sQ(\Rfin_j)$, which commute. 
The universal property of $\varprojlim\sQ(\Rfin_i)$ yields a unique 
homomorphism of sheaves 
$\sQ(\varprojlim \Rfin_i) \lra \varprojlim\sQ(\Rfin_i)$. 
A homomorphism of sheaves is an isomorphism 
if and only if it is an isomorphism on stalks. 
Therefore it remains to show that the stalks $\sQ_{\pnt}:\Artk\lra\Ab$ 
for $\pnt \in X$ are left-exact in $\Rfin\in\Artk$. 
We have 
$\sQ_{\pnt} = \left.\Gm(\sK_{X,\pnt}\tens_k\lul)\right/\Gm(\sO_{X,\pnt}\tens_k\lul)$. 
The tensor product over a field $\sA\tens_k\lul:\Artk$ $\lra\Algk$ 
is exact for any $k$-algebra $\sA$. 
Also the sheaf $\Gm:\Algk\lra\Ab$ is left-exact. 
Therefore the formal $k$-group functors $\Gm(\sK_{X,\pnt}\tens_k\lul)$ 
and $\Gm(\sO_{X,\pnt}\tens_k\lul)$ are formal $k$-groups. 
Since the category $\Gfk$ of formal $k$-groups is abelian 
($\seecite$\cite[$\textrm{VII}_{\textrm{B}}$, 2.4.2]{SGA3}), 
the quotient $\sQ_{\pnt}$ of these two formal $k$-groups is again 
a formal $k$-group. 
\ePf 

%\bPrp 
%\label{Divf(k), Lie(Divf)} 
%\[ \Lie\lrpn{\Divf_X} \,=\, 
%\Gamma\lrpn{X,\left.\sK_{X}\right/\sO_{X}} 
%\] 
%\ePrp 

\bDef 
\label{Support} 
Let $\Rfin$ be a finite $k$-algebra. 
If $D \in \bigpn{\Divfc{X}}_{\et}\pn{\Rfin}$, then $\Supp\pn{D}$ 
denotes the locus of zeroes and poles of local sections 
$\lrpn{f_{\alpha}}_{\alpha}$ of $\Gm\lrpn{\sK_{X}\tens\Rfin_{\red}}$ 
representing 
$D \in \Gam\bigpn{\left.\Gm\lrpn{\sK_{X}\tens\Rfin_{\red}}\right/
  \Gm\lrpn{\sO_{X}\tens\Rfin_{\red}}}$. \\ 
%$ = \bigpn{\Divfc{X}}_{\et}\pn{\Rfin}$. 
If $\delta \in \bigpn{\Divfc{X}}_{\inf}\pn{\Rfin}$, 
then $\Supp\pn{\delta}$ denotes the locus of poles of local sections 
$ \lrpn{g_{\alpha}}_{\alpha}$ of $\Upf_{\Rfin}(\sK_X)$ 
representing 
$\delta \in \Gamma\bigpn{\left.\Upf_{\Rfin}(\sK_X)\right/ 
  \Upf_{\Rfin}(\sO_X)} = \bigpn{\Divfc{X}}_{\inf}\pn{\Rfin}$, 
where $\Upf_{\Rfin} = \ker\bigpn{\WR{\Gm}{\Rfin} \lra \WR{\Gm}{\Rfin_{\red}}}$ 
is the unipotent part of $\Lin_{\Rfin}$ from No.~\ref{LinGroup_Ring}.) 
\eDef 

\bDef 
\label{Supp(F)} 
Let $\fmlG$ be a formal subgroup of $\Divf_X$. 
The support of $\fmlG$ is defined to be 
\[ \Supp\pn{\fmlG} = %\bigcup_{\Rfin\in\Artk} 
   \bigcup_{\substack{\Rfin\in\Artk \\ \sD\in\fmlG\pn{\Rfin}}}\Supp\pn{\sD} 
\] 
where we use the decomposition $\fmlG = \fmlG_{\et} \tms \fmlG_{\inf}$ 
and Definition \ref{Support}. 
%where $\Supp\pn{ D }$ is the support of a Cartier divisor on $X$ or of 
%a deformation of the zero divisor ($\see$Definition \ref{Support}). 
\eDef

Suppose now that $X$ is a geometrically irreducible 
smooth proper variety over a perfect field $k$. 
%an algebraically closed field. 
Then the Picard functor $\Picf_X$ is represented 
by a separated algebraic space $\Pic_X$, 
whose identity component $\Pic_X^0$ is a proper scheme over $k$ 
($\seecite$\cite[No.~8.4, Thm.~3]{BLR}). 
The underlying reduced scheme $\Picor{X}$ of $\Pic_X^0$ 
is an abelian variety, called the \emph{Picard variety} of $X$. 
The subfunctor of $\Picf_X$ that is represented by $\Picor{X}$ 
will be denoted by $\Picorf{X}$. 

There is a natural transformation 
\[ \cl:\Divf_X\lra\Picf_X \laurin 
\] 
We define $\Divor{X}$
%\[ \Divor{X}:\Algk\lra\Ab 
%\] 
to be the subfunctor of $\Divf_X$ given by 
%\[ \Divor{X}\pn{R}=\cl^{-1}\lrpn{\Picorf{X}\pn{R}} 
%\] 
%for each $k$-algebra $R$. 
\[ \Divor{X} = \Divf_X \tms_{\Picf_X} \Picorf{X} \laurin 
\]

\subsection{Categories of rational maps to torsors} 
\label{sub:Categories-of-Rational} 

Let $X$ be a smooth proper variety over a perfect field $k$. 
%an algebraically closed field $k$ of arbitrary characteristic. 
Let $\clfld$ be an algebraic closure of $k$. 
In this note, algebraic groups and formal groups 
are commutative by definition 
($\see$No.~\ref{AlgGroups-FmlGroups}), 
and torsors are always torsors for commutative algebraic groups. 
%Algebraic groups are always assumed to be smooth and connected, 
%unless stated otherwise. 

\subsubsection{Induced transformation} 
\label{sec:indTrafo}

Let $G$ be a smooth connected algebraic group, 
and let $0\to L\to G\to A \to 0$ be the canonical decomposition of $G$, 
where $A$ is an abelian variety and $L$ an affine smooth connected algebraic group 
($\see$Theorem of Chevalley). %(\ref{Thm Chevalley}). 
Write $L = \Trs \tms_k \Upt$ where $\Trs$ is a torus 
and $\Upt$ is unipotent 
($\seecite$\cite[XVII, 7.2.1]{SGA3}). 
If $k$ is algebraically closed, 
$\Trs \cong (\Gm)^t$ for some $t\in\Nat$. 
If $k$ is of characteristic $0$, one has 
$\Upt \cong (\Ga)^s$ for some $s\in\Nat$ 
($\see$\cite[IV, \S~2, 4.2]{DG}). 
%Fix such isomorphism. 
If $k$ is of characteristic $p>0$, 
the unipotent group $\Upt$ is embedded into a finite direct sum 
$\lrpn{\Witt_{r}}^s$ of Witt vector groups for some $r,s\in\Nat$ 
($\see$\cite[V, \S~1, 2.5]{DG}). 
%Fix such an embedding.

%Let $K$ be the function field od $X$. 
Since $\H^1_{\fppf}\bigpn{\Spec(\sO_{X,\pnt}),\Gm} = 0$ 
and $\H^1_{\fppf}\bigpn{\Spec(\sO_{X,\pnt}),\Upt} = 0$ 
for any point $\pnt$ of $X$, we have exact sequences 
\begin{align*} 
 0 \lra L\lrpn{\sK_{X,\pnt}} \lra G\lrpn{\sK_{X,\pnt}} \lra 
   A\lrpn{\sK_{X,\pnt}} \lra 0 \phantom{\laurin} \\ 
 0 \lra L\lrpn{\sO_{X,\pnt}} \lra G\lrpn{\sO_{X,\pnt}} \lra 
   A\lrpn{\sO_{X,\pnt}} \lra 0 \laurin 
\end{align*} 
%Since $A$ is proper and $\sO_{X,\pnt}$ is a discrete valuation ring 
%for every point $\pnt$ of codimension 1 in $X$, 
Since a rational map to an abelian variety 
is defined at every smooth point 
($\seecite$\cite[II, \S~1, Thm.~2]{La}), 
we have $A\lrpn{\sK_{X,\pnt}} = A\lrpn{\sO_{X,\pnt}}$ 
for every point $\pnt$ of $X$. 
Hence the canonical map 
\[ L\lrpn{\sK_{X,\pnt}}/L\lrpn{\sO_{X,\pnt}} \lra 
   G\lrpn{\sK_{X,\pnt}}/G\lrpn{\sO_{X,\pnt}} 
\] 
is bijective. 
By Cartier-duality, we have a pairing 
\[ \lrpair{\llul,\lull} : \Ld \tms \Gam\bigpn{L(\sK_X)/L(\sO_X)} 
   \lra \Gam\bigpn{\Gm(\ul{\sK_X})/\Gm(\ul{\sO_X})} 
\] 
where $\ul{\sK_X} := \sK_X \tens \blank$ and $\ul{\sO_X} := \sO_X \tens \blank$. 
  
\bDef 
\label{induced Trafo}
Let $\phe: X \dra G$ be a rational map to a smooth connected 
algebraic group $G$, %$G\in\Ext\lrpn{A,L}$. 
let $L$ be the affine part of $G$. 
Then $\trafo_{\phe}: \Ld \lra \Divfc{X}$ denotes 
the induced transformation given by 
$\lrpair{\llul,\lin_{\phe}}$, 
where $\lin_{\phe}$ is the image of $\phe\in G\lrpn{\sK_X}$ in 
$\Gam\bigpn{G(\sK_X)/G(\sO_X)} \iso 
  \Gam\bigpn{L(\sK_X)/L(\sO_X)}$. 
By construction, $\trafo_{\phe}$ is a homomorphism of formal $k$-group functors. 
\eDef 

\bLem 
\label{im(tau)} 
Let $G$ be a smooth connected algebraic group, 
let $L$ be the affine part of $G$. 
Let $\phe:X\dra G$ be a rational map. 
Let $\trafo_{\phe}:\Ld\lra\Divfc{X}$ be the induced transformation. 
Then $\im\pn{\trafo_{\phe}}$ is a dual-algebraic formal group. 
\eLem 

\bPf 
$\Divfc{X}$ is a formal group by Proposition \ref{Divf formal group}, 
and $\Gfk$ is a full subcategory of $\Fctr(\Artk,\Ab)$. 
Therefore $\trafo_{\phe}:\Ld\lra\Divfc{X}$ 
is a homomorphism of formal groups. 
Since $\Gfk$ is an abelian category,  kernel and image of 
the homomorphism $\trafo_{\phe}$
are formal groups. 
Since $L$ is algebraic, $\Ld$ is dual-algebraic 
and hence $\im\pn{\trafo_{\phe}}$, as a quotient of $\Ld$, 
is dual-algebraic ($\see$Lemma \ref{subFmlGroup_algebraic}). 
\ePf 

\bLem 
\label{Ld->Picor} 
Let $G\in\Extabk(A,L)$ be a smooth connected algebraic group. 
%with affine part $L$. 
Let $\phe:X\dra G$ a rational map. 
Let $\trafo_{\phe}:\Ld\lra\Divfc{X}$ be the induced transformation. 
Then $\im\pn{\trafo_{\phe}}$ is contained in the completion of $\Divor{X}$. 
\eLem 

\bPf 
As $A$ is an abelian variety, 
the composition $X \overset{\phe} \dra G \overset{\rho} \lra A$ 
extends to a morphism $\ol{\phe}:X \lra A$. 
The description of the induced transformation $\trafo_{\phe}$ 
in terms of local sections into principal fibre bundles 
as given in \cite[No.~2.2]{Ru1} shows: 
the composition 
\bDpl 
\Ld \overset{\trafo_{\phe}}\lra \Divf_X \overset{\cl}\lra \Picf_X 
\eDpl 
is given by 
\;$\lam \lmt \lam_*G_X$, \;
where $\lam_*G$ is the push-out of $G\in\Extabk(A,L)$ 
via $\lam \in \Ld(R) = \Homabk(L,\Lin_R)$, 
and $G_X = G \tms_A X$ is the fibre-product of $G$ and $X$ over $A$. 
Hence it comes down to show that 
for each $R\in\Artk$, each $\lam\in\Ld(R)$ 
the $\Lin_R$-bundle $\lam_* G_X$ yields an element of $\Picor{X}(R)$. 

The universal mapping property of the classical Albanese $\Alba{X}$ 
yields that $\ol{\phe}$ factors through $\Alba{X}$. 
Hence the pull-back $G_X=G\tms_A X$ over $X$ %of $G$ to $X$ 
is a pull-back of $G_{\Alb} = G\tms_A \Alba{X}$ over $\Alba{X}$. 
%\[ \xymatrix{ 
%   G_X \ar[dd] \ar[rr] && 
%   G_{\Alb} \ar[dd] \ar[rr]  &&  G \ar[dd]^{\rho} \\  
%   && && \\  
%   X \ar@/_1.4pc/[uu]_{\phe_{X}} \ar@/_0.7pc/[uurr]_(0.65){\phe_{\Alb}} \ar[rr] 
%   \ar@/_/[uurrrr]_(0.62){\phe}  &&  \Alb(X) \ar[rr] && A } 
%\] 
%%\[ \xymatrix{ 
%%   L \ar[dd] \ar@{=}[rr] && L \ar[dd] \ar@{=}[rr] && L \ar[dd] \\  
%%   && && \\  
%%   G_X \ar@<-0.7ex>[dd] \ar[rr] && 
%%   G_{\Alb} \ar[dd] \ar[rr]  &&  G \ar[dd]^{\rho} \\  
%%   && && \\  
%%   X \ar@{-->}@<-0.7ex>[uu]_{\phe_{X}} \ar@{-->}[uurr]_{\phe_{\Alb}} \ar[rr] 
%%   \ar@{-->}[uurrrr]_(0.45){\phe}  &&  \Alb(X) \ar[rr] && A } 
%%\] 
Then for each $\lam\in\Ld(R)$ the $\Lin_R$-bundle $\lam_* G_{\Alb}$ over 
$\Alba{X}$ is an element of $\Extabk(\Alba{X},\Lin_R)$, 
hence gives an element of $\Pic^0_{\Alba{X}}(R)$. 
Since $\Alba{X} = \big(\Picor{X}\big)^{\vee}$ 
is the dual abelian variety of $\Picor{X}$, we have an isomorphism 
$\Pic^0_{\Alba{X}} \overset{\sim}\lra \Picor{X}$, 
$P \lmt P_X = P\tms_{\Alba{X}}X$. 
%\begin{eqnarray*}
%\Pic^0_{\Alba{X}} & \overset{\sim}\lra & \Picor{X} \\
%P & \lmt & P_{X} = P \tms_{\Alba{X}} X 
%\end{eqnarray*} 
As $\lam_* G_X = \lam_* G_{\Alb} \tms_{\Alba{X}} X$, it holds 
$\lam_* G_X \in \Picor{X}(R)$. 
\ePf 

\bLem 
\label{induced ratl map}
Let $L$ be an affine algebraic group 
and $\trafo: \Ld \lra \Divorc{X}$ a homomorphism of formal groups. 
Let $G \in \Extabk\bigpn{\Alba{X},L}$ be the extension corresponding to 
\;$\cl \circ \,\trafo: \Ld \lra \Divor{X} \lra \Picorf{X}$\; 
under the bijection $\Pop$ from Theorem~\ref{Ext(A,L)=Hom(Ld,Ad)}. 
There is a rational map $\phe: X \dra G$ whose 
induced transformation is $\trafo$, 
and $\phe$ is determined uniquely up to translation by a constant $g \in G(k)$. 
\eLem 

\bPf 
By Lemma \ref{L subset L_S} we may choose an embedding 
$\lam: L \inj \Lin_S$ for some finite ring $S \in \Artk$. 
Let $\sD \in \Divor{X}(S)$ be the image of 
$\id_{\Lin_S} \in \Homabk\pn{\Lin_S,\Lin_S} = \Lin_S^{\vee}(S)$ 
under the composition 
$ \trafo \circ \lam^{\vee}: \Lin_S^{\vee} \lra \Ld \lra \Divf_X 
\laurin 
$ 
%The proof of Theorem \ref{Ext(A,L)=Hom(Ld,Ad)} 
Remark \ref{Pop=conHom} shows that 
$\sO_{X \tens S}(\sD) \in \Picor{X}(S)$ is the line bundle 
corresponding to %$G$ under the map 
$G \in \Extabk\pn{\Alb(X),L}$ under the map 
\[ \Extabk\pn{\Alb(X),L} \lra \Extabk\pn{\Alb(X),\Lin_S} 
= \Pic_{\Alb(X)}^0(S) \iso \Picor{X}(S) 
\laurin 
\] 
Let $\Gp_{\Lin_S}(\sD)$ be the image of $G$ in $\Extabk\pn{\Alb(X),\Lin_S}$, 
%under $\Extabk\pn{\Alb(X),L} \lra \Extabk\pn{\Alb(X),\Lin_S}$, 
then the fibre product $\Pb_{\Lin_S}(\sD) := \Gp_{\Lin_S}(\sD) \tms_{\Alba{X}} X$ 
is the $\Lin_S$-bundle on $X$ associated to $\sO_{X \tens S}(\sD)$, 
and $G$ resp.\ $P := G \tms_{\Alba{X}} X$ 
are reductions of the $\Lin_S$-bundles 
$\Gp_{\Lin_S}(\sD)$ resp.\ $\Pb_{\Lin_S}(\sD)$ 
to the fibre $L$. 
The canonical 1-section of $\sO_{X\tens\Sfin}(\sD)$ 
yields a section $X \dra P$, 
and composition with $P \lra G$ 
yields the desired rational map $\phe: X \dra G$, 
which by construction satisfies $\trafo_{\phe} = \trafo$. 
Then %element 
$\lin_{\phe} \in 
\Gam\bigpn{G(\sK_X)/G(\sO_X)} \cong \Gam\bigpn{L(\sK_X)/L(\sO_X)} 
\subset \Gam\bigpn{\Lin_S(\sK_X)/\Lin_S(\sO_X)} 
= \Gam\bigpn{\Gm(\sK_X \tens S) / \Gm(\sO_X \tens S)}$ 
corresponding to $\sD$ 
is uniquely determined by $\trafo$. 
The rational map $\phe \in G(\sK_X)$, as a lift of $\ell_{\phe}$, 
is determined up to a constant $g \in G(k) = \Gam\bigpn{G(\sO_X)}$, 
according to the exact sequence 
\bDpl 
0 \lra \Gam\bigpn{G(\sO_X)} \lra \Gam\bigpn{G(\sK_X)} 
   \lra \Gam\bigpn{G(\sK_X)/G(\sO_X)}\laurin 
\eDpl 
%\hfill 
\ePf

\subsubsection{Definition of a category of rational maps} 
\label{DefCatRatMaps}

\bDef 
\label{CatMr} 
A \emph{category} $\Mr$ \emph{of rational maps from} $X$
\emph{to torsors} is a category satisfying the following conditions: 
The objects of $\Mr$ are rational maps $\phe: X \dra P$, 
where $P$ is a torsor for a smooth connected algebraic group. 
%and $0 \in \im\phe$. 
The morphisms of $\Mr$ between two objects $\phe: X \dra P$
and $\psi: X \dra Q$ are given by the set of 
those homomorphisms of torsors 
%(homomorphisms of algebraic groups composed with a translation) 
$h:P \lra Q$ 
such that $h \circ \phe = \psi$. 
%i.e.\ the following diagram commutes: 
%\[  \xymatrix{  & X \ar%@{-->}
%[dl]_{\phe} \ar%@{-->}
%[dr]^{\psi} &  \\  
%                G \ar[rr]^{h} & & H \laurin 
%             } 
%\] 
\eDef 

\bRmk 
\label{EquCatMr} 
Let $\phe:X\dra P$ and $\psi:X\dra Q$ be two rational maps 
from $X$ to torsors. Then Definition \ref{CatMr} 
implies that for any category $\Mr$ of rational maps 
from $X$ to torsors containing $\phe$ and $\psi$ as objects 
the set of morphisms $\Hom_{\Mr}(\phe,\psi)$ is the same.  
Therefore two categories $\Mr$ and $\Mr'$ of rational maps from $X$ to 
torsors are equivalent if every object of $\Mr$ is isomorphic to 
an object of $\Mr'$. 
\eRmk 

\bRmk 
\label{k-ratl pt of torsor}
If a $\fld$-torsor $P$ for an algebraic $\fld$-group $G$ 
admits a $\fld$-rational point, 
then $P$ may be identified with $G$. 
%Suppose $\anyfld$ is an extension of finite type of $\fld$ in $\unidom$ 
%and let $\clanyfld$ be the algebraic closure (in $\unidom$). 
Thus for a rational map $\phe: \Xo \dra P$ 
it makes sense to consider the base changed map 
$\phe\tens_{\fld}\clfld: \Xo\tens_{\fld}\clfld \dra P\tens_{\fld}\clfld$ 
as a rational map from $\Xo\tens_{\fld}\clfld = \Xb$ 
to an algebraic $\clfld$-group 
$P\tens_{\fld}\clfld \cong G\tens_{\fld}\clfld$. 
\eRmk 

\bDef 
\label{Mav} The category of rational maps from $X$ to abelian varieties 
is denoted by $\Mav$. 
%Since a rational map to an abelian variety 
%is defined at every smooth point ($\seecite$\cite[II, \S~1, Thm.~2]{La}) 
%and $X$ is smooth, 
%$\Mav$ is equal to the category of morphisms from $X$ to abelian varieties. 
\eDef 

\bRmk 
The objects of $\Mav$ are in fact \emph{morphisms} from $X$ to abelian varieties, 
since a rational map from a smooth variety $X$ 
to an abelian variety $A$ extends to a morphism from $X$ to $A$ 
($\seecite$\cite[II, \S~1, Thm.~2]{La}). 
\eRmk 

\bDef 
\label{Mr_F} 
Let $\fmlG$ be a dual-algebraic formal $k$-subgroup of $\Divf_{X}$. 
If $k$ is algebraically closed, 
then $\MrF$ denotes the category of all 
those rational maps $\phe:X\dra G$ from $X$ to algebraic $k$-groups 
for which the image of the induced transformation 
$\trafo_{\phe}:\Ld\lra\Divf_{X}$
($\see$Definition \ref{induced Trafo}) lies in $\fmlG$, 
i.e.\ which induce a homomorphism of formal groups $\Ld\lra\fmlG$,
where $L$ is the affine part of $G$. 
%\[ \MrF = \{ \phe:X\dra G \;|\; \im \trafo_{\phe} \subset \fmlG \}
%\]
For general $k$, $\MrF$ denotes the category of all 
those rational maps $\phe:X\dra P$ from $X$ to $k$-torsors 
for which the base changed map $\phe \tens \clfld$ 
is an object of $\Mr_{\fmlG \tens \clfld}$. 
(Here we use Rmk.~\ref{k-ratl pt of torsor}.) 
%i.e.\ that a torsor over $\clfld$ 
%can be identified with the algebraic group acting on it.) 
\eDef

\subsection{Universal objects} 
\label{sub:Universal-Objects}

Let $X$ be a smooth proper variety over a perfect field $k$. 
%which in this No.\ is assumed to be algebraically closed. 
Algebraic groups are always assumed to be smooth and connected, 
and torsors are those for smooth connected algebraic groups, 
unless stated otherwise.

\subsubsection{Existence and construction} 
\label{subsub:Exist+Construct} 

\bDef 
\label{univObj} 
Let $\Mr$ be a category of rational maps from $X$ to torsors. 
%algebraic groups. 
Then \,$\lrpn{u:X\dra\sU}\in\Mr$\, is called
a \emph{universal object for} $\Mr$ if it admits 
the universal mapping property in $\Mr$: 
For all \,$\lrpn{\phe:X\dra P}\in\Mr$\, there is a unique 
homomorphism of torsors \,$h:\sU\lra P$\, such that \,$\phe=h\circ u$. 
\eDef 

\bRmk 
\label{univObj_unique}
Universal objects are uniquely determined up to (unique) isomorphism. 
\eRmk 

Now assume that the base field $k$ is algebraically closed. 
(Arbitrary perfect base field is considered from 
No.~\ref{descentBaseField} on.) 
%and No.~\ref{Functoriality}.) 
In this case we may identify a torsor with the algebraic group acting on it 
($\see$Rmk.~\ref{k-ratl pt of torsor}), 
and a homomorphism of torsors becomes a homomorphism of algebraic groups 
composed with a translation (which is an isomorphism of torsors). 

For the category $\Mav$ of morphisms from $X$ to abelian varieties
($\see$Def.~\ref{Mav}) there exists a universal object,
the \emph{Albanese mapping} to the \emph{Albanese variety,} denoted
by $\alb:X\lra\Alba{X}$. 
This is a classical result ($\seecite$\cite{La}, \cite{Ms}, \cite{S2}). 
The Albanese variety $\Alba{X}$ is an abelian variety, 
dual to the Picard variety $\Picor{X}$. 

\vspace{\vs} 
%\newpage 

In the following 
we consider categories $\Mr$ of rational maps from
$X$ to algebraic groups satisfying the following conditions: 

\begin{tabular}{rl}
$\lrpn{ \diamondsuit \; 1 }$ & 
$\Mr$ contains the category $\Mav$. \\ 
$\lrpn{ \diamondsuit \; 2 }$ & 
If $\lrpn{\phe:X\dra G}\in\Mr$ 
and $h: G \lra H$ is a homomorphism of \\ 
 & torsors for smooth connected algebraic groups, then $h\circ\phe \in \Mr$. 
\end{tabular} 

\newpage 

\bThm 
\label{Exist univObj}  
Let $\Mr$ be a category of rational maps from $X$ to algebraic groups 
that satisfies $\lrpn{\diamondsuit \; 1,2 }$. 
Then the following conditions are equivalent: 

\begin{tabular}{rl}
  (i) & For $\Mr$ there exists a universal object $\lrpn{u:X\dra\sU}\in\Mr$\laurin \\ 
 (ii) & There is a dual-algebraic formal subgroup $\fmlG$ of $\Divor{X}$ \\ 
       & such that $\Mr$ is equivalent to $\MrF$\laurin \\ 
(iii) & The formal group $\fmlG_{\Mr} \subset \Divor{X}$ 
          generated by $\bigcup_{\phe \in \Mr} \im\pn{\trafo_{\phe}}$ \\ 
      & is dual-algebraic and $\Mr = \Mrt{\fmlG_{\Mr}}{X}$\laurin 
\end{tabular} 

\noindent 
Here $\MrF$ is the category of rational maps that induce a 
homomorphism of formal groups to $\fmlG$ 
($\see$Definition \ref{Mr_F}). 
\eThm 

\bPf  
(ii)$\Lra$(i) 
Assume that $\Mr$ is equivalent to $\MrF$, 
where $\fmlG$ is a dual-algebraic formal group in $\Divor{X}$.
The first step is the construction of an algebraic group $\sU$ and
a rational map $u:X\dra\sU$. In a second step the universality
of $u:X\dra\sU$ for $\MrF$ will be shown. 

\textbf{Step 1:} Construction of $u:X\dra\sU \laurin$ \\
$X$ is a smooth proper variety over $k$, thus the functor $\Picf_X^0$
is represented by an algebraic group $\Pic_{X}^{0}$ 
whose underlying reduced scheme $\Picor{X}$, 
the Picard variety of $X$, is an abelian variety. 
The class map $\Divf_X\lra\Picf_X$
induces a homomorphism  $\fmlG\lra\Picor{X}$. 

We obtain a 1-motive $M=\big[\fmlG\lra\Picor{X}\big]$. 
%Let $M^{\vee}$ be the dual 1-motive of $M$. 
Since $\Picor{X}$ is an abelian variety, 
the dual 1-motive of $M$ is of the form $\Md = \lrbt{0 \ra G}$, 
where $G$ is a smooth connected algebraic group. 
Then define $\sU$ to be this algebraic group. 
%i.e.\ $\left[0\ra\sU\right]$ is the dual 1-motive of $\big[\fmlG\lra\Picor{X}\big]$.
The canonical decomposition $0\ra\sL\ra\sU\ra\sA\ra0$
is the extension of $\sA=\Alba{X}=\big(\Picor{X}\big)^{\vee}$ 
by $\sL=\fmlG^{\vee}$ 
induced by the homomorphism $\fmlG\lra\Picor{X}$ 
($\see$Theorem \ref{Ext(A,L)=Hom(Ld,Ad)}). 
%where $\sL=\fmlG^{\vee}$ is the Cartier-dual of $\fmlG$
%and $\sA=\Picor{X}^{\vee}$ is the dual abelian variety of $\Picor{X}$, 
%which is $\Alba{X}$. 
%We have an isomorphism 
%\begin{eqnarray*}
%\Ext\lrpn{\sA,\Gm}\:\simeq\:
%\Pic_{\sA}^{0} & \overset{\sim}\lra & \Picor{X} \\
%P & \lmt & P_X = P\tms_{\sA}X 
%\end{eqnarray*} 

%By Lemma \ref{L subset L_S} there exists a finite $k$-algebra $\Sfin$ 
%and an injective homomorphism of affine algebraic groups 
%$\lam: \sL \lra \Lin_{\Sfin}$. 
%Assume first that $\sL = \Lin_{\Sfin}$. 
%The homomorphism $\fmlG \lra \Picor{X}$ on $\Sfin$-valued points 
%$\Lin_{\Sfin}^{\vee}(\Sfin) \lra \sA^{\vee}(\Sfin)$  
%has values in $\Extabk(\sA,\Lin_{\Sfin})$, 
%according to Proposition \ref{Xeil}. 
%The proof of Theorem \ref{Ext(A,L)=Hom(Ld,Ad)} shows that 
%$\sU \in \Extabk(\sA,\Lin_{\Sfin})$ is the image of 
%$\iota := \id_{\Lin_{\Sfin}} \in \Homabk(\Lin_{\Sfin},\Lin_{\Sfin}) = 
%  \Lin_{\Sfin}^{\vee}(\Sfin)$. 
%%under the homomorphism $\Lin_{\Sfin}^{\vee} \lra \sA^{\vee}$. 
%The case $\sL \subset \Lin_{\Sfin}$ is achieved by a descent argument. 
%The composition of the pull-back $\lam^{\vee}: \Lin_{\Sfin}^{\vee} \lra \sL^{\vee}$ 
%with $\sL^{\vee} = \fmlG \lra \Picor{X} = \sA^{\vee}$ 
%yields a homomorphism $\phi: \Lin_{\Sfin}^{\vee} \lra \sA^{\vee}$. 
%As $\phi$ factors through $\sL^{\vee}$, the image $\sU = \phi\lrpn{\iota}$ 
%of $\iota \in \Lin_{\Sfin}^{\vee}(\Sfin)$ lies actually in $\Extabk(\sA,\sL)$, 
%cf.\ Step 2 in the proof of Theorem \ref{Ext(A,L)=Hom(Ld,Ad)}. 
Define the rational map $u:X\dra\sU$ by the condition that 
the induced transformation $\trafo_u: \fmlG \lra \Divor{X}$ 
from Definition \ref{induced Trafo} is the inclusion. 
According to Lemma \ref{induced ratl map}, 
$u: X \dra \sU$ is determined up to a constant $c \in \sU(k)$. 

Note that $u:X\dra\sU$ generates $\sU$: 
Let $H$  %:= \lrpair{\im \phe} \subset \sU$ 
be the subgroup generated by the image of $u$, 
let $\Lam \subset \sL$ be the affine part of $H$. 
As $u: X \lra \sU$ factors through $H$, 
the induced transformation $\trafo_{u}: \sL^{\vee} \lra \Divfc{X}$ 
factors through the quotient $\sL^{\vee} \lsur \Lamd$. 
As $\trafo_{u}$ is injective, this yields $\Lamd \cong \sL^{\vee}$, 
hence $\Lam \cong \sL$. 
As the composition $X \overset{u} \dra \sU \lra \sA$ generates $\sA$, 
the abelian quotient of $H$ is $\sA$. 
These two conditions imply that $H \cong \sU$ by the Five Lemma. 

\textbf{Step 2:} Universality of $u:X\dra\sU \laurin$ \\
Let $\phe:X\dra G$ be a rational map to a smooth connected algebraic group 
$G$ with canonical decomposition $0\ra L\ra G\overset{\rho}\ra A\ra0$, 
inducing a homomorphism of formal groups 
$\trafo_{\phe}: \Ld \lra \fmlG \subset \Divor{X}$,
$\lam \lmt \lrpair{\lam,\lin_{\phe}}$ %for $\lam\in\Ld(R)$ 
($\see$Definition \ref{induced Trafo}). 
Let $l := (\trafo_{\phe})^{\vee} : \sL \lra L$ be the dual homomorphism 
of affine groups. 
The composition $X\overset{\phe}\dra G\overset{\rho}\lra A$
extends to a morphism from $X$ to an abelian variety. Translating 
$\phe$ by a constant $g\in G(k)$, if necessary, we may assume that 
$\rho\circ\phe$ factors through $\sA=\Alba{X}$. 
%\[  \xymatrix{  X \ar[dr]_{\alb} \ar[rr]^{\rho\,\circ\,\phe} & & A \\  
%                &  \Alb(X) \laurin \ar[ur] & \\ } \]
We are going to show that we have a commutative diagram as follows: 
\[ \xymatrix{ 
   & \sL \ar[d] \ar[r]^{l} & L \ar[d] \ar@{=}[r] 
                                            & L \ar[d] \ar@{=}[r] & L \ar[d] \\ 
   & \sU \ar[d] \ar[r]^h & l_*\,\sU \ar[d] \ar[r]^{\sim} 
                              & G_{\sA} \ar[d] \ar[r] & G \ar[d]^{\rho} \\ 
   X \ar[r] \ar@{-->}[urrrr]_(.55){\phe} \ar@{-->}[urrr]^(.5){\phe_{\sA}} 
               \ar@{-->}[ur]^u & \sA  \ar@{=}[r] & \sA  \ar@{=}[r] & \sA \ar[r] & A 
            } 
\]
%\[ \xymatrix{ 
%   && \sL \ar[dd] \ar[rr]^{l} && L \ar[dd] \ar@{=}[rr] 
%                                            && L \ar[dd] \ar@{=}[rr] && L \ar[dd] \\ 
%   &&  &&  &&  &&  \\                                          
%   && \sU \ar[dd] \ar[rr]^h && l_*\,\sU \ar[dd] \ar[rr]^{\sim} 
%                              && G_{\sA} \ar[dd] \ar[rr] && G \ar[dd]^{\rho} \\ 
%   &&  &&  &&  &&  \\                                          
%   X \ar[rr] \ar[uurrrrrrrr]_(.55){\phe} \ar[uurrrrrr]^(.5){\phe_{\sA}} 
%               \ar[uurr]^u && \sA  \ar@{=}[rr] && \sA  \ar@{=}[rr] && \sA \ar[rr] && A 
%            } 
%\]
where $G_{\sA}=G\times_{A}\sA$ is the fibre product, 
$l_*\,\sU = \sU \amalg_{\sL} L$ the amalgamated sum 
and $h: \sU \lra l_*\,\sU$ the map obtained from the amalgamated sum. 

\noindent 
If $\ol{l}: \Gam\bigpn{\sL(\sK_X)/\sL(\sO_X)} \lra \Gam\bigpn{L(\sK_X)/L(\sO_X)}$ 
denotes the map induced by $l: \sL \lra L$, 
then $\lin_{h \circ u} = \ol{l}(\lin_u)$. 
This yields 
\begin{eqnarray*}
\trafo_{h \circ u} 
& = & \lrpair{\llul,\lin_{h \circ u}} \; = \; \lrpair{\llul,\ol{l}(\lin_u)} 
\; = \; \lrpair{\llul \circ l,\lin_u} \\ 
& = & \trafo_u \circ \ld \; = \; \trafo_{\phe} 
\end{eqnarray*}
since $\trafo_u: \fmlG \lra \Divor{X}$ is the inclusion by construction of $u$. \\ 
This implies that $l_*\,\sU_X$ and $G_X$ are isomorphic $L$-bundles over $X$. 
Then $l_*\,\sU$ and $G_{\sA}$ are isomorphic as extensions of $\sA$ by $L$, 
using the isomorphism $\Pic_X^0 \iso \Pic_{\sA}^0$. 
Thus $\trafo_{h \circ u} = \trafo_{\phe}$ shows that 
$h \circ u$ and $\phe_{\sA}$ coincide up to translation. 
As $u: X \lra \sU$ generates $\sU$, each $h':\sU \lra G_{\sA}$ fulfilling 
$h' \circ u = \phe_{\sA}$ coincides with $h$. Hence $h$ is unique. 

(i)$\Lra$(iii) 
Assume that $u:X\dra\sU$ is universal for $\Mr$. 
Let $0\ra\sL\ra\sU\ra\sA\ra0$ be the canonical decomposition of $\sU$, 
and let $\fmlG$ be the image of the induced transformation 
$\trafo_u : \sL^{\vee} \lra \Divor{X}$.
For $\lam\in\sL^{\vee}(R)$  the uniqueness of the homomorphism 
$h_{\lam}:\sU\lra\lam_{*}\,\sU$ fulfilling $u_{\lam}=h_{\lam}\circ u$ 
implies that the rational maps $u_{\lam}:X\dra\lam_{*}\,\sU$ are 
non-isomorphic to each other for distinct $\lam \in \sL^{\vee}(R)$. 
Hence $\dv_{R}\lrpn{u_{X,\nu}} \neq \dv_{R}\lrpn{u_{X,\lam}}$ 
for $\nu \neq \lam \in \sL^{\vee}(R)$. 
Therefore $\sL^{\vee}\lra\fmlG$ is injective, hence an isomorphism. 

Let $\phe:X\dra G$ be an object of $\Mr$ and $0\ra L\ra G\ra A\ra0$
be the canonical decomposition of $G$. 
Translating $\phe$ by a constant $g\in G(k)$, if necessary, 
we may assume that $\phe:X\dra G$ 
factors through a unique homomorphism $h:\sU\lra G$.
The restriction of $h$ to $\sL$ gives a homomorphism of affine groups
$l:\sL\lra L$. 
Then the dual homomorphism $\ld:\Ld\lra\fmlG$ 
yields a factorization of $\Ld\lra\Divor{X}$ through $\fmlG$. 
The properties $\lrpn{\diamondsuit \; 1,2 }$ 
and the existence of a universal object 
guarantee that $\Mr$ contains all rational maps 
that induce a transformation to $\fmlG$, 
hence $\Mr$ is equal to $\MrF$.

(iii)$\Lra$(ii) 
is evident. 
\ePf 

\bNot  
%The universal object for a category $\Mr$ of rational maps
%from $X$ to algebraic groups, 
%if it exists, is denoted by 
%$\alb_{\Mr}:X\dra\Alb_{\Mr}\pn{X}$. \\
If $\fmlG$ is a dual-algebraic formal subgroup of $\Divor{X}$, 
then the universal object for $\MrF$ is denoted by 
$\alb_{\fmlG}:X\dra\AlbF{X}$. 
%\\
%For $\fmlG=0$ the universal object for $\Mr_{0}$ is usually simply
%denoted by $\alb:X\lra\Alba{X}$. 
\eNot  

\bRmk 
\label{Alb_constr} 
By construction, $\AlbF{X}$ is generated by $X$. 
Since $X$ is reduced, $\AlbF{X}$ is reduced as well, thus smooth. 
In the proof of Thm.\ \ref{Exist univObj}
we have seen that $\AlbF{X}$ is an extension of the
abelian variety $\Alba{X}$ by the affine group $\fmlG^{\vee}$. 
More precisely, $\bigbt{0\lra\AlbF{X}}$
is the dual 1-motive of $\bigbt{\fmlG\lra\Picor{X}}$. 
The rational map 
$\bigpn{\alb_{\fmlG}:X\dra\AlbF{X}} \in\MrF$
is characterized by the fact that the transformation 
$\trafo_{\alb_{\fmlG}}:\Ld\lra\Divor{X}$ 
is the identity $\fmlG\overset{\id}\lra\fmlG$. 
\eRmk

\subsubsection{Descent of the base field} 
\label{descentBaseField}

Let $\fld$ be a perfect field. 
%Let $\unidom$ be a universal domain, 
%i.e.\ $\unidom$ is an algebraically closed extension of $\fld$ 
%of infinite transcendence degree over $\fld$. 
Let $\clfld$ be an algebraic closure of $\fld$. 
Let $\Xo$ be a smooth proper variety defined over $\fld$, 
write $\Xb = \Xo \tens_{\fld} \clfld$. 
Let $\fmlG$ be a formal $\fld$-subgroup of $\Divf^0_{X}$, 
write $\fmlGb = \fmlG \tensc \clfld$. 
%and suppose that $\fmlG$ is defined over $\fld$: 
%
%\bDef 
%\label{fmlG/k_0} 
%A formal group $\fmlG \subset \Divf^0_{\Xb}$ is 
%\emph{defined over} $\fld$, 
%if for each finitely generated $\clfld$-algebra $R$ 
%the following condition is satisfied: \\ 
%If $\sD\in\fmlG(R)$ 
%then for each $\gal\in\Gal\pn{\clfld/\fld}$ 
%it holds $\sD^{\gal}\in\fmlG(R)$, 
%where $\sD^{\gal}$ is the conjugate of $\sD$ by means of $\gal$. 
%\eDef 

The wish is to show that the universal object 
$\albFb:\Xb\dra\AlbFb{\Xb}$ for the category $\MrFb$ 
can be defined over $\fld$. 
This will be accomplished by a Galois descent, 
as described in \cite[V, \S~4]{S}. 
%A short sketch of this procedure is the following: 
%Let $\erwfld$ be a (finite) Galois extension of $\fld$. 
%Let $\Ve$ be a $\erwfld$-variety living in a category $\Cat$ 
%of $\erwfld$-varieties (e.g.\ algebraic $\erwfld$-groups or $\erwfld$-torsors). 
%Suppose we are given $\Cat$-isomorphisms $\hg{\gal}:\Ve\lra\Ve^{\gal}$ 
%between $\Ve$ and its conjugate $\Ve^{\gal}$ 
%for each $\gal\in\Gal\pn{\erwfld/\fld}$. 
%If $\Ve$ is a homogeneous space for an algebraic $\erwfld$-group 
%or satisfies certain other criteria (see \cite[V, No.~20]{S}) 
%and the $\hg{\gal}$ satisfy the identity
%\bDpl 
%\hg{\gal\gall} = \lrpn{\hg{\gal}}^{\gall} \circ \hg{\gall} \laurink 
%\eDpl 
%then there exists a $\fld$-variety $\Vo$ and a $\Cat$-isomorphism 
%$ \viso:\Vo\tens_{\fld}\erwfld \iso \Ve $. 
%Here $\Vo$ inherits the structure of $\Ve$ preserved under the $\hg{\gal}$. 
%In this case it holds $\hg{\gal} = \viso^{\gal} \circ \viso^{-1}$. 
%
%\bRmk 
%\label{fmlGroup_descent}
Due to Cartier duality between formal groups and affine groups 
($\see$Thm.~\ref{Cartier-dual}), 
Galois descent applies to formal groups as well. 
%\eRmk 

When one does not assume that $\Xo$ is endowed with a $\fld$-rational 
point, one is led to two different descents of $\AlbFb{\Xb}$: 

(1) The universal mapping property of 
$\albFb:\Xb\dra\AlbFb{\Xb}$ gives for every $\gal\in\Gal\pn{\clfld/\fld}$
transformations 
$\hgg{1}{\gal}: \AlbFb{\Xb} \lra \AlbFb{\Xb}^{\gal}$ 
between $\AlbFb{\Xb}$ and its conjugate $\AlbFb{\Xb}^{\gal}$, 
which are homomorphisms of torsors. 
%i.e.\ compositions of a homomorphism by a translation. 
Therefore the descent of $\AlbFb{\Xb}$ by means of the $\hgg{1}{\gal}$ 
yields a $\fld$-torsor $\AlbbF{1}{\Xo}$. 

(2) In order to avoid translations or the reference to base points, 
one may reformulate the universal mapping property, 
replacing rational maps $\phe:\Xb\dra G$ from $\Xb$ to algebraic groups 
by its associated ``difference maps'' 
$\phe^{(-)}:\Xb\tms\Xb \dra G$,\; $(p,q) \lmt \phe(q)-\phe(p)$. 
In this way translations are eliminated and one obtains transformations 
$\hgg{0}{\gal}: \AlbFb{\Xb} \lra \AlbFb{\Xb}^{\gal}$ 
which are homomorphisms of algebraic groups. 
Then the descent of $\AlbFb{\Xb}$ by means of the 
$\hgg{0}{\gal}$ yields an algebraic $\fld$-group $\AlbbF{0}{\Xo}$. 
This is the $\fld$-group acting on the $\fld$-torsor $\AlbbF{1}{\Xo}$. 

\bNot 
\label{Difference map} 
If $\phe:\Xo\dra P$ is a rational map to a torsor 
%(= principal homogeneous space) 
$P$ for an algebraic group $G$, 
then 
\bDpl 
\phe^{(-)} : \Xo\tms\Xo \dra G 
\eDpl 
denotes the rational map to the algebraic group $G$ that assigns 
for $S \in \Algk$ to $(p,q)\in\Xo(S)\tms\Xo(S)$ the unique $g\in G(S)$ 
such that $g\cdot\phe(p)=\phe(q)$. 
\eNot 

\bNot 
\label{phe^(i)} 
If $\phe:\Xo\dra P$ is a rational map to a torsor, then set 
\bDpl 
\phe^{(1)} := \phe \,\laurink  \hspace{10pt}  \phe^{(0)} := \phe^{(-)} \laurin 
\eDpl 
%\begin{eqnarray*}
%\phe^{(1)} & := & \phe \\ 
%\phe^{(0)} & := & \phe^{(-)} 
%\end{eqnarray*}
\eNot 

\bThm 
\label{descent}
There exists a $\fld$-torsor $\AlbbF{1}{\Xo}$ for an algebraic 
$\fld$-group $\AlbbF{0}{\Xo}$ and rational maps defined over $\fld$ 
\[ \albbF{i}: \Xo^{2-i} \dra \AlbbF{i}{\Xo} 
\] 
for $i=1,0$, satisfying the following universal property: 

%Let $\anyfld$ be an extension of finite type of $\fld$ in $\unidom$ 
%and let $\clanyfld$ be the algebraic closure (in $\unidom$). 
If $\phe: \Xo \dra \Trr{1}$ 
is a rational map defined over $\fld$ to a $\fld$-torsor $\Trr{1}$ 
for an algebraic $\fld$-group $\Trr{0}$ 
which is an object of $\MrF$, 
%such that $\phe\tens_{\fld}\clfld$ is an object of $\MrF\lrpn{\Xb}$, 
then there exist a unique homomorphism of $\fld$-torsors 
$\homm{1}: \AlbbF{1}{\Xo} \lra \Trr{1}$ 
and a unique homomorphism of algebraic $\fld$-groups 
$\homm{0}: \AlbbF{0}{\Xo} \lra \Trr{0} \laurink$ 
defined over $\fld$, such that 
\bDpl \phe^{(i)} = \homm{i} \circ \albbF{i} 
\eDpl 
for $i=1,0$. 

The algebraic group $\AlbbF{0}{X}$ 
is dual to the 1-motive $\bigbt{\fmlG\lra\Picor{X}}$. 
\eThm 

\bPf 
Galois descent. 
The same arguments as given in 
\cite[V, No.~22]{S} work in our situation. 
\ePf

\subsubsection{Functoriality} 
\label{Functoriality} 

%\bPrp 
\label{pull-back_of_1-motives} 
Let $\fmlG \subset \Divor{X}$ be a dual-algebraic formal $\fld$-group. 
Let $\morphism: Y \lra X$ be a $\fld$-morphism of smooth proper $\fld$-varieties, 
such that no irreducible component of $\morphism\pn{Y}$ 
is contained in $\Supp(\fmlG)$. 
%Let $\lul\cut Y:\Decf_{X,Y} \lra \Divf_Y$ 
%be the pull-back of Cartier divisors from $X$ to $Y$. 
%($\see$Definition \ref{Dec_X,Y}, Proposition \ref{_.Y}). \\ 
For each dual-algebraic formal $\fld$-group $\fmlGr \subset \Divor{Y}$ 
containing $\morphism^*\fmlG$ 
the pull-back of relative Cartier divisors and of line bundles 
induces a homomorphism of 1-motives 
\[ \lrbt{\fmlGr \ra \Picor{Y}} \lla \lrbt{\fmlG \ra \Picor{X}} \laurin 
\] 
%\[ \left[ \begin{array}{c} 
%          \fmlGr \\ \downarrow \\ \Picor{Y} \\
%          \end{array} 
%   \right]
%   \lla 
%   \left[ \begin{array}{c} 
%          \fmlG \\ \downarrow \\ \Picor{X} \\
%          \end{array} 
%   \right] \laurin 
%\] 
%\ePrp 

According to the construction of universal objects over $\clfld$ 
($\see$Remark \ref{Alb_constr}), 
we obtain via dualization of 1-motives and by Galois descent 

\bPrp 
\label{alb_F(morphism)}
%Let $\fmlG \subset \Divor{X}$ be a dual-algebraic formal group. 
%Let $\morphism: Y \lra X$ be a morphism of smooth proper varieties, 
%such that no irreducible component of $\morphism\pn{Y}$ 
%is contained in $\Supp(\fmlG)$. 
%%such that $\morphism\pn{Y}$ intersects $\Supp(\fmlG)$ properly. 
Using the assumptions from above, 
then $\morphism$ induces 
for every formal group $\fmlGr \subset \Divor{Y}$ 
containing $\morphism^*\fmlG$ 
%satisfying $\fmlGr \supset \fmlG\cut Y$ 
a homomorphism of $\fld$-torsors 
$\Albb{1}^{\fmlG}_{\fmlGr}\pn{\morphism}$ 
and a homomorphism of algebraic $\fld$-groups 
$\Albb{0}^{\fmlG}_{\fmlGr}\pn{\morphism}$, 
\[ \Albb{i}_{\fmlGr}^{\fmlG}\pn{\morphism}: 
   \Alb^{(i)}_{\fmlGr}\pn{Y} \lra \Alb^{(i)}_{\fmlG}\pn{X}  
\hspace{8mm} \textrm{for } i = 1,0 \laurin 
\] 
%for $i = 1,0$. 
\ePrp

%\newpage 

\section{Albanese with modulus} 
\label{sec:Modulus}

Let $X$ be a smooth proper variety over a perfect field $k$. 
%an algebraically closed field of arbitrary characteristic. 
Let $\mdl$ %= \sum n_i \mdl_i$ 
be an effective divisor on $X$ 
(with multiplicity). 
The Albanese $\Albbm{1}{X}{\mdl}$ of $X$ of modulus $\mdl$ is 
a higher dimensional analogue 
to the generalized Jacobian with modulus of Rosenlicht-Serre. 
$\Albbm{1}{X}{\mdl}$ is defined by the universal mapping property 
for morphisms from $X\setminus\mdl$ to torsors 
of modulus $\leq\mdl$ 
($\see$Def.~\ref{DefMod}).  
Our definition of the modulus of rational maps %into algebraic groups 
coincides with the classical definition from \cite[III, \S~1]{S} in the curve case. 
Therefore the Albanese with modulus agrees 
with the Jacobian with modulus of Rosenlicht-Serre for curves, 
which we review in Subsection \ref{subsec:JacMod}. 

In Subsection \ref{subsec:ChowMod} we consider a Chow group 
$\CHmo{X}{\mdl}$ of 0-cycles relative to the modulus $\mdl$ 
($\see$Definition \ref{DefCHoMod}). 
We give an alternative characterization of $\Albm{X}{\mdl}$ 
as a universal quotient of $\CHmo{X}{\mdl}$, 
when the base field is algebraically closed 
($\see$Theorem \ref{ThmCHoMod}).

\subsection{Filtrations of the Witt group} 
\label{sec:filWitt}

Here we present a global version of some basic notions from \cite{KR2} 
that are needed for the construction of the Albanese with modulus. 

Let $(K,\val)$ be a discrete valuation field of characteristic $p>0$ 
with residue field $k$. 
The group of Witt vectors of length $r$ is denoted by $\Witt_r$ 
($\see$Exm.~\ref{WittVectors}). 
%We define filtrations $\fil_{n}\Witt_r(K)$ and $\filf_{n}\Witt_r(K)$, $n \in \Nat$, 
%on the group $\Witt_r(K)$ of Witt vectors of length $r$. 

\bDef 
\label{fil_n Witt}
Let $\fil_{n}\Witt_r(K)$ be the subgroup of $\Witt_r(K)$ 
from \cite[No.~1, Prop.~1]{Br}\laurind  
\[ \fil_{n}\Witt_r(K) = 
   \left\{ \lrpn{f_{r-1},\ldots,f_{0}} 
   \;\left|\;  
   \begin{array}{l} 
   f_i \in K, \quad \val\lrpn{f_i} \geq -n/p^i \\ 
   \forall \,0 \leq i \leq r-1 
   \end{array}
   \right.\right\} \laurin 
\] 
%\eDef 
%\bDef 
%\label{filf_n Witt}
%This filtration appears in 
Let $\filf_{n}\Witt_r(K)$ be the subgroup of $\Witt_r(K)$ 
generated by $\fil_{n}\Witt_r(K)$ by means of the Frobenius $\Frob$, 
see \cite[2.2]{KR2}\laurink 
\[ \filf_{n}\Witt_r(K) = 
   \sum_{\nu\geq 0} \Frob^{\nu} \fil_{n}\Witt_r(K) \laurin 
\] 
\eDef 

Let $X$ be a variety over $k$, 
regular in codimension 1. 
Let $\mdl=\sum_{\pnt\in S} n_{\pnt} \mdl_{\pnt}$ 
be an effective divisor on $X$, 
where $S$ is a finite set of points of codimension 1 in $X$, 
where $\mdl_{\pnt}$ are the prime divisors associated to $\pnt\in S$ 
and $n_{\pnt}$ are integers $\geq 1$ for $\pnt\in S$. 

\bDef 
\label{fil_D Witt}
Let $\fil_{\mdl}\Witt_r(\sK_X)$ (resp.\ $\filf_{\mdl}\Witt_r(\sK_X)$) 
be the sheaf of subgroups of 
$\Witt_r(\sK_X)$ formed by the groups 
\[ \bigpn{\fil_{\mdl}\Witt_r(\sK_X)}(U) = 
   \left\{ w \in \Witt_r(\sK_X) \;\left|\; 
   \begin{array}{l}
   w \in \fil_{n_{\pnt}}\Witt_r(\sK_{X,\pnt}) \quad \forall\pnt\in S\cap U \\
   w \in \Witt_r(\sO_{X,\pntt}) \hspace{31pt} \forall\pntt\in U\setminus S 
   \end{array} 
   \right.\right\}
\] 
resp.\ 
\[ \bigpn{\filf_{\mdl}\Witt_r(\sK_X)}(U) = 
   \left\{ w \in \Witt_r(\sK_X) \;\left|\; 
   \begin{array}{l}
   w \in \filf_{n_{\pnt}}\Witt_r(\sK_{X,\pnt}) \quad \forall\pnt\in S\cap U \\
   w \in \Witt_r(\sO_{X,\pntt}) \hspace{31pt} \forall\pntt\in U\setminus S 
   \end{array} 
   \right.\right\}
\] 
for open $U \in X$, 
where $\fil_{n}\Witt_r(\sK_{X,\pnt})$ (resp.\ $\filf_{n}\Witt_r(\sK_{X,\pnt})$) 
denotes the filtration associated to the valuation $\val_{\pnt}$ 
attached to the point $\pnt\in S$. 
\eDef 

\bPrp 
\label{finite fil_D}
Suppose $X$ is a projective variety over $k$ 
and $\mdl$ an effective divisor on $X$. 
Then $\Gam\bigpn{X,\fil_{\mdl}\Witt_r(\sK_X)}$ 
is a finite $\Witt_r(k)$-module. 
\ePrp 

\bPf 
The Verschiebung $\Ver: \Witt_{r-1} \lra \Witt_r$ 
yields an exact sequence 
$$0 \lra \fil_{\mdl}\Witt_{r-1}(\sK_X) \lra \fil_{\mdl}\Witt_{r}(\sK_X) \lra 
  \fil_{\lfloor\mdl/p^{r-1}\rfloor}\Witt_{1}(\sK_X) \lra 0$$ 
where $\lfloor\mdl/p^{r-1}\rfloor = \sum_{\pnt\in S} \lfloor n_{\pnt}/p^{r-1}\rfloor \mdl_{\pnt}$. 
This induces the exact sequence   
$$0 \lra \Gam\bigpn{\fil_{\mdl}\Witt_{r-1}(\sK_X)} \lra 
  \Gam\bigpn{\fil_{\mdl}\Witt_{r}(\sK_X)} \lra 
  \Gam\bigpn{\fil_{\lfloor\mdl/p^{r-1}\rfloor}\Witt_{1}(\sK_X)} \laurin 
$$ 
By induction over $r\geq 1$ and since $\Witt_1(k) = k$ is noetherian, 
it is sufficient to show the statement for $r=1$. 
Now $\fil_{\mdl}\Witt_1(\sK_X) = \sO_X(\mdl)$ is a coherent sheaf, 
hence $\Gam\bigpn{X,\fil_{\mdl}\Witt_1(\sK_X)}$ is a finite module 
over $\Witt_1(k) = k$.  
\ePf 

%\vspace{\vs} 

%Let $p$ be a prime number. 
\bDef 
\label{R[F]}
Let $R$ be a commutative ring over $\bF_p$. 
We let $R[\Frob]$ be the non-commutative polynomial ring defined by 
\[ R[\Frob] = \lrst{\left. \sum_{i=1}^n \Frob^i r_i \; \right| \; 
  \begin{array}{c} 
  r_i \in R \\  
  n \in \Nat 
  \end{array}
  } \laurink 
  \hspace{4em} \Frob r = r^p \Frob \hspace{2em} \forall r \in R \laurin  
\]
\eDef 

\bDef 
\label{pder}
If $\Kot_{\sK_X} = \Kot_{\sK_X/k}$ 
is the module of differentials of $\sK_X$ over $k$, 
we let $\pder$ be the homomorphism 
\[ \pder: \Witt_r(\sK_X) \lra \Kot_{\sK_X} \laurink
   \hspace{4em} 
   \lrpn{f_{r-1},\ldots,f_0} \lmt \sum_{i=0}^{r-1} f_i^{p^i-1} \der f_i \laurin 
\] 
\eDef 

\bDef 
\label{Kot(log)}
If $\Es$ is a reduced effective divisor on $X$ with normal crossings, 
we let $\Kot_{X}\pn{\log\Es}$ be the sheaf of differentials on $X$ 
with $\log$-poles along $\Es$, 
i.e.\ the $\sO_X$-module generated locally by $\der f$, $f \in \sO_X$ 
and $\der\log t = t^{-1} \der t$, where $t$ is a local equation for $\Es$. 
\eDef 

\bPrp 
\label{fil-quot}
Suppose $\mdl_{\red}$ is a normal crossing divisor. 
The homomorphism $\pder$ from Definition \ref{pder} 
induces injective homomorphisms 
\begin{eqnarray*}
\qfil_{\mdl}: \: 
\filf_{\mdl} \Witt_r(\sK_X) / \filf_{\lfloor\mdl/p\rfloor} \Witt_r(\sK_X) 
 & \linj & \Qfil_{\mdl} \\ 
\qfill_{\mdl}: \: 
\filf_{\mdl} \Witt_r(\sK_X) / \filf_{\mdl-\mdl_{\red}} \Witt_r(\sK_X) 
 & \linj & \Qfill_{\mdl} 
\end{eqnarray*}
where $\Qfil_{\mdl}$ resp.\ $\Qfill_{\mdl}$ are the $\sO_X$-modules 
\begin{eqnarray*}
\Qfil_{\mdl} & = & k[\Frob] \tens_k \bigpn{\Kot_{X}\pn{\log\mdl_{\red}} 
\tens_{\sO_X} \sO_X(\mdl)/\sO_X(\lfloor\mdl/p\rfloor)} \laurink \\ 
\Qfill_{\mdl} & = & k[\Frob] \tens_k \bigpn{\Kot_{X}\pn{\log\mdl_{\red}} 
\tens_{\sO_X} \sO_X(\mdl)/\sO_X\pn{\mdl-\mdl_{\red}}} \laurin 
\end{eqnarray*}
and $\lfloor\mdl/p\rfloor$ means the largest divisor $\mdll$ 
such that $p\mdll \leq \mdl$. 
\ePrp

\bPf 
This is the global formulation of \cite[4.6]{KR2}. 
\ePf 

\bDef 
\label{flat-fil_D}
Let $\Qfillb_{\mdl}$ be the image in $\Qfill_{\mdl}$ 
of the $\sO_X$-module \\ 
\;$k[\Frob] \tens_k \bigpn{\Kot_{X} 
\tens_{\sO_X} \sO_X(\mdl)/\sO_X\pn{\mdl-\mdl_{\red}}}$\; 
(without log-poles). 
Then 
\[ \Qfillb_{\mdl} \cong k[\Frob] \tens_k \bigpn{\Kot_{\mdl_{\red}} 
\tens_{\sO_{\mdl_{\red}}} \sO_X(\mdl)/\sO_X\pn{\mdl-\mdl_{\red}}}
\] 
since $t^{-n_{\pnt}} \der t = t^{1-n_{\pnt}} \der\log t$ vanishes in $\Qfill_{\mdl}$ 
for any local equation $t$ of $\mdl_{\red}$. 
Then we let 
\,$\filfb_{\mdl} \Witt_r(\sK_{X}) \subset \filf_{\mdl} \Witt_r(\sK_{X})$\, 
be the inverse image of $\Qfillb_{\mdl}$ under the map $\qfill_{\mdl}$ 
from Proposition \ref{fil-quot}. 
According to \cite[4.7]{KR2}, this is a global version 
of the following alternative definition: 
\eDef 

\bDef 
\label{flat-fil_n}
Let $\filb_{n}\Witt_r(K)$ be the subgroup of $\fil_{n}\Witt_r(K)$ 
consisting of all elements $\lrpn{f_{r-1},\ldots,f_0}$ 
satisfying the following condition: 
If the $p$-adic order $\nu$ of $n$ is $< r$, 
then $p^{\nu} \val(f_{\nu}) > -n$. 
Then $\filfb_{n}\Witt_r(K)$ is the subgroup of $\Witt_r(K)$ 
generated by $\filb_{n}\Witt_r(K)$ by means of the Frobenius $\Frob$\laurink 
\[ \filfb_{n}\Witt_r(K) = 
   \sum_{\nu\geq 0} \Frob^{\nu} \,\filb_{n}\Witt_r(K) \laurin 
\] 
\eDef 

\bLem 
\label{fil-restriction}
Let $\morphism: Y \lra X$ be a morphism of varieties over $k$, 
such that $\morphism\pn{Y}$ intersects $\Supp\pn{\mdl}$ properly. 
Let $\mdl \cut Y$ denote the pull-back of $\mdl$ to $Y$. 
Suppose that $\mdl_{\red}$ and $\pn{\mdl \cut Y}_{\red}$ 
are normal crossing divisors. 
There is a commutative diagram of homomorphisms 
with injective rows 
\[ \xymatrix{ 
\filf_{\mdl} \Witt_r(\sK_X) / \filf_{\lfloor\mdl/p\rfloor} \Witt_r(\sK_X) 
\ar[r]^-{\qfil_{X,\mdl}} \ar[d] & \Qfil_{X,\mdl} \ar[d] \\ 
\filf_{\mdl\cut Y} \Witt_r(\sK_Y) / \filf_{\lfloor\mdl\cut Y/p\rfloor} \Witt_r(\sK_Y) 
\ar[r]^-{\qfil_{Y,\mdl\cut Y}} & \Qfil_{Y,\mdl\cut Y}
}
\] 
where the vertical arrows are the obvious pull-back maps from $X$ to $Y$. 
\eLem 

\bPf 
Straightforward. 
\ePf

\subsection{Albanese with modulus} 
\label{subsec:AlbMod}

\subsubsection{Existence and construction} 
\label{Alb(X,D)exist}

Let $X$ be a %irreducible 
smooth proper variety over a perfect field $k$. 
%which we assume to be algebraically closed in this No. 
%Arbitrary perfect base field is considered in 
%No.~\ref{Alb(X,D)descent} and No.~\ref{Alb(X,D)functorial}. 

\bDef 
\label{DefMod}
First assume $k$ is algebraically closed. 
Let $\phe:X\dra G$ be a rational map from $X$ 
to a smooth connected algebraic group $G$. 
Let $L$ be the affine part of $G$ 
and $\Upt$ the unipotent part of $L$. 
The \emph{modulus of $\phe$} from \cite[\S~3]{KR2} 
is the following effective divisor 
\[ \modu\lrpn{\phe} = \sum_{\hgt(\pnt)=1} \modu_{\pnt}(\phe) \;\mdl_{\pnt} 
\] 
where $\pnt$ ranges over all points of codimension 1 in $X$, 
and $\mdl_{\pnt}$ is the prime divisor associated to $\pnt$.  
For each $\pnt\in X$ of codimension 1, 
the canonical map 
$L\lrpn{\sK_{X,\pnt}}/L\lrpn{\sO_{X,\pnt}} \lra 
  G\lrpn{\sK_{X,\pnt}}/G\lrpn{\sO_{X,\pnt}}$ 
is bijective, cf.\ \ref{sec:indTrafo}. 
Take an element $l_{\pnt} \in L\lrpn{\sK_{X,\pnt}}$ 
whose image in $G\lrpn{\sK_{X,\pnt}}/G\lrpn{\sO_{X,\pnt}}$ 
coincides with the class of $\phe \in G\lrpn{\sK_{X,\pnt}}$. 
If $\chr(k) = 0$, let $(u_{\pnt,i})_{1\leq i\leq s}$ be the image of $l_{\pnt}$ 
in $\Ga(\sK_{X,\pnt})^s$ under $L\to\Upt\cong (\Ga)^s$. 
If $\chr(k) = p > 0$, 
let $(u_{\pnt,i})_{1\leq i\leq s}$  be the image of $l_{\pnt}$ in 
$\Witt_r(\sK_{X,\pnt})^s$ under  $L \to \Upt \subset \lrpn{\Witt_r}^s$. 
\[ \modu_{\pnt}(\phe) = \left\{ 
   \begin{array}{ll} 
      0 & \textrm{if } \; \phe \in G\lrpn{\sO_{X,\pnt}} \\ 
      1+\max\left\{\nty_{\pnt}(u_{\pnt,i}) \;|\; 1\leq i\leq s \right\} & 
         \textrm{if } \; \phe \notin G\lrpn{\sO_{X,\pnt}} 
   \end{array} 
   \right. 
\] 
where for $u\in\Ga\pn{\sK_{X,\pnt}}$ resp. $\Witt_r\pn{\sK_{X,\pnt}}$ 
\[ \nty_{\pnt}(u) = \left\{ 
   \begin{array}{ll} 
      -\val_{\pnt}(u) & \textrm{if} \; \chr(k) = 0 \\ 
      \min\left\{n\in\Nat \;|\; u\in \filf_{n} \Witt_{r}\pn{\sK_{X,\pnt}}\right\} & 
         \textrm{if} \; \chr(k) = p > 0 \laurin 
   \end{array} 
   \right. 
\] 
\eDef 

\noindent 
Note that $\modu_v(\phe)$ is independent of the choice of the isomorphism $\Upt \cong (\Ga)^s$
resp.\ of the embedding $\Upt \subset \lrpn{\Witt_r}^s$, 
see \cite[Thm.~3.3]{KR2}. 

For arbitrary perfect base field $k$ 
and $G$ a torsor for a smooth connected algebraic group 
we obtain $\modu\lrpn{\phe}$ 
by means of a Galois descent from $\modu\lrpn{\phe \tens_k \clfld}$, 
where $\clfld$ is an algebraic closure of $k$, 
see \cite[No.~3.4]{KR2}, cf.\ also Rmk.~\ref{k-ratl pt of torsor}. 

\bDef 
\label{DefAlb(X,D)}
Let $\mdl$ be an effective divisor on $X$. 
Then $\Mrm{X}{\mdl}$ denotes the category of those rational maps 
$\phe$ from $X$ to torsors 
such that \;$\modu\lrpn{\phe} \leq \mdl$. 
The universal object of $\Mrm{X}{\mdl}$ (if it exists) 
is called the \emph{Albanese of $X$ of modulus $\mdl$} 
and denoted by $\Albbm{1}{X}{\mdl}$, 
or just $\Albm{X}{\mdl}$, if it admits a $\fld$-rational point 
(cf.\ Rmk.~\ref{k-ratl pt of torsor}). 
\eDef 

\bRmk 
\label{saturationFrob}
In the definition of $\modu(\phe)$ ($\see$Def.~\ref{DefMod}) 
we used, instead of the original filtration $\fil_{\bullet} \Witt$ of Brylinski, 
the saturation $\filf_{\bullet} \Witt$ of $\fil_{\bullet} \Witt$ w.r.t.\ the Frobenius. 
This is motivated as follows: 
Let $\modun(\phe)$ be the modulus of a rational map $\phe$ 
using the filtration $\fil_{\bullet} \Witt$ instead of $\filf_{\bullet} \Witt$. 
If $\phe: X \dra \Ga$ is a non-constant rational map, 
i.e.\ the multiplicity of $\modu(\phe) =: \mdl$ is $> 1$, then 
$\modun(\Frob^{\nu} \circ \,\phe) = p^{\nu}\pn{\mdl - \mdl_{\red}} + \mdl_{\red}$, 
where $\mdl_{\red}$ is the reduced part of $\mdl$. 
On the other hand, if $u: X \dra \sU$ is a universal object 
for a certain category of rational maps $\Mr$, 
then clearly $u$ satisfies the universal mapping property as well 
for all maps of the form $\Frob^{\nu} \circ \,\phe$, $\phe \in \Mr$ 
(cf.\ condition $\lrpn{ \diamondsuit \; 2 }$ before Thm.~\ref{Exist univObj}). 
This shows that $\modun(\phe)$ is not compatible 
with the notion of universal objects. 
\eRmk 

\bDef 
\label{DefFm(X,D)}
Let $\mdl$ be an effective divisor on $X$, 
let $\mdl_{\red}$ be the reduced part of $\mdl$. 
Then $\Fm{X}{\mdl}$ denotes the formal subgroup of $\Divf_X$ 
characterized by 
\[ \lrpn{\Fm{X}{\mdl}}_{\et} = 
   \left\{ B\in\Divf_{X}(k) \;\big|\; \Supp(B) \subset \Supp(\mdl) \right\}
\] 
and for $\chr(k)=0$ 
\[ \lrpn{\Fm{X}{\mdl}}_{\inf} = \exp \lrpn{ \Gac \tens_k 
   \Gam\bigpn{\sO_{X}\lrpn{\mdl-\mdl_{\red}} \big/ \sO_{X}}} 
\] 
\phantom{and} for $\chr(k)=p>0$ 
\[ \lrpn{\Fm{X}{\mdl}}_{\inf} = 
   \Expah \lrpn{ \sum_{r > 0} \Wcfl{r} \tens_{\Witt_r(k)} 
   \Gam\Bigpn{\filf_{\mdl-\mdl_{\red}} \Witt_r(\sK_{X}) \Big/ \Witt_r(\sO_{X})},1}  \laurin 
\] 
Let $\Fmor{X}{\mdl} = \Fm{X}{\mdl} \tms_{\Divf_X} \Divor{X}$ 
be the intersection of $\Fm{X}{\mdl}$ and $\Divor{X}$. 
\eDef 

\bPrp 
\label{Fm(X,D)algebraic}
The formal groups $\Fm{X}{\mdl}$ and $\Fmor{X}{\mdl}$ are dual-algebraic. 
\ePrp 

\bPf 
The statement is obvious for $\chr(k)=0$, 
therefore we suppose $\chr(k)=p>0$. 
The proof is done in two steps. 
Let $\Fmn{X}{\mdl}$ be the formal subgroup of $\Divf_X$ 
defined in the same way as $\Fm{X}{\mdl}$, 
but using the filtration $\fil_{\mdl-\mdl_{\red}} \Witt_r(\sK_{X})$ 
instead of $\filf_{\mdl-\mdl_{\red}} \Witt_r(\sK_{X})$. 
In the first step, we show that for any effective divisor $\mdl$ 
the formal group $\Fmn{X}{\mdl}$ is dual-algebraic. 
In the second step, we show that for any $\mdl$ 
there exists $\mdl' \geq \mdl$ such that $\Fm{X}{\mdl}$
is contained in the image of $\Fmn{X}{\mdl'}$ in $\Fm{X}{\mdl'}$. 
Thus $\Fm{X}{\mdl}$ is a formal subgroup of a quotient of a 
dual-algebraic formal group, hence dual-algebraic 
by Lemma \ref{subFmlGroup_algebraic}. 
Then also the formal subgroup $\Fmor{X}{\mdl}$ of $\Fm{X}{\mdl}$ 
is dual-algebraic. 

\textbf{Step 1:} 
Let $\mdl$ be an effective divisor on $X$. 
Write $\mdl = \sum_{\hgt(\pnt)=1} n_{\pnt} \, \mdl_{\pnt}$, 
where $\pnt$ ranges over all points of codimension 1 in $X$, 
and $\mdl_{\pnt}$ is the prime divisor associated to $\pnt$.  
Let $S$ be the finite set of those $\pnt$ with $n_{\pnt} > 0$. 
Let $m = \min\left\{r \;|\; p^r > n_{\pnt} - 1 \; \forall \pnt \in S \right\}$. 
Hence for $r \geq m$, if 
$\lrpn{f_{r-1},\ldots,f_0} \in \fil_{\mdl-\mdl_{\red}} \Witt_r(\sK_{X})$, 
then $f_i \in \sO_X$ for $r > i \geq m$, 
according to Definition \ref{fil_D Witt}. 
%of $\fil_{\mdl-\mdl_{\red}} \Witt_r(\sK_{X})$. 
Then the Verschiebung 
$\Ver^{r-m}: \Witt_m(\sK_X) \lra \Witt_r(\sK_X)$ 
yields a surjective homomorphism 
%on the filtered subgroups modulo integral elements, 
$\fil_{\mdl-\mdl_{\red}} \Witt_m(\sK_{X}) / \Witt_m(\sO_{X}) \lsur 
  \fil_{\mdl-\mdl_{\red}} \Witt_r(\sK_{X}) / \Witt_r(\sO_{X})$. 
Thus $\lrpn{\Fmn{X}{\mdl}}_{\inf}$ is already generated 
by a finite sum via $\Expah$: 
\[ \lrpn{\Fmn{X}{\mdl}}_{\inf} = 
   \Expah \lrpn{ \sum_{1 \leq r \leq m} \Wcfl{r} \tens_{\Witt_r(k)} 
\Gam\Bigpn{\fil_{\mdl-\mdl_{\red}} \Witt_r(\sK_{X}) \Big/ \Witt_r(\sO_{X})},1} 
\laurin 
\] 
Each 
$\Gam\bigpn{X,\fil_{\mdl-\mdl_{\red}} \Witt_r(\sK_{X}) \big/ \Witt_r(\sO_{X})}$ 
is a finitely generated $\Witt_r(k)$-module, 
by the same proof as for Proposition \ref{finite fil_D}. 
Hence $\lrpn{\Fmn{X}{\mdl}}_{\inf}$ is a quotient of the direct sum 
of finitely many $\Wcfl{r}$. %, $r\in\Nat$. 

Moreover, 
$\lrpn{\Fmn{X}{\mdl}}_{\et} = \lrpn{\Fm{X}{\mdl}}_{\et}$ 
is an abelian group of finite type, 
since $\mdl$ has only finitely many components. 
Thus $\Fmn{X}{\mdl}$ is dual-algebraic, 
according to Proposition \ref{dual-algebraic}. 

\textbf{Step 2:} 
We show that for any effective divisor $\mdl$ there exists 
an effective divisor 
$\mdl' \geq \mdl$ such that 
$\Fm{X}{\mdl}$ is generated by $\Fmn{X}{\mdl'}$. 
We will find an effective divisor 
$\mdl'\geq\mdl$ such that 
$\Gam\bigpn{\filf_{\mdl-\mdl_{\red}} \Witt_r(\sK_{X}) \big/ \Witt_r(\sO_{X})}$
is generated by $\sum_{\nu\geq 0} \Frob^{\nu} 
\Gam\bigpn{\fil_{\mdl'-\mdl'_{\red}} \Witt_r(\sK_{X}) \big/ \Witt_r(\sO_{X})}$. 
This is sufficient because 
$\Expah \lrpn{v \tens \sum_i \Frob^{\nu_i}\oma_i,1} = 
\Expah \lrpn{\sum_i \Ver^{\nu_i}v \tens \oma_i,1}$.  
Since the homomorphism %induced by the Verschiebung 
$\Ver^{r-m}: \filf_{\mdl-\mdl_{\red}} \Witt_m(\sK_{X}) / \Witt_m(\sO_{X}) 
  \lra  \filf_{\mdl-\mdl_{\red}} \Witt_r(\sK_{X}) / \Witt_r(\sO_{X})$ 
is surjective for $r \geq m$, we only need to consider $r=m$. 

The exact sequence 
\[ 0 \lra \Witt_r(\sO_X) \lra \Witt_r(\sK_X) \lra 
   \Witt_r(\sK_X) \big/ \Witt_r(\sO_X) \lra 0
\] 
yields the exact sequence 
\[ \Gam\bigpn{\Witt_r(\sK_X)} \lra 
   \Gam\bigpn{\Witt_r(\sK_X) \big/ \Witt_r(\sO_X)} \lra 
   \H^1\bigpn{\Witt_r(\sO_X)} \lra 0 \laurin 
\] 
Here $\H^1\bigpn{\Witt_r(\sK_X)} = 0$ since $\Witt_r(\sK_X)$ 
is a flasque sheaf. 
Since $\H^1\bigpn{\Witt_r(\sO_X)}$ is a finite $\Witt_r(k)$-module, 
there is an effective divisor $\mdll$ such that the map 
$\Gam\bigpn{\fil_{\mdll}\Witt_r(\sK_X)\big/\Witt_r(\sO_X)} \lra 
  \H^1\bigpn{\Witt_r(\sO_X)}$ 
is surjective. 
Hence for any $\sig \in 
\Gam\bigpn{\filf_{\mdl-\mdl_{\red}} \Witt_r(\sK_{X}) \big/ \Witt_r(\sO_{X})}$ 
there is $\rho \in \Gam\bigpn{\fil_{\mdll}\Witt_r(\sK_X)\big/\Witt_r(\sO_X)}$ 
such that $\sig-\rho$ lies in the image of $\Gam\bigpn{\Witt_r(\sK_X)}$, 
hence in the image of $\Gam\bigpn{\filf_{\mdll'}\Witt_r(\sK_X)}$, 
where $\mdll' = \max\lrst{\mdll,\mdl-\mdl_{\red}}$. 
Therefore we are reduced to showing that for any $\mdl$ there exists 
$\mdl'\geq\mdl$ such that 
$\Gam\bigpn{\filf_{\mdl} \Witt_r(\sK_{X})}$
is generated by $\sum_{\nu\geq 0} \Frob^{\nu} 
\Gam\bigpn{\fil_{\mdl'} \Witt_r(\sK_{X})}$. 

Consider the exact sequence 
\[ 0 \lra \dsum_{\nu\geq 0} \fil_{\lfloor\mdl/p\rfloor} \Witt_r(\sK_X) \lra 
   \dsum_{\nu\geq 0} \fil_{\mdl} \Witt_r(\sK_X) \lra 
   \filf_{\mdl} \Witt_r(\sK_X) \lra 0
\] 
where the third arrow is 
$\lrpn{w_{\nu}}_{\nu} \lmt \sum_{\nu} \Frob^{\nu}w_{\nu}$, 
and the second arrow is 
$\lrpn{w_{\nu}}_{\nu} \lmt \lrpn{\Frob w_{\nu} - w_{\nu-1}}_{\nu}$, 
where we set $w_{-1} = 0$. 
Here $\lfloor\mdl/p\rfloor$ means the largest divisor $\mdll$ 
such that $p\mdll \leq \mdl$. 
This yields an exact sequence 
\[ \dsum_{\nu\geq 0} \Gam\bigpn{\fil_{\mdl} \Witt_r(\sK_X)} \lra 
   \Gam\bigpn{\filf_{\mdl} \Witt_r(\sK_X)} \lra 
   \dsum_{\nu\geq 0} \H^1\bigpn{\fil_{\lfloor\mdl/p\rfloor} \Witt_r(\sK_X)} 
   \laurin 
\] 
$\Witt_r(\sK_X)$ is the inductive limit of $\fil_{\mdll} \Witt_r(\sK_X)$, 
where $\mdll$ ranges over all effective divisors on $X$, hence 
%\begin{eqnarray*}
\[ 
0  =  \H^1\bigpn{\Witt_r(\sK_X)} 
    =  \H^1\biggpn{\varinjlim_{\mdll} \fil_{\mdll} \Witt_r(\sK_X)} 
    =  \varinjlim_{\mdll} \H^1\bigpn{\fil_{\mdll} \Witt_r(\sK_X)} \laurin 
\]  
%\end{eqnarray*}
As $\H^1\bigpn{\fil_{\lfloor\mdl/p\rfloor} \Witt_r(\sK_X)}$ 
is a finite $\Witt_r(k)$-module, there is an effective divisor 
$\mdl' \geq \mdl$ such that the image of 
$\H^1\bigpn{\fil_{\lfloor\mdl/p\rfloor} \Witt_r(\sK_X)}$ in 
$\H^1\bigpn{\fil_{\lfloor\mdl'/p\rfloor} \Witt_r(\sK_X)}$ is $0$. 
Thus the image of $\Gam\bigpn{\filf_{\mdl} \Witt_r(\sK_X)}$ in 
$\Gam\bigpn{\filf_{\mdl'} \Witt_r(\sK_X)}$ is contained in 
$\sum_{\nu\geq 0} \Frob^{\nu} 
  \Gam\bigpn{\fil_{\mdl'} \Witt_r(\sK_X)}$. 
\ePf 

\bLem 
\label{mod(phe)-im(trafo_phe)}
Let $\phe:X\dra G$ be a rational map from $X$ to a smooth connected 
algebraic group $G$. Then the following conditions are equivalent: 

\begin{tabular}{rl} 
\hspace{\kukwad} (i) & $\modu\lrpn{\phe} \leq \mdl$\laurink \\ 
\hspace{\kukwad} (ii) & $\im\lrpn{\trafo_{\phe}} \subset \Fm{X}{\mdl}$\laurin 
\end{tabular} 
\eLem 

\bPf 
Write $\mdl = \sum_{\hgt(\pnt)=1} n_{\pnt} \, \mdl_{\pnt}$, 
where $\pnt$ ranges over all points in $X$ of codimension 1, 
and $\mdl_{\pnt}$ is the prime divisor associated to $\pnt$.  
Condition $(i)$ is thus expressed by the condition that 
for all $\pnt \in X$ of codimension 1 it holds 

\begin{tabular}{rl} 
%\hspace{1pt} 
\hspace{\kukwad} $(i)_{\pnt}$ & $\modu_{\pnt}(\phe) \leq n_{\pnt}$\laurin  
\end{tabular} 

\noindent 
Using the canonical splitting of a formal group into an \'etale 
and an infinitesimal part, condition $(ii)$ is equivalent 
to the condition that the following $(ii)_{\et}$ and $(ii)_{\inf}$ 
are satisfied: 

\begin{tabular}{ll} 
\hspace{\kukwad} $(ii)_{\et}$ & 
$\im\lrpn{\trafo_{\phe,\et}} \subset \lrpn{\fmlG_{X,\mdl}}_{\et}$\laurink \\ 
\hspace{\kukwad} $(ii)_{\inf}$ & 
$\im\lrpn{\trafo_{\phe,\inf}} \subset \lrpn{\fmlG_{X,\mdl}}_{\inf}$\laurin 
\end{tabular} 

\noindent 
Let $L$ be the affine part of $G$. 
Remember from \ref{sec:indTrafo}
%Definition \ref{induced Trafo} 
that the transformation $\trafo_{\phe}: \Ld \lra \Divfc{X}$ is given by 
$\lrpair{\llul,\lin_{\phe}}$, where $\lin_{\phe}$ 
is the image of $\phe\in G\lrpn{\sK_X}$ in 
$\Gam\bigpn{G(\sK_X)/G(\sO_X)} \iso 
  \Gam\bigpn{L(\sK_X)/L(\sO_X)}$, 
and the pairing 
\[ \lrpair{\llul,\lull} : \Ld \tms \Gam\bigpn{L(\sK_X)/L(\sO_X)} 
   \lra \Gam\bigpn{\Gm(\ul{\sK_X})/\Gm(\ul{\sO_X})} 
\] 
is obtained from Cartier duality. 
Write $L = \Trs \tms_k \Upt$ as a product of 
a torus $\Trs$ and a unipotent group $\Upt$. 
Fix an isomorphism $\Trs \cong \lrpn{\Gm}^m$ 
and an isomorphism $\Upt \cong \lrpn{\Ga}^a$ if $\chr(k) = 0$,  
resp.\ an embedding $\Upt \subset \lrpn{\Witt_r}^a$ if $\chr(k) = p > 0$.  
Let $\pn{t_j}_{1\leq j \leq m}$ be the image of $\lin_{\phe}$ under 
\[ \Gam\bigpn{L(\sK_X)/L(\sO_X)} \lra 
   \Gam\bigpn{\Trs(\sK_X)/\Trs(\sO_X)} \lra 
   \Gam\bigpn{\Gm(\sK_X)/\Gm(\sO_X)}^m 
\] 
and $\pn{u_i}_{1\leq i \leq a}$ be the image of $\lin_{\phe}$ under 
\[ \Gam\bigpn{L(\sK_X)/L(\sO_X)} \lra 
   \Gam\bigpn{\Upt(\sK_X)/\Upt(\sO_X)} \lra 
   \left\{ 
   \begin{array}{l}
   \Gam\bigpn{\Ga(\sK_X)/\Ga(\sO_X)}^a \\ 
   \Gam\bigpn{\Witt_r(\sK_X)/\Witt_r(\sO_X)}^a \laurin 
   \end{array}
   \right. 
\] 
The \'etale part of $\trafo_{\phe}$ is 
\begin{eqnarray*}
   \trafo_{\phe,\et}: \hspace{9pt} \Zint^m \hspace{9pt} 
   & \lra & \Gam\bigpn{\Gm(\sK_X)/\Gm(\sO_X)} \\ 
   \pn{e_j}_{1\leq j \leq m} & \lmt & \prod_{j=1}^m t_j^{e_j} \laurin 
\end{eqnarray*}
The image of the infinitesimal part of $\trafo_{\phe}$ is 
given by the image of 
\begin{eqnarray*}
 \wt{\trafo}_{\phe,\inf}: 
   \left. 
   \begin{array}{l} 
   \bigpn{\Gac}^a \\ 
   \bigpn{\Wcfl{r}}^a 
   \end{array} 
   \right\} 
   & \lra & 
   \Gam\bigpn{\Gm(\ul{\sK_X})/\Gm(\ul{\sO_X})} \\ 
%   \Gam\bigpn{\Upf_{\ghost}(\sK_X)/\Upf_{\ghost}(\sO_X)} \\ 
   \pn{v_i}_{1\leq i \leq a} 
   & \lmt & 
   \left\{ 
   \begin{array}{l}
   \prod_{i=1}^a \exp\pn{v_i \, u_i} \\ 
   \prod_{i=1}^a \Expah\pn{v_i \cdot u_i, 1} \laurink 
   \end{array}
   \right. 
\end{eqnarray*}
cf.\ Example~\ref{WittVectors} %\cite[V, \S~4, 4.5]{DG} 
for the pairing $\Wcfl{r} \tms \Witt_r \lra \Gm$. 
For each $\pnt\in X$ of codimension 1 
let $\lrpn{t_{\pnt,i}}_{1\leq j\leq m}$ be a representative in 
$\Gm\pn{\sK_{X,\pnt}}^m$ of the image of $\pn{t_j}_{1\leq j \leq m}$ 
under 
\[ \Gam\bigpn{\Gm(\sK_X)/\Gm(\sO_X)}^m \lra 
   \Gm(\sK_{X,\pnt})^m/\Gm(\sO_{X,\pnt})^m \laurink  
\] 
and let $\lrpn{u_{\pnt,i}}_{1\leq i\leq a}$ be a representative in 
$\Ga\pn{\sK_{X,\pnt}}^a$ resp.\ $\Witt_r\pn{\sK_{X,\pnt}}^a$ 
of the image of $\pn{u_i}_{1\leq i \leq a}$ under 
\[ \Gam\bigpn{\Ga(\sK_X)/\Ga(\sO_X)}^a \lra 
   \Ga(\sK_{X,\pnt})^a/\Ga(\sO_{X,\pnt})^a 
\] 
resp.\ 
\[ \Gam\bigpn{\Witt_r(\sK_X)/\Witt_r(\sO_X)}^a \lra 
   \Witt_r(\sK_{X,\pnt})^a/\Witt_r(\sO_{X,\pnt})^a \laurin 
\] 
Then $(ii)_{\et}$ is equivalent to the following condition 
being satisfied for every point $\pnt \in X$ of codimension 1: 

\begin{tabular}{ll} 
\hspace{\kukwad} $(ii)_{\et,\pnt}$ & 
If $n_{\pnt} = 0$, 
then $t_{\pnt,j} \in \Gm\pn{\sO_{X,\pnt}}$ for $1\leq j\leq m$. 
\end{tabular} 

\noindent 
On the other hand, $(ii)_{\inf}$ is equivalent to the following condition 
being satisfied for every point $\pnt \in X$ of codimension 1, 
according to Definition \ref{DefFm(X,D)} of $\Fm{X}{\mdl}$: 

\begin{tabular}{ll} 
\hspace{\kukwad} $(ii)_{\inf,\pnt}$ & 
If $n_{\pnt} = 0$, 
then $u_{\pnt,i} \in \Ga\pn{\sO_{X,\pnt}}$ resp.\ $\Witt_r\pn{\sO_{X,\pnt}}$
for $1\leq i\leq a$. \\ 
 & If $n_{\pnt} > 0$, 
then $\nty_{\pnt}\pn{u_{\pnt,i}} \leq n_{\pnt} - 1$. 
\end{tabular} 

\noindent 
Note that $\phe \in G\pn{\sO_{X,\pnt}}$ 
if and only if $t_{\pnt,j} \in \Gm\pn{\sO_{X,\pnt}}$ for $1\leq j\leq m$ 
and $u_{\pnt,i} \in \Ga\pn{\sO_{X,\pnt}}$ resp.\ $\Witt_r\pn{\sO_{X,\pnt}}$ 
for $1\leq i\leq a$. 
By Definition \ref{DefMod}, for each $\pnt\in X$ of codimension 1 

\begin{tabular}{ll} 
$(i)_{\pnt}$ & $\modu_{\pnt}\pn{\phe} \leq n_{\pnt}$ 
\end{tabular} 

\noindent 
if and only if $(ii)_{\et,\pnt}$ and $(ii)_{\inf,\pnt}$ are satisfied. 
\ePf 

%\newpage 
\vspace{\vs}

Now assume $k$ is algebraically closed. 
Arbitrary perfect base field is considered in 
No.~\ref{Alb(X,D)descent} and No.~\ref{Alb(X,D)functorial}. 

\bThm 
\label{Mdl-category}
The category $\Mrm{X}{\mdl}$ of rational maps of modulus $\leq\mdl$ 
is equivalent to the category $\Mrt{\Fm{X}{\mdl}}{X}$ 
of rational maps which induce a transformation to $\Fm{X}{\mdl}$. 
\eThm 

\bPf 
According to the definitions of $\Mrm{X}{\mdl}$ and $\Mrt{\Fm{X}{\mdl}}{X}$, 
the statement is due to Lemma \ref{mod(phe)-im(trafo_phe)}. 
\ePf 

\bThm 
\label{AlbMod-construction}
The Albanese $\Albm{X}{\mdl}$ of \,$X$ of modulus $\mdl$ exists 
and is dual (in the sense of 1-motives) to the 1-motive 
$\bigbt{\Fmor{X}{\mdl} \lra \Picor{X}}$. 
\eThm 

\bPf 
By Theorem \ref{Mdl-category}, $\Albm{X}{\mdl}$ is the universal object 
of the category $\Mrt{\Fm{X}{\mdl}}{X}$ (if it exists). 
A rational map from $X$ to an algebraic group 
induces a transformation to $\Fm{X}{\mdl}$ 
if and only if it induces a transformation to $\Fmor{X}{\mdl}$, 
by Lemma \ref{Ld->Picor}. 
Since $\Fmor{X}{\mdl}$ is dual-algebraic 
($\see$Proposition \ref{Fm(X,D)algebraic}), 
the category $\Mrt{\Fm{X}{\mdl}}{X}$ admits a universal object 
($\see$Theorem \ref{Exist univObj}), 
and this universal object is dual to 
$\big[ \Fmor{X}{\mdl} \lra \Picor{X} \big]$
($\see$Remark \ref{Alb_constr}). 
\ePf 

%\bRmk 
%For a rational map $\phe: X \dra G$ from $X$ to a smooth connected 
%algebraic group, 
%let  $\modun(\phe)$ be the modulus of $\phe$ defined by using the 
%filtration  $\fil_n \Witt_r$ instead of $\filf_n \Witt_r$. 
%Let $\Mrmn{X}{\mdl}$  be the category of those rational maps 
%$\phe$ such that  $\modun(\phe) \leq X$. 
%
%According to Theorem \ref{Mdl-category}, 
%the category $\Mrm{X}{\mdl}$ corresponds to the formal group 
%$\Fm{X}{\mdl}$. 
%The same argument as in Lemma \ref{mod(phe)-im(trafo_phe)}
%shows that the category $\Mrmn{X}{\mdl}$ is contained in 
%the category corresponding to $\Fmn{X}{\mdl}$. 
%We have only an implication in one direction in this case. 
%From the categorical point of view, this can be seen as follows: 
%In Theorem \ref{Exist univObj}
%it is shown that categories which satisfy the stability conditions 
%$\lrpn{\diamondsuit \; 1-3}$ on page~\pageref{Exist univObj} 
%correspond to formal subgroups $\sF$ of $\Divor{X}$. 
%Condition $\lrpn{\diamondsuit \; 2}$ says that the category needs to be 
%stable under certain homomorphisms between affine algebraic groups. 
%This condition is not satisfied by $\Mrmn{X}{\mdl}$, 
%because this category does not take the Frobenius into account. 
%$\Fmn{X}{\mdl}$ is the smallest formal group whose 
%associated category contains $\Mrmn{X}{\mdl}$. 
%This means $\Mr_{\Fmn{X}{\mdl}}$ 
%is the saturation of $\Mrmn{X}{\mdl}$ 
%with respect to these stability conditions. 
%\eRmk 

\subsubsection{Descent of the base field} 
\label{Alb(X,D)descent}

Let $\fld$ be a perfect field. 
Let $\clfld$ be an algebraic closure of $\fld$. 
Let $\Xo$ be a smooth proper variety defined over $\fld$, 
and let $\mdl$ be an effective divisor on $\Xo$ rational over $\fld$. 

\bThm 
\label{descentAlbm(X,D)}
There exists a $\fld$-torsor $\Albbm{1}{\Xo}{\mdl}$ for an algebraic 
$\fld$-group $\Albbm{0}{\Xo}{\mdl}$ and rational maps defined over $\fld$ 
\[ \albbm{i}{\Xo}{\mdl}: \Xo^{2-i} \dra \Albbm{i}{\Xo}{\mdl} 
\] 
for $i=1,0$, satisfying the following universal property: 

If $\phe: \Xo \dra \Trr{1}$ 
is a rational map defined over $\fld$ to a $\fld$-torsor $\Trr{1}$ 
for an algebraic $\fld$-group $\Trr{0}$, such that 
$\modu\lrpn{\phe} \leq \mdl$, 
then there exist a unique homomorphism of $\fld$-torsors 
$\homm{1}: \Albbm{1}{\Xo}{\mdl} \lra \Trr{1}$ 
and a unique homomorphism of algebraic $\fld$-groups 
$\homm{0}: \Albbm{0}{\Xo}{\mdl} \lra \Trr{0} \laurink$ 
defined over $\fld$, such that 
\bDpl \phe^{(i)} = \homm{i} \circ \albbm{i}{\Xo}{\mdl}
\eDpl 
for $i=1,0$. 

Here $\Albbm{0}{X}{\mdl}$ 
is dual to the 1-motive $\bigbt{\Fmor{X}{\mdl}\lra \Picor{X}}$. 
\eThm 

\bPf 
Follows directly from Thm.~\ref{descent} 
and the definition of the modulus via Galois descent (Def.~\ref{DefMod}). 
\ePf 

\bCor 
\label{every rat map}
For every rational map $\phe: X \dra P$ from $X$ to a torsor $P$ 
%smooth connected algebraic group $G$ 
there exists an effective divisor $\mdl$, 
namely $\mdl = \modu\lrpn{\phe}$, 
such that $\phe$ factors through $\Albbm{1}{X}{\mdl}$. 
\eCor 

\bPrp 
\label{exhaust Div_X^0}
Let $\fmlG$ be a formal $k$-subgroup of $\Divor{X}$. 
Then $\fmlG$ is dual-algebraic if and only if 
there exists an effective divisor $\mdl$, 
rational over $k$, 
such that $\fmlG \subset \Fm{X}{\mdl}$. 
\ePrp 

\bPf 
$(\Lla)$ A formal subgroup of a dual-algebraic group is also 
dual-algebraic, according to Lemma \ref{subFmlGroup_algebraic}. 

$(\Lra)$ 
By Galois descent %(cf.\ Remark.~\ref{fmlGroup_descent}) 
(possible for formal groups due to Cartier duality) 
we may assume that $k$ is algebraically closed. 
Let $\mdl = \modu\lrpn{\alb_{\fmlG}}$ be the modulus of the 
universal rational map \;$\alb_{\fmlG}: X \lra \AlbF{X}$\,
associated to $\fmlG \subset \Divor{X}$. 
Then by Lemma \ref{mod(phe)-im(trafo_phe)} it holds that 
$\fmlG = \im\lrpn{\trafo_{\alb_{\fmlG}}} \subset \Fm{X}{\mdl}$. 
\ePf

\subsubsection{Functoriality} 
\label{Alb(X,D)functorial}

We specialize the results from No.~\ref{Functoriality} 
to the case of Albanese varieties with modulus. 

\bPrp 
\label{alb_X,D(morphism)} 
Let $\morphism: Y \lra X$ be a morphism of smooth proper varieties. 
Let $\mdl$ be an effective divisor on $X$ 
intersecting $\morphism\pn{Y}$ properly. 
%with $Y$ decident ($\see$Definition \ref{decident}) to $\Supp(\mdl)$. 
Then $\morphism$ induces 
a homomorphism of torsors $\Albb{1}^{X,\mdl}_{Y,\mdll}\pn{\morphism}$ 
and a homomorphism of algebraic groups 
$\Albb{0}^{X,\mdl}_{Y,\mdll}\pn{\morphism}$, 
\[ \Albb{i}^{X,\mdl}_{Y,\mdll}\pn{\morphism}: 
   \Albbm{i}{Y}{\mdll} \lra \Albbm{i}{X}{\mdl} 
\] 
for each effective divisor $\mdll$ on $Y$ satisfying 
$\mdll \geq \lrpn{\mdl-\mdl_{\red}}\cut Y + \lrpn{\mdl\cut Y}_{\red}$, 
where $B \cut Y$ denotes the pull-back of a Cartier divisor $B$ on $X$ to $Y$. 
\ePrp 

\bPf 
According to Proposition \ref{alb_F(morphism)}, 
for the existence of $\Albmm{Y}{\mdll}{X}{\mdl}\pn{\morphism}$ 
it is sufficient to show $\Fm{Y}{\mdll} \supset \Fm{X}{\mdl} \cut Y$. 
Definition \ref{DefFm(X,D)} of $\Fm{X}{\mdl}$ implies that this is 
the case if and only if $\Supp\pn{\mdll} \supset \Supp\pn{\mdl\cut Y}$ 
and $\mdll-\mdll_{\red} \geq \lrpn{\mdl-\mdl_{\red}} \cut Y$. 
But this is equivalent to 
$\mdll \geq \lrpn{\mdl-\mdl_{\red}}\cut Y + \lrpn{\mdl\cut Y}_{\red}$. 
\ePf 

\bCor 
\label{Albm(E>D)}
If $\mdl$ and $\mdll$ are effective divisors on $X$ with $\mdll \geq \mdl$, 
then there are canonical surjective homomorphisms 
$ \Albbm{i}{X}{\mdll} \lsur \Albbm{i}{X}{\mdl} 
$ 
for $i = 1,0$, 
given by $\Albb{i}_{X,\mdll}^{X,\mdl}\pn{\id_X}$. 
\eCor 

\bPf 
If $\mdll \geq \mdl$, it is evident 
that $\Albbm{i}{X}{\mdll}$ generates $\Albbm{i}{X}{\mdl}$, 
thus $\Albb{i}_{X,\mdll}^{X,\mdl}\pn{\id_X}$ is surjective. 
\ePf

%\newpage 

\subsection{Jacobian with modulus} 
\label{subsec:JacMod}

Let $\Crv$ be a smooth proper curve over a perfect field $k$, 
which we assume to be algebraically closed for convenience. 
Let $\mdl = \sum_{\pnt\in S} n_{\pnt} \,\pnt$ be an effective divisor 
on $\Crv$, 
where $S$ is a finite set of closed points on $\Crv$ 
and $n_{\pnt}$ are integers $\geq 1$ for $\pnt\in S$. 
The Jacobian $\Jacm{\Crv}{\mdl}$ of $\Crv$ of modulus $\mdl$ 
is by definition the universal object for the category of those 
morphisms $\phe$ from $\Crv\setminus S$ to algebraic groups 
such that $\phe\big(\dv(f)\big) = 0$ for all $f\in\sK_{\Crv}$ 
with $f\equiv 1 \mod\mdl$. 
Here we used the definition 
$\phe\big(\sum l_j c_j\big) = \sum l_j \,\phe(c_j)$ 
for a divisor $\sum l_j c_j$ on $\Crv$ with $c_j\in\Crv\setminus S$, 
and ``$f\equiv 1 \mod\mdl$'' means $\val_{\pnt}(1-f) \geq n_{\pnt}$ 
for all $\pnt\in S$, 
where $\val_{\pnt}$ is the valuation attached to the point $\pnt\in\Crv$. 

\bThm 
\label{StructureJacm}
The generalized Jacobian $\Jacm{\Crv}{\mdl}$ of $\Crv$ 
of modulus $\mdl$ is an extension 
\[ 0 \lra \Lm{\Crv}{\mdl} \lra \Jacm{\Crv}{\mdl} \lra \Jacc{\Crv} \lra 0 
\] 
of the classical Jacobian $\Jacc{\Crv} \cong \Pic^0_{\Crv}$ of $\Crv$, 
which is an abelian variety,  
by the affine algebraic group $\Lm{\Crv}{\mdl}$, 
which is characterized by 
\[ \Lm{\Crv}{\mdl}(k) =  \frac{\prod_{\pnt\in S} k(\pnt)^*}{k^*} \tms 
       \prod_{\pnt\in S} \frac{1+\fm_\pnt}{1+\fm_\pnt^{n_\pnt}} 
\] 
where $k(\pnt)$ denotes the residue field 
and $\fm_{\pnt}$ the maximal ideal at $\pnt\in\Crv$. 
%\[ \Lm{\Crv}{\mdl} = 
%   \frac{\prod_{\pnt\in S}\GmS{\pnt}}{\Gm} \tms 
%   \prod_{\pnt\in S} \Upf_{\sO_{\Crv,\pnt}\left/\fm_{\Crv,\pnt}^{n_{\pnt}}\right.}
%   %\Upf_{k[[t_{\pnt}]]\left/\lrpn{t_{\pnt}^{n_{\pnt}}}\right.} 
%\] 
%where $\GmS{\pnt}$ denotes a group isomorphic 
%to the multiplicative group and attached to the point $\pnt$, 
%where $\Upf_R = \ker \bigpn{\Lin_R \lra \Lin_{R_{\red}}}$ 
%is the unipotent group associated 
%to a finite $k$-algebra $R$ 
%(here $R_{\red} = R / \Nil(R)$), 
%%from \ref {LinGroup_Ring}, 
%%($\see$Lemma \ref{Vnipot}), 
%and $\lrpn{\sO_{\Crv,\pnt},\fm_{\Crv,\pnt}}$ is the local ring at $\pnt\in S$. 
\eThm 

\bPf 
\cite[V, \S~3]{S}, see also the summary \cite[I, No.~1]{S}. 
\ePf 

%\vspace{\vs} 
%
%We give an illustration of Theorem \ref{StructureJacm}, 
%cf. \cite[I, No.~1]{S}. 
%The classical Jacobian $\Jacc{\Crv}$ is 
%the group of divisor classes on $\Crv$ of degree 0. 
%The Jacobian $\Jacm{\Crv}{\mdl}$ of modulus $\mdl$ is identified to 
%the group of classes of divisors prime to $S$ 
%modulo principal divisors $\dv\pn{f}$ with $f \equiv 1 \mod \mdl$. 
%Then the affine part $\Lm{\Crv}{\mdl}$ of $\Jacm{\Crv}{\mdl}$, 
%as the kernel of the canonical surjection 
%$\Jacm{\Crv}{\mdl} \lsur \Jacc{\Crv}$, 
%is characterized by %the following group of $k$-valued points 
%
%\begin{eqnarray*} 
%\Lm{\Crv}{\mdl}\pn{k} 
% & = & \frac{\left\{\dv(f)\;|\;f\in\sO_{C,\pnt}^* 
%                            \quad\forall \pnt\in S\right\}}
%            {\left\{\dv(f)\;|\;\val_{\pnt}\lrpn{1-f}\geq n_{\pnt}
%                            \quad\forall \pnt\in S\right\}} \\ 
% & = & \frac{\left\{f\in\sK_{C}\;|\;f\in\sO_{C,\pnt}^* 
%                            \quad\forall \pnt\in S\right\}}
%            {k^* \tms \left\{f\in\sK_{C}\;|\;\val_{\pnt}\lrpn{1-f}\geq n_{\pnt}
%                            \quad\forall \pnt\in S\right\}} \\ 
% & = & \frac{\prod_{\pnt\in S} \sO_{C,\pnt}^*}
%            {k^* \tms \prod_{\pnt\in S}\lrpn{1+\fm_\pnt^{n_\pnt}}} \\ 
% & = & \frac{\prod_{\pnt\in S} k(\pnt)^*}{k^*} \tms 
%       \prod_{\pnt\in S} \frac{1+\fm_\pnt}{1+\fm_\pnt^{n_\pnt}} 
%\end{eqnarray*}
%where $k(\pnt)$ denotes the residue field 
%and $\fm_{\pnt}$ the maximal ideal at $\pnt\in\Crv$. 

\bThm 
\label{dual_Jacm}
The Jacobian with modulus $\Jacm{\Crv}{\mdl}$ is dual 
(in the sense of 1-motives) to the 1-motive 
$ \left[\Fmo{\Crv}{\mdl} \lra \Pic^0_{\Crv}\right]
$, 
where $\Fmo{\Crv}{\mdl} = \Fmor{\Crv}{\mdl}$ 
is the formal subgroup of $\Divf^0_{\Crv}$ 
from Definition \ref{DefFm(X,D)}, 
and $\Fmo{\Crv}{\mdl} \lra \Pic^0_{\Crv}$ is the homomorphism 
induced by the class map $\Divf^0_{\Crv} \lra \Picf^0_{\Crv}$. 
\eThm 

\bPf 
We have to ensure that the category 
for which $\Jacm{\Crv}{\mdl}$ is universal 
is characterized by the formal group $\Fm{\Crv}{\mdl}$.  
%\[ \Mrm{\Crv}{\mdl} = 
%   \{ \phe:\Crv\setminus S \lra G \;\left|\; 
%   \phe\bigpn{\dv(f)}=0 \;\textrm{ if }\; f\equiv 1\mod\mdl \right.\} 
%\] 
The Jacobian $\Jacm{\Crv}{\mdl}$ of modulus $\mdl$ 
is by definition the universal object for morphisms 
$\phe$ from $\Crv\setminus S$ to algebraic groups 
satisfying 

\begin{tabular}{rl} 
(i) & $\phe\bigpn{\dv(f)}=0$  \hspace{40pt}
$\forall f \in\sK_{\Crv}$ with $f\equiv 1\mod\mdl$\laurin 
\end{tabular} 

\noindent 
Condition (i) is equivalent to 

\begin{tabular}{rl} 
(ii) & $\lrpn{\phe,f}_\pnt = 0$  \hspace{40pt} 
$\forall \pnt \in S,\, \forall f \in \sK_{\Crv}$ with $f\equiv 1\mod\mdl$ at $q$
\end{tabular} 

\noindent 
where $\lrpn{\phe,\lull}_{\lul}: \sK_{\Crv}^* \tms \Crv \lra G(k)$ 
is the local symbol associated to the morphism 
$\phe:\Crv\setminus S\lra G$, 
according to \cite[I, No.~1, Thm.~1 and III, \S~1]{S}. 
It is shown in \cite[No.~6.1--6.3]{KR2} that condition (ii) 
is equivalent to 

\begin{tabular}{rl} 
(iii) & $\modu\lrpn{\phe} \leq \mdl$\laurin 
\end{tabular} 

\noindent 
Then the assertion is due to 
Theorems \ref{Mdl-category} and \ref{AlbMod-construction}. 
\ePf

%\newpage 

\subsection{Relative Chow group with modulus} 
\label{subsec:ChowMod}

Let $X$ be a smooth proper variety over an algebraically closed field $k$, 
let $\mdl$ be an effective divisor on $X$ 
and $\mdl_{\red}$ the reduced part of $\mdl$. 

\bNot 
If $C$ is a curve in $X$, then $\nu:\Ct\lra C$ denotes the normalization. 
For $f \in \sK_C$, we write $\wt{f} := \nu^{*}f$ for the image of $f$ in $\sK_{\Ct}$. 
If $\phe: X \dra G$ is a rational map, we write $\phe|_{\Ct} := \phe|_C \circ \nu$ 
for the composition of $\phe$ and $\nu$. 
If $B$ is a Cartier divisor on $X$ intersecting $C$ properly, 
then $B\cut\Ct$ denotes the pull-back of $B$ to $\Ct$. 
\eNot 

\bDef 
\label{DefCHoMod}
Let $\Z_{0}\lrpn{X\setminus\mdl}$ 
be the group of 0-cycles on $X\setminus\mdl$, set 
\[ \Rlm{X}{\mdl} = \left\{ \lrpn{C,f}\left|\begin{array}{l}
C\textrm{ a curve in $X$ intersecting } \Supp\pn{\mdl} \textrm{ properly, } \\ 
f \in \sK_C^* \textrm{ s.t.\ } 
\wt{f} \equiv 1 \mod \lrpn{\mdl-\mdl_{\red}}\cut \Ct + \bigpn{\mdl\cut\Ct}_{\red} 
\end{array}\right.\right\} %\laurink 
\] 
and let $\R_{0}\lrpn{X,\mdl}$ be the subgroup of 
$\Z_{0}\lrpn{X\setminus\mdl}$ 
generated by the elements $\dv\pn{f}_{C}$ with 
$\lrpn{C,f}\in\Rlm{X}{\mdl}$.
Then define 
\[ \CHm{X}{\mdl} = \Z_{0}(X\setminus\mdl) \big/ \R_{0}(X,\mdl) \laurin 
\] 
Let $\CHmo{X}{\mdl}$ be the subgroup of $\CHm{X}{\mdl}$
of cycles $\zeta$ with $\deg\zeta|_{W}=0$ for all irreducible components
$W$ of $X\setminus\mdl$. %$W\in\Cp\lrpn{X\setminus\mdl}$ 
\eDef  

\bDef 
\label{MrCH(X,D)} 
Let $\MrCHm{X}{\mdl}$ be the category of rational maps 
from $X$ to algebraic groups defined as follows: 
the objects of $\MrCHm{X}{\mdl}$ are morphisms 
$\phe:X\setminus\mdl \lra G$ 
whose associated map 
on 0-cycles of degree zero 
\bDpl 
  \Z_0(X\setminus\mdl)^0  \lra  G(k) \laurink \hspace{1mm} 
  \sum l_i \, p_i  \lmt  \sum l_i \, \phe(p_i) \laurink 
  l_i \in \Zint \laurink 
\eDpl 
factors through a homomorphism of groups $\CHmo{X}{\mdl}\lra G(k)$. 
The morphisms are the ones as in Definition \ref{CatMr}. 
%\footnote{A category of rational maps to algebraic groups is defined 
%already by its objects.%, according to Remark \ref{EquCatMr}. 
%} 
We refer to the objects of $\MrCHm{X}{\mdl}$ as 
rational maps from $X$ to algebraic groups 
\emph{factoring through $\CHmo{X}{\mdl}$}. 
\eDef 

\bThm 
\label{ThmCHoMod}
The category $\Mrm{X}{\mdl}$ of rational maps of modulus $\leq\mdl$ 
is equivalent to the category $\MrCHm{X}{\mdl}$ of rational maps 
factoring through $\CHmo{X}{\mdl}$. 
In particular, the Albanese $\Albm{X}{\mdl}$ of $X$ of modulus $\mdl$ 
is the universal quotient of $\CHmo{X}{\mdl}$. 
\eThm 

\bPf 
According to the definitions of $\Mrm{X}{\mdl}$ and $\MrCHm{X}{\mdl}$ 
the task is to show that for a morphism 
$\phe:X\setminus\mdl \lra G$ from $X\setminus\mdl$ to a smooth 
connected algebraic group $G$ the following conditions are equivalent: 

\begin{tabular}{rl} 
\hspace{\kukwad} (i) & $\modu\lrpn{\phe} \leq \mdl$\laurink \\ 
\hspace{3pt} (ii) & $\phe\lrpn{\dv\pn{f}_C} = 0 
\hspace{50pt} \forall\, \lrpn{C,f} \in \Rlm{X}{\mdl}$\laurin 
\end{tabular} 

\noindent 
Since $\phe\lrpn{\dv\pn{f}_C} = \phe|_{\Ct} \bigpn{\dv\bigpn{\wt{f}}_{\Ct}}$ 
($\seecite$\cite[Lemma~3.32]{Ru1}), 
condition (ii) is equivalent to the condition 

\begin{tabular}{rl} 
\hspace{\kukwad} (iii) & $\modu\lrpn{\phe|_{\Ct}} \leq 
\lrpn{\mdl-\mdl_{\red}}\cut\Ct + \bigpn{\mdl\cut\Ct}_{\red}$ \\ 
 & for all curves $C$ in $X$ intersecting $\Supp\pn{\mdl}$ properly\laurink 
\end{tabular} 

\noindent 
as was seen in the proof of Theorem \ref{dual_Jacm}, 
substituting $D$ by $\lrpn{\mdl-\mdl_{\red}}\cut\Ct + \bigpn{\mdl\cut\Ct}_{\red}$. 
The equivalence of (i) and (iii) is the content of 
Lemma \ref{Mdl_restr_to_curves}. 
\ePf 

\bLem 
\label{Mdl_restr_to_curves}
Let $\phe:X\dra G$ be a rational map from $X$ 
to a smooth connected algebraic group $G$. 
Then the following conditions are equivalent: \\ 
\begin{tabular}{rl} 
\hspace{\kukwad} (i) & $\modu\lrpn{\phe} \leq \mdl$\laurink \\ 
\hspace{\kukwad} (ii) & $\modu\lrpn{\phe|_{\Ct}} \leq 
\lrpn{\mdl-\mdl_{\red}}\cut\Ct + \bigpn{\mdl\cut\Ct}_{\red}$ \\ 
 & for all curves $C$ in $X$ intersecting $\Supp\pn{\mdl}$ properly\laurin 
%\hspace{\kukwad} (ii) & $\modu\lrpn{\phe|_C} \leq 
%\lrpn{\mdl-\mdl_{\red}}\cut C + \lrpn{\mdl\cut C}_{\red} 
%\hspace{13pt} \forall\, C \textrm{ decident to }\Supp\pn{\mdl}$\laurin 
\end{tabular}
\eLem 

\bPf 
(i)$\Lra$(ii)
Let $C$ be a curve in $X$ intersecting $\mdl$ properly. 
As $\phe$ is regular away from $\mdl$, 
the restriction $\phe|_{\Ct}$ of $\phe$ to $\Ct$ 
is regular away from $\mdl \cut\Ct$. 
Hence 
$\Supp\bigpn{\modu\lrpn{\phe|_{\Ct}}} \subset \Supp\bigpn{\mdl\cut\Ct} 
  = \Supp\bigpn{\pn{\mdl-\mdl_{\red}}\cut\Ct + \bigpn{\mdl\cut\Ct}_{\red}}$.  
%If a divisor $\Es$ is locally given by a function $f \in \sK_{X}^*$, 
%then by definition $\Es \cut\Ct$ is locally given by $\wt{f|_C} \in \sK_{\Ct}^*$. 
%With this, 
According to Definition \ref{DefMod} of the modulus, 
it is easy to see that 
$\modu\lrpn{\phe} \leq \mdl = \lrpn{\mdl-\mdl_{\red}} + \mdl_{\red}$ 
implies 
\;$\modu\lrpn{\phe|_{\Ct}} \leq 
\lrpn{\mdl-\mdl_{\red}}\cut\Ct + \bigpn{\mdl\cut\Ct}_{\red}$. 

\medskip 
(ii)$\Lra$(i) 
Let $\mdll := \modu\lrpn{\phe}$ 
and $\pnt\in\Supp\pn{\mdll}$ be a point of codimension 1 in $X$. 
We are going to construct a family of smooth curves $\lrst{C_{e}}_{e}$ 
%with the property of $\lrpn{\dagger }$. 
intersecting $\mdll$ in a fixed point $x \in \Es_{\pnt} = \ol{\lrst{\pnt}}$ 
such that 
\[ \lim_{e\ra\infty} \frac{\modu_{x}\lrpn{\phe|_{C_{e}}}}
{\multy_{x}\bigpn{\lrpn{\mdll-\mdll_{\red}} \cut C_{e}} + 1} = 1
\] 
where $\multy_{x}\lrpn{\mdll \cut \Crv}$ denotes the intersection multiplicity 
of $\mdll$ and $\Crv$ at $x$. 

After the construction we will show that 
the existence of such a family of curves for each 
$\pnt\in\Supp\pn{\mdll}$ of codimension 1 in $X$ 
yields the implication (ii)$\Lra$(i). 

If $\chr(k) = 0$, it is easy to see that a general curve $C$ in $X$ 
intersecting $\Es_{\pnt}$ in a point $x$ satisfies 
$\modu_{x}\lrpn{\phe|_C} = 
  \multy_{x}\bigpn{\lrpn{\mdll-\mdll_{\red}}\cut C} + 1$. 
Therefore we suppose that $\chr(k)=p>0$. 
Using the notation of Definition \ref{DefMod}, let 
%$\lrpn{u_{\pnt,i}}_{1\leq i\leq a} \in \Ga\pn{\sK_{X,\pnt}}^a$ resp. 
$\lrpn{u_{\pnt,i}}_{1\leq i\leq a} \in \Witt_r\pn{\sK_{X,\pnt}}^a$
be a representative of the unipotent part of the class of 
$\phe\in G\pn{\sK_{X,\pnt}}$ in $G\pn{\sK_{X,\pnt}}/G\pn{\sO_{X,\pnt}} 
= L\pn{\sK_{X,\pnt}}/L\pn{\sO_{X,\pnt}}$. 
Then $\modu_{\pnt}\pn{\phe} = 1 + \nty_{\pnt}\pn{u_{\pnt,i}}$ 
for some $1 \leq i \leq a$. 
Set $n := \nty_{\pnt}\pn{u_{\pnt,i}}$. 
Let $t \in \fm_{X,\pnt}$ be a uniformizer at $\pnt$. 
Let $\sum_{\nu} \Frob^{\nu} \tens \,\oma_{\nu} \tens t^{-n} \in 
k[\Frob] \tens_k \Kot_{X,\pnt}\pn{\log\pnt} \tens_{\sO_{X,\pnt}} \fm_{X,\pnt}^{-n}$ 
be a representative of $\qfill_{n\pnt} \pn{u_{\pnt,i}} \in \Qfill_{n\pnt}$ 
($\see$Definition \ref{fil-quot}). 
Choose a regular closed point $x \in \Es_{\pnt}$ 
such that $t$ is a local equation for $\Es_{\pnt}$ at $x$
and $\oma_{\nu}$ is regular and $\neq 0$ at $x$ for some $\nu$. 
We may assume that $\dim X = 2$ 
via cutting down by hyperplanes through $x$ transversal to $\Es_{\pnt}$.  
Let $s \in \fm_{X,x}$ be a local parameter at $x$ 
that gives a uniformizer of $\sO_{\Es_{\pnt},x}$. 
Define a curve $C_{e}$ locally around $x$ 
by the equation \;$t = s^{e}$ \; for $e \geq 1$. 
Note that $\mdll-\mdll_{\red}$ is locally defined by the equation $t^{n} = 0$. 
Then 
\[ \multy_{x}\bigpn{\lrpn{\mdll-\mdll_{\red}}\cut C_{e}} \,=\, 
   \dim_k \frac{\sO_{X,x}}{\pn{t^{n}, t - s^e}} \,=\, ne \laurin 
\]
We can write \;$\oma_{\nu} = g\, \der s + h\, \der\log t$\; 
with $g,h \in \sO_{X,\pnt}$ 
and the values at $x$ are 
$g(x) \neq 0$ if $\qfill_{n\pnt} \pn{u_{\pnt,i}} \in \Qfillb_{n\pnt}$, 
and $h(x) \neq 0$ if $\qfill_{n\pnt} \pn{u_{\pnt,i}} \in 
 \Qfill_{n\pnt} \setminus \Qfillb_{n\pnt}$ 
and $x$ in general position (what we assume), 
for some $\nu$. 
The restriction of $t^{-n} \oma_{\nu}$ to $C_{e}$ is 
\begin{eqnarray*} 
t^{-n} \oma_{\nu}|_{C_{e}} 
& = & s^{-ne} g\, \der s +  s^{-ne} h\, \der\log s^e \\ 
& = & s^{1-ne} g\, \der\log s + e \,s^{-ne} h\, \der\log s \laurink 
\end{eqnarray*} 
and the class of $t^{-n} \oma_{\nu}|_{C_{e}}$ is $\neq 0$ in 
%with $0 \neq \lrbt{t^{-n} \oma_{\nu}|_{C_{e}}}$ 
\[ %\in
\left\{ 
\begin{array}{cl}
\Kot_{C_{e},x}(\log x) \tens_{\sO_{C_{e},x}} 
\fm_{C_{e},x}^{-ne} \big/ \fm_{C_{e},x}^{1-ne}
 & \; \textrm{ if }\; \qfill_{n\pnt} \pn{u_{\pnt,i}} \in 
 \Qfill_{n\pnt} \setminus \Qfillb_{n\pnt} \textrm{ and } p \nmid e \\
\Kot_{C_{e},x}(\log x) \tens_{\sO_{C_{e},x}} 
\fm_{C_{e},x}^{1-ne} \big/ \fm_{C_{e},x}^{2-ne}
 & \; \textrm{ if }\; \qfill_{n\pnt} \pn{u_{\pnt,i}} \in \Qfillb_{n\pnt} \laurin 
\end{array}
\right. 
\] 
Lemma \ref{fil-restriction} assures that the modulus of $\phe|_{C_{e}}$ 
is computed from the restriction (of a representative) of 
$\qfill_{n\pnt} \pn{u_{\pnt,i}}$ to $C_{e}$, 
for $e$ large enough such that $ne-1 > \lfloor ne/p\rfloor$ 
(this is satisfied for $e > 2$). 
Thus we have 
\[ \nty_{x} \lrpn{u_{\pnt,i}|_{C_{e}}} = 
\left\{ 
\begin{array}{cl}
ne & \textrm{ if }\; \qfill_{n\pnt} \pn{u_{\pnt,i}} \in 
 \Qfill_{n\pnt} \setminus \Qfillb_{n\pnt} \textrm{ and } p \nmid e \\
ne-1 & \textrm{ if }\; \qfill_{n\pnt} \pn{u_{\pnt,i}} \in \Qfillb_{n\pnt} \laurink 
\end{array}
\right. 
\]  
\[ \modu_{x} \lrpn{\phe|_{C_{e}}} = 
\left\{ 
\begin{array}{cl}
ne+1 & \textrm{ if }\; \qfill_{n\pnt} \pn{u_{\pnt,i}} \in 
 \Qfill_{n\pnt} \setminus \Qfillb_{n\pnt} \textrm{ and } p \nmid e \\
ne & \textrm{ if }\; \qfill_{n\pnt} \pn{u_{\pnt,i}} \in \Qfillb_{n\pnt} \laurin 
\end{array}
\right. 
\]  
Then 
\[ \lim_{\substack{e \ra \infty \\ p \nmid e}} \frac{\modu_{x}\lrpn{\phe|_{C_{e}}}}
{\multy_{x}\bigpn{\lrpn{\mdll-\mdll_{\red}} \cut C_{e}} + 1} = 1 
\laurin 
\] 

\medskip 
Now we show $\neg (i) \Lra \neg (ii)$. 
Suppose $\mdll := \modu\lrpn{\phe} \not \leq \mdl$. 
Then there is a point $\pnt \in \Supp\pn{\mdll}$ of codimension 1 in $X$ 
such that $\multy_{\pnt}(\mdll) > \multy_{\pnt}(\mdl)$, 
where $\multy_{\pnt}$ is the multiplicity at $\pnt$. 
By the construction above there is a sequence of curves 
$\lrst{C_{e}}_{e}$ in $X$ intersecting $\mdll$ in a fixed point 
$x \in \Es_{\pnt}$ such that 
\[ \underset{\substack{e \ra \infty \\ p \nmid e}} \lim 
\frac{\modu_{x}\lrpn{\phe|_{C_{e}}}}
{\multy_{x}\lrpn{\lrpn{\mdll-\mdll_{\red}} \cut C_{e}} + 1} = 1 
\laurin 
\] 
If $\multy_{\pnt}(\mdl) \neq 0$, then since 
\[ \underset{e\geq 0} \sup 
\frac{\multy_{x}\lrpn{\lrpn{\mdl-\mdl_{\red}} \cut C_{e}} + 1}
{\multy_{x}\lrpn{\lrpn{\mdll-\mdll_{\red}} \cut C_{e}} + 1} < 1 
\] 
there is %a sufficiently large 
$e$ such that 
\;$\modu_{x}\lrpn{\phe|_{C_{e}}} > 
\multy_{x}\bigpn{\lrpn{\mdl-\mdl_{\red}} \cut C_{e}} + 1$. 
If $\multy_{\pnt}(\mdl) = 0$, then 
$0 \neq \modu\lrpn{\phe|_{C_{e}}}_{x} > 
\multy_{x}\bigpn{\lrpn{\mdl-\mdl_{\red}}\cut C_{e} + 
\lrpn{\mdl\cut C_{e}}_{\red}} = 0$. 
Thus \;$\modu\lrpn{\phe|_{C_{e}}} \not \leq 
\lrpn{\mdl-\mdl_{\red}}\cut C_{e} + \lrpn{\mdl\cut C_{e}}_{\red}$. 
\ePf

{\scshape
\begin{flushright}
\begin{tabular}{l}
Freie Universit\"at Berlin, Germany \\ 
Institut f\"ur Mathematik \\ 
{\upshape e-mail: \texttt{henrik.russell@math.fu-berlin.de}}\\
\end{tabular}
\end{flushright}
}


\newpage 

\begin{thebibliography}{KS}
%\addcontentsline{toc}{section}{References} 

%\bibitem[At]{A} M. F. Atiyah, I. G. MacDonald, 
%\emph{Introduction to Commutative Algebra}, 
%Addison-Wesley 1969. 

%\bibitem[BaS]{BS} L. Barbieri-Viale, V. Srinivas, 
%\emph{Albanese and Picard 1-motives},
%M\'emoire SMF \textbf{87}, Paris, 2001. 

\bibitem[BB]{BaBe} L. Barbieri-Viale, A. Bertapelle, 
\emph{Sharp de Rham realization}, 
Advances in Mathematics \textbf{222} (2009), pp.~1308--1338.  

\bibitem[Ba]{B}I. Barsotti, \emph{Structure theorems for group varieties}, 
Annali di Matematica pura et applicata, Serie IV, T. \textbf{38} (1955),
pp.~77--119. 

%\bibitem[BiS]{BiS} J. Biswas, V. Srinivas, 
%\emph{Roitman's theorem for singular projective varieties}, 
%Compositio Mathematica \textbf{119} (1999), No.~2, pp.~213--237. 

\bibitem[BLR]{BLR}S. Bosch, W. L\"utkebohmert, M. Raynaud, 
\emph{N\'eron models}, 
Ergebnisse der Mathematik \textbf{21}, 3. Folge, Springer-Verlag 1990. 

\bibitem[Br]{Br} J.-L. Brylinski, 
\emph{Th\'eorie du corps de classes de Kato et rev\^etements ab\'eliens 
de surfaces}, 
Annales de l'Institut Fourier, tome \textbf{33}, No.~3 (1983), 
pp.~23--38. 

\bibitem[Ch]{C} C. Chevalley, \emph{La th\'eorie de groupes alg\'ebriques}, 
Proceedings of the International Congress of Mathematicians 1958, 
Cambridge University Press, Cambridge 1960, pp.~53--68. 

\bibitem[Dl]{D2} P. Deligne, \emph{Th\'eorie de Hodge II} and \emph{III}, 
Publications Math\'ematiques de l'IH\'ES 
\textbf{40} (1971), pp.~5--57 and \textbf{44} (1974), pp.~5--78. 

%\bibitem[EGA2]{EGA2}J. Dieudonne, A. Grothendieck, 
%\emph{El\'ements de G\'eom\'etrie Alg\'ebrique II}, 
%Publications Math\'ematiques de l'IH\'ES \textbf{8} (1961), pp.~5--222. 

\bibitem[Dm]{D} M. Demazure, \emph{Lectures on p-Divisible Groups}, 
Lecture Notes in Mathematics \textbf{302}, Springer-Verlag 1972. 

\bibitem[DG]{DG} M. Demazure, P. Gabriel, \emph{Groupes Alg\'ebriques, 
Tome I: G\'eometrie alg\'ebrique - g\'en\'eralities, Groupes commutatifs}, 
North-Holland Publishing Company - Amsterdam 1970. 

\bibitem[ESV]{ESV} H. Esnault, V. Srinivas, E. Viehweg, \emph{The universal 
regular quotient of the Chow group of points on projective varieties}, 
Inventiones Mathematicae \textbf{135} (1999), pp.~595--664. 

%\bibitem[FGA]{FGA} A. Grothendieck, 
%\emph{Fondements de la G\'eom\'etrie Alg\'ebrique}, 
%[Extraits du S\'eminaire Bourbaki 1957-1962], S\'ecr\'etariat
%math\'ematique, 11 rue Pierre Curie, Paris 5e, 1962. 

\bibitem[FW]{FW} G. Faltings, G. W\"ustholz, \emph{Einbettungen 
kommutativer algebraischer Gruppen und einige ihrer Eigenschaften}, 
Journal f\"ur die reine und angewandte Mathematik \textbf{354} (1984), 
pp.~175--205. 

%\bibitem[Fl]{Fl} H. Flenner, 
%\emph{Die S\"atze von Bertini f\"ur lokale Ringe}, 
%Mathematische Annalen \textbf{229} (1977), pp.~97--111. 

\bibitem[Fo]{Fo} J.-M. Fontaine, 
\emph{groupes p-divisibles sur les corps locaux}, 
Ast\'erisque \textbf{47--48} (1977). 

%\bibitem[Fr]{Fr} P. Freyd, 
%\emph{Abelian Categories, an Introduction to the Theory of Functors}, 
%Harper \& Row, New York 1964. 

%\bibitem[Fu]{F} W. Fulton, \emph{Intersection Theory}, 
%Ergebnisse der Mathematik und ihrer Grenzgebiete, 
%3. Folge, Band 2, Springer-Verlag 1984. 

%\bibitem[GM]{GM} G. van der Geer, B. Moonen, \emph{Abelian Varieties}, 
%online script (2006), 
%\texttt{http://staff.science.uva.nl/\~{}bmoonen/boek/BookAV.html}. 

%\bibitem[GZ]{GZ} P. Gabriel, M. Zisman, 
%\emph{Calculus of fractions and homotopy theory}, 
%Ergebnisse der Mathematik und ihrer Grenzgebiete, 
%Springer-Verlag 1967. 

%\bibitem[Ha]{H} R. Hartshorne, \emph{Algebraic Geometry}, 
%Graduate Texts in Mathematics \textbf{52}, Springer Verlag 1977. 

%\bibitem[Jo]{J} A. J. de Jong, 
%\emph{Smoothness, semistability and alterations}, 
%Publications Math\'ematiques de l'IH\'ES \textbf{83} (1996), pp.~51--93. 

%\bibitem[KwS]{KwS} M. Kashiwara, P. Schapira, \emph{Categories and sheaves}, 
%Grundlehren der mathematischen Wissenschaften \textbf{332}, 
%Springer Verlag 2006. 

\bibitem[KR1]{KR1} K. Kato, H. Russell, 
\emph{Albanese varieties with modulus and Hodge theory}, 
Annales de l'Institut Fourier  \textbf{62}, No.~2 (2012), pp.~783--806. 

\bibitem[KR2]{KR2} K. Kato, H. Russell, 
\emph{Modulus of a rational map into a commutative algebraic group}, 
Kyoto Journal of Mathematics \textbf{50}, No.~3 (2010), 
pp.~607--622. 

\bibitem[KS]{KS} K. Kato, S. Saito, 
\emph{Two Dimensional Class Field Theory}, 
in ``Galois Groups and their Representations'', 
Advanced Studies in Pure Mathematics \textbf{2} (1983), pp.~103--152. 

\bibitem[Kl]{K} S. L. Kleiman, \emph{The Picard scheme}, 
preprint alg-geom/0504020v1 (2005). 

\bibitem[Lg]{La} S. Lang, \emph{Abelian Varieties}, 
Interscience Publisher, New York 1959. 

\bibitem[Ln]{L} G. Laumon, \emph{Transformation de Fourier g\'en\'eralis\'ee},
preprint alg-geom/9603004 (1996). 

%\bibitem[LW]{LW} M. Levine, C. A. Weibel, 
%\emph{Zero cycles and complete intersections on singular varieties}, 
%Journal f\"ur die reine und angewandte Mathematik \textbf{359} (1985), 106--120. 

%\bibitem[Mmu]{Mm} H. Matsumura, \emph{Commutative ring theory}, 
%Cambridge University Press 1986. 

\bibitem[Mi]{Mi} J. S. Milne, 
\emph{\'Etale Cohomology}, 
Princeton University Press 1980. 

\bibitem[Ms]{Ms} T. Matsusaka, 
\emph{On the algebraic construction of the Picard variety II}, 
Japanese Journal of Mathematics \textbf{22} (1952), pp.~51--62. 

%\bibitem[Mu1]{M2} D. Mumford, \emph{Lectures on Curves on an Algebraic Surface}, 
%Annals of Mathematics Studies \textbf{59}, Princeton University Press 1966
%\bibitem[Mu2]{M} D. Mumford, \emph{Abelian varieties}, 
%Tata Institute of Fundamental Research, Bombay, 1970, 
%and Oxford University Press 1985. 

\bibitem[\"On]{On} H. \"Onsiper, 
\emph{Generalized Albanese varieties for surfaces in characteristic $p>0$}, 
Duke Mathematical Journal, Vol.~\textbf{59}, No.~2 (1989), 
pp.~359--364. 

\bibitem[Oo]{O} F. Oort, \emph{Commutative group schemes}, 
Lecture Notes in Mathematics \textbf{15}, Springer-Verlag 1966. 

%\bibitem[Po]{P} A. Polishchuk, 
%\emph{Abelian varieties, theta functions and the Fourier transform}, 
%Cambridge University Press 2003. 

%\bibitem[Ra1]{Ra1} N. Ramachandran, 
%\emph{Duality of Albanese and Picard 1-motives},
%$K$-Theory \textbf{22} (2001), pp.~271--301. 

%\bibitem[Ra2]{Ra2} N. Ramachandran, 
%\emph{One-motives and a conjecture of Deligne},
%Journal of Algebraic Geometry \textbf{13}, No.~1 (2004), pp.~29--80. 

\bibitem[Ru]{Ru1} H. Russell, 
\emph{Generalized Albanese and its dual}, 
Journal of Mathematics of Kyoto University, Vol.~\textbf{48}, No.~4 (2008), pp.~907--949. 

%\bibitem[Ru2]{Ru2} H. Russell, 
%\emph{1-Motives with Unipotent Part}, 
%preprint arXiv:0902.2533 (2009). 

\bibitem[Ro]{R} M. Rosenlicht, 
\emph{Some basic theorems on algebraic groups}, 
American  Journal of Mathematics, Vol. \textbf{78}, No.~2 (1956), pp.~401--443. 

\bibitem[Se1]{S2} J.-P. Serre, 
\emph{Morphisme universels et vari\'et\'e d'Albanese}, 
in ``Vari\'et\'es de Picard'' ENS, S\'eminaire C. Chevalley, No.~10 (1958/59). 

\bibitem[Se2]{S3} J.-P. Serre, 
\emph{Morphisme universels et diff\'erentielles de troisi\`eme esp\`ece,} 
in ``Vari\'et\'es de Picard'' ENS, S\'eminaire C. Chevalley, No.~11 (1958/59). 

\bibitem[Se3]{S} J.-P. Serre, \emph{Groupes alg\'ebriques et corps de classes}, 
Hermann 1959. 

%\bibitem[Se4]{Se} J.-P. Serre, \emph{Corps locaux}, 
%Hermann 1966. 

%\bibitem[SC]{SC} S\'eminaire Chevalley, \emph{Groupes de Lie alg\'ebriques}, 
%1956/58. 

\bibitem[SGA3]{SGA3} M. Demazure, A. Grothendieck, 
\emph{Propri\'et\'es G\'en\'erales des Sch\'emas en Groupes}, 
SGA3, Lecture Notes in Mathematics \textbf{151},
Springer-Verlag 1959. 

\bibitem[Wa]{W} W. C. Waterhouse, \emph{Introduction to Affine Group Schemes}, 
Springer-Verlag New York 1979. 

\bibitem[We]{Wei} C. A. Weibel, 
\emph{An Introduction to Homological Algebra}, 
Cambridge University Press 1994. 

\end{thebibliography}
\end{document}